\documentclass{article}

\usepackage{PRIMEarxiv}

\usepackage[italicComments=true,indLines=true,noEnd=false]{algpseudocodex}
\usepackage{algorithm}
\usepackage{blindtext}
\usepackage[english]{babel}
\usepackage[utf8]{inputenc} % allow utf-8 input
\usepackage[T1]{fontenc}    % use 8-bit T1 fonts
\usepackage{hyperref}       % hyperlinks
\usepackage{url}            % simple URL typesetting
\usepackage{booktabs}       % professional-quality tables
\usepackage{amsfonts}       % blackboard math symbols
\usepackage{nicefrac}       % compact symbols for 1/2, etc.
\usepackage{microtype}      % microtypography
\usepackage{lipsum}
\usepackage{amssymb}
\usepackage{adjustbox}% graphics
\usepackage{xcolor}
\usepackage{amsmath}
\usepackage{subfigure}
\usepackage{bm}
\usepackage{caption}
\usepackage{fancyhdr}       % header
\usepackage{graphicx}       % graphics
\graphicspath{{media/}}     % organize your images and other figures under media/ folder
\newtheorem{theorem}{Theorem}
\newtheorem{lemma}[theorem]{Lemma}
\title{Constructing coarse-scale bifurcation diagrams from spatio-temporal observations of microscopic simulations: A parsimonious machine learning approach
}

\author{
  Galaris Evangelos \\
  Dipartimento di Matematica e Applicazioni ``Renato Caccioppoli''\\
  Universit\'a degli Studi di Napoli Federico II, Naples, Italy\\
  \And
  Fabiani Gianluca \\
  Scuola Superiore Meridionale\\ Universit\'a degli Studi di Napoli Federico II, Naples, Italy \\
  \And
  Gallos Ioannis \\
  School of Applied Mathematical and Physical Sciences \\
  National Technical University of Athens, Greece\\
  \And
  Kevrekidis Ioannis \\
  Department of Chemical and Biomolecular Engineering,\\
  Department of Applied Mathematics and Statistics, \\
  Department of Medicine\\
  Johns Hopkins University \\
  Baltimore, Maryland, USA\\
  \And
  Siettos Constantinos\thanks{corresponding author: constantinos.siettos@unina.it} \\
  Dipartimento di Matematica e Applicazioni ''Renato Caccioppoli'',\\
   Scuola Superiore Meridionale\\
  Universit\'a degli Studi di Napoli Federico II, Naples, Italy\\
  \\
}

\begin{document}

\maketitle

\begin{abstract}
We address a three-tier data-driven approach to solve the inverse problem in complex systems modelling from spatio-temporal data produced by microscopic simulators using machine learning. In the first step, we exploit manifold learning and in particular parsimonious Diffusion Maps using leave-one-out cross-validation (LOOCV) to both identify the intrinsic dimension of the manifold where the emergent dynamics evolve and for feature selection over the parametric space. In the second step, based on the selected features, we learn the right-hand-side of the effective partial differential equations (PDEs) using two machine learning schemes, namely shallow Feedforward Neural Networks (FNNs) with two hidden layers and single-layer Random Projection Networks(RPNNs) which basis functions are constructed using an appropriate random sampling approach. Finally, based on the learned black-box PDE model, we construct the corresponding bifurcation diagram, thus exploiting the numerical bifurcation analysis toolkit. For our illustrations, we implemented the proposed method to construct the one-parameter bifurcation diagram of the 1D FitzHugh-Nagumo PDEs from data generated by $D1Q3$ Lattice Boltzmann simulations. The proposed method was quite effective in terms of numerical accuracy regarding the construction of the coarse-scale bifurcation diagram. Furthermore, the proposed RPNN scheme was $\sim$ 20 to 30 times less costly regarding the training phase than the traditional shallow FNNs, thus arising as a promising alternative to deep learning for solving the inverse problem for high-dimensional PDEs.
\end{abstract}

% keywords can be removed

\section{Introduction}
The  discovery of physical laws and the solution of the inverse problem in complex systems modelling, i.e. the construction of Partial Differential Equations (PDEs) for the emergent dynamics from data and consequently the systematic analysis of their dynamics with established numerical analysis techniques is a holy grail in the study of complex systems and has been the focus of intense research efforts over the the last years \cite{karniadakis2021physics,schmidt2009distilling,SifanWang,kovachki2021neural}. From the early '90s, exploiting both theoretical and technological advances, researchers employed machine learning algorithms for system identification using macroscopic observations, i.e. assuming that we already know the set of coarse variables to model the underlying dynamics and the derivation of normal forms (\cite{hudson1990, Krischer1993,Rico-Martinez1994, Anderson1996, Gonzalez-Garcia1998,alexandridis2002modelling}). More recently, Bongard and Lipson \cite{bongard2007automated} proposed a method for generating symbolic equations for nonlinear dynamical systems that can be described by ordinary differential equations (ODEs) from time series. Brunton et al. ~\cite{brunton2016discovering} addressed the so-called sparse identification of nonlinear dynamics (SINDy) method to obtain {\em explicit} data-driven PDEs when the variables are known, and construct normal forms for bifurcation analysis.
Wang et al. \cite{SifanWang} addressed a physics-informed machine learning scheme based on deep learning to learn the solution operator of arbitrary PDEs. Kovachki et al. \cite{kovachki2021neural} addressed the concept of Neural Operators, mesh-free, infinite dimensional operators with neural networks, to learn surrogate functional maps for the solution operators of PDEs.\\
However, for complex systems, such ``good'' macroscopic observables that can be used effectively for modelling the dynamics of the emergent patterns are not always directly available. Thus, such an appropriate set of ``hidden'' macroscopic variables have to be identified from data. Such data can be available either directly from experiments or from detailed simulations using for example molecular dynamics, agent-based models, and Monte-Carlo methods. Hence, all in all, we confront with two major problems: (a) the identification of the appropriate variables that define (parametrize) the emerging (coarse-gained) dynamics, (b) the construction of models based on these variables. In the early 2000's, the Equation-Free and Variable-Free multiscale framework \cite{Kevrekidis2003,Kevrekidis2004,Makeev2002,Siettos2003,erban2007variable} provided a systematic framework for the numerical analysis (numerical bifurcation analysis, design of controllers, optimization, rare-events analysis) of the emergent dynamics as well as for the acceleration of microscopic simulations, by bridging the microscale where the physical laws may be known and the macroscopic scale where the emergent dynamics evolve. This bridging is achieved  via the concept of the ``coarse time steppers'', i.e. the construction of a black-box map on the macroscopic scale. By doing so, one can perform multiscale numerical analysis, even for microscopically large-scale systems tasks by exploiting the algorithms (toolkit) of matrix-free methods in the Krylov subspace \cite{Kevrekidis2003,vandekerckhove2009efficient,samaey2008newton,samaey2010analysis,Siettos2012,erban2007variable}, thus bypassing the need to construct explicitly models in the form of PDEs. In the case when the macroscopic variables are not known a-priori, one can resort to non-linear manifold learning algorithms such as Diffusion maps ~\cite{Coifman2005,Coifman2006,nadler2006diffusion,singer2009detecting} to identify the intrinsic dimension of the slow manifold where the emergent dynamics evolve.\\
Over the last few years, efforts have been focused on developing physics-informed machine learning methods for solving both the forward and inverse problems, i.e. the numerical solution of high-dimensional multiscale problems described by PDEs, and that of discovering the hidden physics \cite{karniadakis2021physics,lee2017resilient,raissi2017inferring,raissi2017machine,Lee2020}, thus both identifying the set to coarse observables and based on them to learn the effective PDEs. Lee et al. \cite{Lee2020} addressed a methodology to find the right-hand-side of macroscopic PDEs directly from microscopic data using Diffusion maps and Automatic Relevance Determination for selecting a good set of macroscopic variables, and Gaussian processes and artificial neural networks for modelling purposes. The approach was applied to learn a ``black-box'' PDE from data produced by Lattice Boltzmann simulations of the FitzHugh-Nagumo model at a specific value of the bifurcation parameter where sustained oscillations are observed.\\
%Reduced scaled models were constructed for large-scale simulations (\cite{Noid2013}). Most of the times, the useful information of these systems refers to some macro-scale features and if one is able to clarify them properly, studying the systems becomes now computationally effective. However, the proper definition of the macroscopic equations is often a hard and time-demanding issue. 
%Takunaga et al (\cite{Tokunaga1994}) proposed a method of constructing a BD only from data coming from time-series (varying the bifurcation parameter) using a neural network and reconstructed the BDs of the Henon map and coupled logistic-delayed logistic maps without observational noise as well as the BD of the R{\"o}ssler equation with observations noise (\cite{Tokuda1996}). Later, similar methods were used for reconstructing the BDs of the cubic map and FitzHugh-Nagumo equations \cite{Bagarinao2000}, R{\"o}ssler equations (\cite{Small2001}) and more. More recently, Itoh et al. (\cite{Itoh2020}) used an extreme learning machine as a time series predictor in order to reconstruct the BDs from time-series data generated by electronic circuits. 
In this paper, building on previous efforts \cite{Lee2020}, we exploit machine learning to perform numerical bifurcation analysis from spatio-temporal data produced by microscopic simulators. For the discovery of the appropriate set of coarse-gained variables, we used parsimonious Diffusion maps \cite{Dsilva2018,Holiday2019}, while for the identification of the right-hand side of the emergent coarse-grained PDEs, we used shallow Feedforward Neural Networks (FNNs) and Random Projection Neural Networks (RPNNs), thus proposing an appropriate sampling approach for the construction of the (random) basis functions. For our illustrations, we have used a $D1Q3$ Lattice Boltzmann simulator of the FitzHugh Nagumo (FHN) spatio-temporal dynamics. Upon training, the tracing of the coarse-grained bifurcation diagram was obtained by coupling the machine learning models with the pseuso-arc-length continuation approach. The performance of the machine learning schemes was compared with the reference bifurcation diagram obtained by finite differences of the FHN PDEs. 

\section{Methodology}
The pipeline of our computational framework for constructing the bifurcation diagrams from data produced from detailed microscopic simulations consists of three tasks: (a) the identification of a set of coarse-scale variables from fine-scale spatio-temporal data using manifold learning and in particular parsimonious Diffusion maps using leave-one-out cross-validation (LOOCV), (b) based on the parsimonious coarse-grained set of variables, the reconstruction of the right-hand-side of the effective PDEs using machine learning and, (c) based on the machine learning models, the construction of the coarse-scale bifurcation diagrams of the emergent dynamics using the numerical bifurcation analysis toolkit.
%For the presentation of the proposed method we assume two concentration fields (coarse-scale observables): one of the activator $u$ and one of the inhibitor $v$. Applying the LBM model, one can compute these corse-scale variables by averaging the particle distribution functions ($f_i$) on a given lattice point (see section 3.2). Acquiring the values of the coarse-scale observables, their time and space derivatives can be easily estimated. We estimate with finite differences the first order time derivatives and up to second order spatial derivatives ($u_t,v_t,u_x,v_x,u_{xx},v_{xx}$). We collect observations for different values of the parameter $\varepsilon$ (i.e. the bifurcation parameter). 

The assumption here is that the emergent dynamics of the complex system under study on a domain $\Omega\times[t_0,t_{end}] \subseteq \mathbb{R}^d\times \mathbb{R}$ can be modelled by a system, of say  $m$ (parabolic) PDEs in the form of:
\begin{equation}\label{PDEs}
\begin{aligned}
    \frac{\partial u^{(i)}(\bm{x},t)}{\partial t}\equiv u_t^{(i)}=
     F^{(i)}(t,\bm{x},\bm{u}(\bm{x},t),\mathcal{D}\bm{u}(\bm{x},t),\mathcal{D}^{\bm{2}}\bm{u}(\bm{x},t),\dots,\mathcal{D}^{\bm{\nu}}\bm{u}(\bm{x},t),\bm{\varepsilon}),\\
     (\bm{x},t)\in \Omega\times[t_0,t_{end}], \qquad i=1,2,\dots,m\\
    \end{aligned}
\end{equation}
where $\bm{u}(\bm{x},t)=[u^{(1)}(\bm{x},t),\dots,u^{(m)}(\bm{x},t)]$, $F^{(i)},$ $i=1,2,\dots m$ is a non-linear operator, $\mathcal{D}^{\bm{\nu}}\bm{u}(\bm{x},t)$ is the generic multi-index $\bm{\nu}$-th order spatial derivative at time $t$ i.e.:
\begin{equation*}
    \mathcal{D}^{\bm{\nu}}\bm{u}(\bm{x},t):=\left \{ \frac{\partial^{|\bm{\nu}|}\bm{u}(\bm{x},t)}{\partial x_1^{\nu_1}\cdots\partial x_d^{\nu_d}} \bigg| |\bm{\nu}|=\nu_1+\nu_2+\dots+\nu_d,\, \nu_1,\dots,\nu_d\ge0 \right \},
\end{equation*}
and $\bm{\varepsilon}$ denotes the (bifurcation) parameters of the system.\\
The boundary conditions read:
\begin{equation}
    B_l^{(i)}(u^{(i)}(\bm{x},t))=h^{(i)}_l(\bm{x},t)\qquad \bm{x}\in \partial\Omega_l,
\end{equation}
where $\{\partial\Omega_l\}$ denotes an $l$ partition of the boundary of $\Omega$, and initial conditions
\begin{equation}
    u^{(i)}(\bm{x},t_0)=u^{(i)}_0, \quad \bm{x}\in\Omega.
\end{equation}\\
The right-hand-side of the $i$-th PDE depend on say $\gamma^{(i)}$ number of variables and bifurcation parameters from the set of variables

\[\mathcal{S}^{(i)}=\{\bm{x},\bm{u}(\bm{x},t),\mathcal{D}\bm{u}(\bm{x},t),\mathcal{D}^{\bm{2}}\bm{u}(\bm{x},t),\dots,\mathcal{D}^{\bm{\nu}}\bm{u}(\bm{x},t),\bm{\varepsilon}\}.\] 

\begin{figure}
    \centering
    \includegraphics[width=0.9 \textwidth]{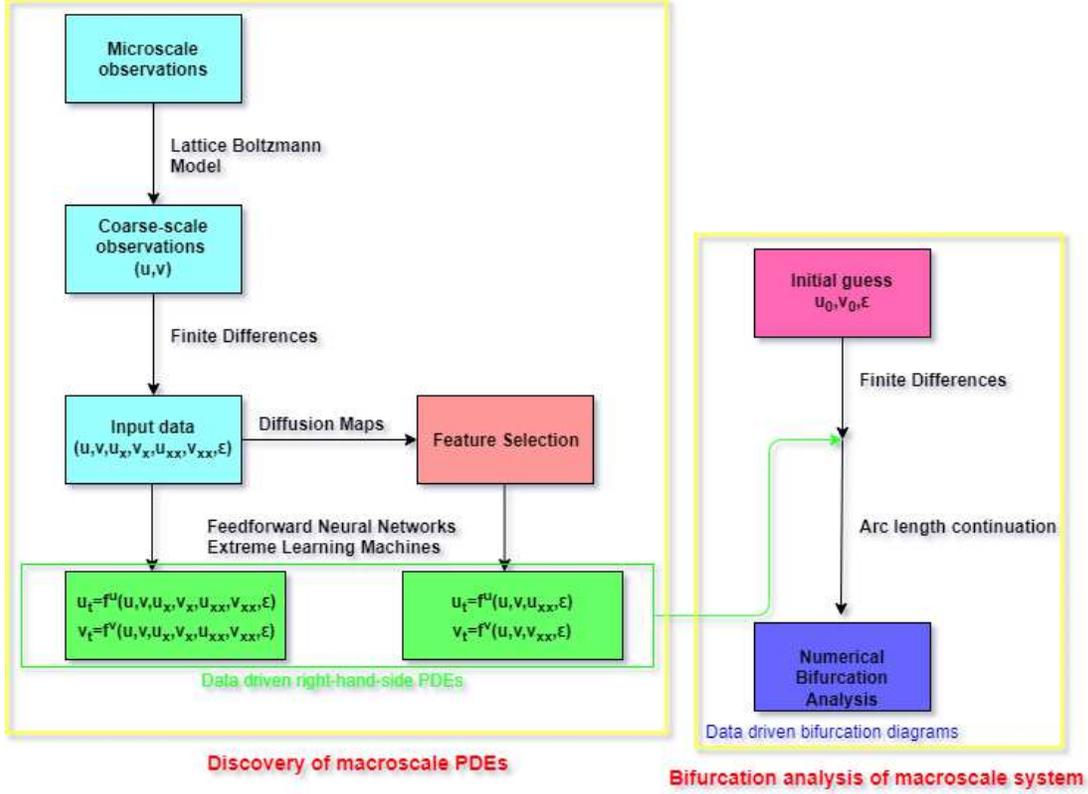}
    \caption{Schematic of the three-stage workflow for constructing coarse-grained bifurcation diagrams from fine scale observations using the paradigm of two parabolic PDEs: 1) Identify a set of parsimonious coarse-grained observables using Diffusion Maps from microscopic simulations (here $D1Q3$ Lattice Boltzmann simulations) and compute their spatial and time derivatives using finite differences, 2) ``learn'' the right-hand-side of the PDEs using machine learning (here shallow FNNs with 2 hidden layers and single-layer RPNNs), and, 3) employ the numerical bifurcation analysis toolkit (here the pseudo-arc-length-continuation method) to construct the coarse-grained bifurcation diagrams.}
    \label{overview}
\end{figure}

Let us denote this set as $\mathcal{S}^{(i)}$, with cardinality $|\mathcal{S}^{(i)}|=\gamma(i)$.
Hence, at each spatial point $\bm{x}_q,q=1,2,\dots,M$ and time instant $t_s,s=1,2,\dots,N$ the set of features for the $i$-th PDE can be described by a vector $\boldsymbol{z}_q(t_s)\in \mathbb{R}^{\gamma(i)}$.\\
Here, we assume that such macroscopic PDEs in principle exist but there are not available in a closed-form.
\\
Instead, we assume that we have detailed observations from microscopic simulations from which we can compute the time and spatial derivatives of all the observables in $N$ points in time and $M$ points in space using e.g. finite differences. Thus, we aim to (a) identify the intrinsic dimension of the manifold on which the coarse-grained dynamics evolve, i.e. for each PDE identify $\gamma(i)$, and the coordinates that define the low-dimensional manifold, i.e. the sets $\mathcal{S}^{(i)}$, and based on them (b) identify the right-hand-side (RHS) of the effective PDEs using machine learning.\\
To demonstrate the proposed approach, we have chosen to produce data from $D1Q3$ Lattice Boltzmann (LB) simulations of the coupled FitzHugh-Nagumo PDEs of activation-inhibition dynamics. Using the LB simulator, we produced data in time and space from different initial conditions and values of the bifurcation parameter. For the identification of an appropriate set of coarse-scale variables that define the low-dimensional manifold on which the emergent dynamics evolve, we performed feature selection using parsimonious Diffusion Maps \cite{Dsilva2018,Holiday2019}. Then, we trained the machine learning schemes to learn the right-hand-side of the coarse-grained PDEs on the low-dimensional manifold. Based on the constructed models, we performed numerical bifurcation analysis, employing the pseudo-arc-length continuation method. The performance of the proposed data-driven scheme for constructing the coarse-grained bifurcation diagram was validated against the one computed with the PDEs using finite differences. A schematic overview of the proposed framework for the case of two effective PDEs (as in the problem of the FitzHugh-Nagumo activation-inhibition dynamics) is shown in Figure \ref{overview}.\\
In what follows, we first describe the parsimonious Diffusion Maps algorithm for feature selection. Then, we present the machine learning schemes used for identifying the right-hand-side of the effective PDEs from the microscopic simulations, and then we show how one can couple the machine learning models with the pseudo-arc-length continuation method to construct the coarse-scale bifurcation diagrams. Finally, we present the numerical results and compare the performance of the proposed machine learning schemes. 
\subsection{Parsimonious Diffusion Maps}
\label{sec:dmaps}
Diffusion Maps is a non-linear manifold learning algorithm \cite{Coifman2005,Coifman2006,nadler2006diffusion} that identifies a low-dimensional representation  $\boldsymbol{y}_i\in\mathbb{R}^{\mu}$ of a point $\boldsymbol{z}_i\in \mathbb{R}^{n}, i=1,2,\dots N$ in the high dimensional space  ($\mu<<n$) addressing the diffusion distance among points as the preserved metric \cite{nadler2006diffusion}.
Difussion Maps assume that the high-dimensional data are embedded in a smooth low-dimensional manifold. It can be shown that the eigenvectors of the large normalized kernel matrices constructed from the data converge to the eigenfunctions of the Laplace-Beltrami operator on this manifold at the limit of infinite data \cite{Coifman2005,nadler2006diffusion}. The approximation of this Laplace-Beltrami operator is made by representing the weighted edges connecting nodes $i$ and $j$ commonly by a normalized diffusion kernel matrix $W$ with elements:
\begin{equation}
    w_{ij}=\text{exp}\biggl(-\frac{||\boldsymbol{z}_i-\boldsymbol{z}_j||^2}{\sigma}\biggr).
\end{equation}
Then, one can define the $N\times N$ diffusion matrix $P$ by:
\begin{equation}
    P=D^{-1}W, D=\text{diag} \left( \sum_{j=1}^{N}w_{ij} \right)
\end{equation}
which elements $ p_{ij}$ correspond to the probability of jumping from one point to another in the high-dimensional space.

Taking the power of $t$ of the diffusion matrix $P$ is essentially identical of observing $t$ steps forward of a Markov chain process $Z_t$ on the data points. Thus, the transition probability of moving from point $i$ to point $j$ reads:
 \begin{equation}
     p_{ij}=p(\boldsymbol{z}_i,\boldsymbol{z}_j)=Prob(Z_{t+1}=\boldsymbol{z}_j|Z_t=\boldsymbol{z}_i).
 \end{equation}
On the weighted graph, the random walk can be defined as:
 \begin{equation}
     p_{ij}=p(\boldsymbol{z}_i,\boldsymbol{z}_j)=\frac{w_{ij}}{deg(\boldsymbol{z}_i)},
 \end{equation}
where $deg(z_i)$ denotes the weighted degree of the point $i$ defined as:
 \begin{equation}
     deg(\boldsymbol{z}_i)=\sum_{j}w_{ij}.
 \end{equation}
At the next step, it is easy to compute the graph Laplacian matrix $\Tilde{P}$ defined as:
\begin{equation}
    \Tilde{P}=\Tilde{D}^{1/2}PD^{-1/2}.
\end{equation}
The eigendecomposition of $\tilde{P}$ results to $\tilde{P}=U\Lambda U^{*}$, where $\Lambda$ is a diagonal matrix with the eigenvalues and $U$ is the matrix with columns the eigenvectors of $\tilde{P}$. The eigenvalues of $P$ are the same of $\tilde{P}$ since $P$ is the adjoined of the symmetric matrix $P$, while the left and right eigenvectors of $P$ (say $\Phi$ and $\Psi$) are related to those of $\tilde{P}$ as \cite{Nadler2008Diffusion}:
\begin{equation}
    \begin{aligned}
        \Phi=UD^{1/2},\qquad
        \Psi=UD^{-1/2}.
    \end{aligned}
\end{equation}
The embedding of the manifold in $\mu$ dimensions is realized by taking the first $\mu$ non-trivial/dependent eigenvectors of $\Tilde{P}$:
\begin{equation}
    \boldsymbol{y_i}=(\lambda^t_1\phi_{1,i},\dots,\lambda^t_{\mu}\phi_{\mu,i}),\qquad i=1,\dots,N,
    \label{embedding}
\end{equation}
where $t$ denotes the number of diffusion steps (usually $t=1$) and $\lambda_1,\dots,\lambda_{\mu}$ the descending order eigenvalues. For any pair of points $\boldsymbol{z}_i$ and $\boldsymbol{z}_j$, the diffusion distance at the time step $t$ is defined as:
\begin{equation}
D^2_{t}(\boldsymbol{z}_{i},\boldsymbol{z}_{j})=\sum_{k}{\frac{(p_{t}(\boldsymbol{z}_{i},\boldsymbol{z}_{k})-p_t(\boldsymbol{z}_{j},\boldsymbol{z}_{k}))^2}{\Phi_0(\boldsymbol{z}_k)}},
\end{equation}
where $\Phi_0$ denotes the stationary distribution of the random walk described by the diffusion matrix $P$ \cite{Gao2021}:
\begin{equation}
    \Phi_0(y_i)=\frac{deg(\boldsymbol{z}_i)}{\sum_{\boldsymbol{z}_j\in Y}deg(\boldsymbol{z}_j)}.
\end{equation}
In practice, the embedding dimension $\mu$ is determined by the spectral gap in the eigenvalues of the final decomposition. Such a numerical gap means that the first few eigenvalues would be adequate for the approximation of the diffusion distance between all pairs of points \cite{Coifman2005,Coifman2006}. Here we retain only the $\mu$ parsimonious eigendimensions of the final decomposition as proposed in \cite{Dsilva2018,Holiday2019}.

\subsubsection{Feature selection using Diffusion Maps with leave-one-out cross-validation (LOOCV)}
\label{sec:feature-selection}
Here, by identifying the coarse-scale spatio-temporal behaviour of a system of PDEs, we mean learning their right-hand-sides as a black-box model. Hence, we first have to deal with the task of discovering a set of coarse-grained variables embedded in the high-dimensional input data space. For this task, one can employ various methods for feature selection such as LASSO \cite{Santosa1986, Robert1996} and Random Forests \cite{Tin1995, Tin1998}. In our framework, we used a method that extracts the dominant features based on manifold parametrization through output-informed Diffusion Maps \cite{Dsilva2018}. The core assumption of this method is that given a dataset in a high-dimensional space, then we can parametrize it on a lower-dimensional manifold.\\
%Here, we identify the intrinsic dimension of manifold on which the dynamics evolve. Hence, keeping the leading Diffusion Maps coordinates, we can test the various combination of subsets of the input domain coordinates (spatial derivatives) as to their ability to effectively parametrize the lowest-dimensional embedding of the PDE manifold.
For this purpose, given a set of $\boldsymbol{\phi}_1,\boldsymbol{\phi}_2,\dots,\boldsymbol{\phi}_{k-1} \in \mathbb{R}^N$ Diffusion Maps eigenvectors, for each element $i=1,2\dots,N$ of $\boldsymbol{\phi}_k$, we use a local linear regression model:
 \begin{equation}
     \phi_{k,i}\approx \alpha_{k,i}+\beta_{k,i}^T\Phi_{k-1,i},\quad i=1,2,\dots,N
 \end{equation}
to investigate if $\phi_k$ is an dependent eigendirection; $\Phi_{k-1,i}=[\phi_{1,i},\phi_{2,i},\dots,\phi_{k-1,i}^T]$, $\alpha_{k,i}\in \mathbb{R}$ and $\beta_{k,i}\in \mathbb{R}^{k-1}$. The values of parameters $\alpha_{k,i}$ and $\beta_{k,i}$ are found solving an optimization problem of the form:
 \begin{equation}
     \hat{\alpha}_{k,i},\hat{\beta}_{k,i}=\underset{\alpha,\beta}{\operatorname{argmin}} \sum_{j\neq i}K(\Phi_{k-1,i},\Phi_{k-1,j})(\phi_{k,j}-(\alpha+\beta^T\Phi_{k-1,j}))^2,
 \end{equation}
where $K$ is a kernel weighted function, usually the Gaussian kernel:
 \begin{equation}
     K(\Phi_{k-1,i},\Phi_{k-1,j})=\text{exp}\biggl(-\frac{||\Phi_{k-1,i}-\Phi_{k-1,j}||}{\sigma^2}\biggr),
 \end{equation}
where $\sigma$ is the shape parameter. The final normalized leave-one-out cross-validation (LOOCV) error for this local linear fit is defined as:
 \begin{equation}
     r_k=\sqrt{\frac{\sum_{i=1}^N(\phi_{k,i}-(\hat{\alpha}_{k,i}+\hat{\beta}_{k,i}^T\Phi_{k-1,i})^2}{\sum_{i=1}^\mu(\phi_{k,i})^2}}.
 \end{equation}
For small values of $r_k$, $\boldsymbol{\phi}_k$ is considered to be dependent of the other eigenvectors and hence as a harmonic or repeated eigendirection, while large values of $r_k$, suggest that $\boldsymbol{\phi}_k$ can serve as a new independent eigendirection.\\

In our approach, we provide as inputs to the Diffusion Maps algorithm the combined input-output domain (the observables $\boldsymbol{z}_i$ and their spatial and time derivatives). In principle, any of the subsets that is capable to parametrize the discovered embedding coordinates that were chosen after the above described methodology, can be considered as a new possible input data domain that can be used for learning the right-hand-side of the effective PDEs. Actually, we find the subsets of variables of the input space that minimally parametrize the intrinsic embedding by quantifying it with a total regression loss $L_T$ based on a mean squared error:
 \begin{equation}
    L_T=(\sum^\mu_{k=1}L^2_{\boldsymbol{\phi}_k})^{\frac{1}{2}}.
    \label{totalreglosserror}
\end{equation}
Here, as $L_{\phi_j}$, we define the regression loss for representing the intrinsic coordinate $\phi_j$ when using $s$ out of $n$ selected input features:
\begin{equation}
    L_{\boldsymbol{\phi}_k}=\frac{1}{N}\sum^N_{i=1}(\phi_{k,i}-g(\cdot))^2,
    \label{reglosserror}
\end{equation}
where $g(\cdot)$ is the output of the regressors with inputs the values of the features in the ambient space and target values the eigenvectors $\boldsymbol{\phi}_k$. Note, that in this procedure, we did not include the values of the bifurcation parameter into the dataset. We have chosen to employ the above method separately for every subset of the same value of the bifurcation parameter and finally to select the subset of features with the minimum sum of the total regression loses across all the embedding spaces.

\subsection{Shallow Feedforward Neural Networks}
%FNNs constitute a biologically-inspired powerful machine learning tool with a huge variety of applications. Their capability of performing complicated tasks made them very rapidly one of the most popular machine learning algorithms for regression, classification, forecasting and many more tasks. They usually contain one input, one output and one or more hidden layers. Every layer has several units (neurons), often fully connected by weights ($\omega$) and biases (b), and an activation function ($\psi(\cdot)$). We provide as inputs the macroscopic observables and their spatial derivatives while the corresponding time derivative is obtained by the output layer.
%Many results are available concerning the approximation properties of FFNs. 
%The most important one from the numerical point of view is the Universal Approximation Theorem, for which we refer to the original papers \cite{Cybenko1989,Hornik1989,Hornik1990,park1991universal,leshno1993multilayer}.
It is well known that FNNs are universal approximators of any (piecewise) continuous (multivariate) function, to any desired accuracy \cite{Cybenko1989,Hornik1989,Hornik1990,park1991universal,leshno1993multilayer}. This implies that any failure of a network must arise from an inadequate choice/calibration of weights and biases or an insufficient number of hidden nodes.\par
%The celebrated universal approximation theorem (\cite{Cybenko1989,leshno1993multilayer}) guarantees that for a single hidden layer FNN with a sufficient number of neurons an approximation function $\hat{f}$ can be found for every function $f \in C(\mathbb{R})$. Later on, Kidger et. al. (\cite{Kidger2019}) showed that the same results holds for any FNN with $K$ input units, $D$ output units, arbitrary hidden layers each one consisting of at least $D+K+2$ neurons with any continuously differentiable at at-least one point and non-zero derivative at that point activation function $\psi(\cdot)$.
The output of a FNN with two hidden layers, with $H$ hidden units in each layer, that models the right-hand-side of the $i$-th PDE at each input point $\boldsymbol{z}^{(i)}(\boldsymbol{x}_q,t_s) \in \mathbb{R}^{\gamma(i)}$, evaluated at each point in space $\boldsymbol{x}_q,\, q=1,\dots,M$, and time $t_s, s=1,2,\dots,N$ can be written as:
\begin{equation}
\begin{aligned}
    \hat{u}^{(i)}_t(\boldsymbol{x}_q,t_s)&=\hat{F}^{(i)}(\boldsymbol{z}^{(i)}(\boldsymbol{x}_q,t_s);\mathcal{W}^{1(i)},\mathcal{W}^{2(i)},\boldsymbol{b}^{1(i)},\boldsymbol{b}^{2(i)},\boldsymbol{\omega}^{o(i)},b^{o(i)})\\
    &=\sum_{j=1}^H\omega_j^{o(i)}\psi(\sum^H_{l=1}[\omega_{lj}\psi(\boldsymbol{\omega}_l^{(i)}\cdot \boldsymbol{z}^{(i)}(\boldsymbol{x}_q,t_s)+b^{1(i)}_l)+b^{2(i)}_j)]+b^{o(i)}.
\end{aligned}
\label{eq:FNN}
\end{equation}
$\psi(.)$ is the activation function (based on the above formulation it is assumed to be the same for all nodes in the hidden layers) $\boldsymbol{\omega}^{o(i)}=(\omega_1^{o(i)},\omega_2^{o(i)},\dots,...\omega_H^{o(i)})\in \mathbb{R}^{1 \times H}$ are the external weights connecting the second hidden layer and the output, $b^{o(i)} \in \mathbb{R}$ is the bias of the output node, the matrix $\mathcal{W}^{1(i)} \in \mathbb{R}^{H\times \gamma(i)}$ with rows $\boldsymbol{\omega}_l^{1(i)} \in \mathbb{R}^{\gamma(i)}$ are the weights connecting the input and the first hidden layer, $\boldsymbol{b}^{1(i)}=(b^{1(i)}_1,b^{1(i)}_2,\dots,b^{1(i)}_H) \in \mathbb{R}^H$ are the biases of the nodes of the first hidden layer, the matrix $\mathcal{W}^{2(i)}\in\mathbb{R}^{H\times H}$ contains the weights $\omega_{lj}$ connecting the $l$-th unit of the first hidden layer with the $j$-th unit of the second hidden layer and $\boldsymbol{b}^{2(i)}=(b^{2(i)}_1,b^{2(i)}_2,\dots,b^{2(i)}_H)\in \mathbb{R}^H$ are the biases of the second hidden layer. In the same way, one can easily extend the above formula for more than two hidden layers.
Then, a loss function for each one of the $m$ PDEs can be specified as:
\begin{equation}
    E^{(i)}=\sum_{q=1}^M\sum_{s=1}^N(u_{t}^{(i)}(\boldsymbol{x}_q,t_s)-\hat{u}_{t}^{(i)}(\boldsymbol{x}_q,t_s))^2,
\end{equation}
The main task of a neural network is the generalization. Foresee and Hagan (\cite{Dan1997}) showed that adding the regularization term $E_{\omega}=\sum_{j=1}^H\omega^2_j$ to the cost function will maximize the posterior probability based on Bayes' rule. Hence, the total cost function is:
\begin{equation}
    E_{total}=E+\lambda E_w,
\end{equation}

where $\lambda$ is the regularization parameter that has to be tuned. For our simulations we used the Bayesian regularized back-propagation updating the weight values using the Levenbgerg-Marquadt algorithm.  (\cite{Hagan1994}) as implemented in Matlab.

\subsection{Random Projection Neural Networks}
Random Projection Neural Networks (RPNNs) are a family of neural networks including Random Vector Functional Links (RVFLs) \cite{barron1993universal,igelnik1995stochastic}, Reservoir Computing/ Echo state networks \cite{verstraeten2007experimental,jaeger2001echo}, Extreme Learning Machines \cite{huang2006extreme} and Liquid-State Networks \cite{maass2002real}. The basic idea, which seed can be found in the pioneering work of Rosenblatt back in '60s \cite{rosenblatt1962perceptions}, behind all these approaches is to use FNNs with fixed-weights  between the input and the hidden layer(s), fixed biases for the nodes of the hidden layer(s), and a linear output layer. Based on that configuration, the output is projected linearly onto the functional subspace spanned by the nonlinear basis functions of the hidden layer, and the only remaining unknowns  are the weights between the hidden and the output layer. Their estimation is done by solving a nonlinear regularized least squares problem \cite{fabiani2021numerical,calabro2021extreme}. The universal approximation properties of the RPNNs has been proved in a series of papers (see e.g. \cite{barron1993universal,igelnik1995stochastic,jaeger2001echo,huang2006extreme}).
In general, the universal approximation property of random projections can been can be rationalized by the celebrated Johnson and Lindenstrauss (JL) Theorem \cite{johnson1984extensions}:
\begin{theorem}[Johnson and Lindenstrauss]
Let $\mathcal{Z} \in \mathbb{R}^{n\times N}$
matrix with $N$ points $\boldsymbol{z}_i \in \mathbb{R}^n$. Then, $\forall \epsilon \in (0,1)$ and $\mu \in \mathbb{N}$ such that $\mu \geq O (\frac{\ln{N}}{\epsilon^2})$, there exists a map $G: \mathbb{R}^n \rightarrow \mathbb{R}^{\mu}$ such that $\forall \boldsymbol{z}_i, \boldsymbol{z}_j \in \mathcal{Z}$:
\begin{equation} \label{eqn:jl}
    (1-\epsilon)\lVert \boldsymbol{z}_i-\boldsymbol{z}_j \rVert ^2  \leq \lVert G(\boldsymbol{z}_i)-G(\boldsymbol{z}_j) \rVert ^2 \leq (1+\epsilon) \lVert \boldsymbol{z}_i-\boldsymbol{z}_j \rVert ^2.
\end{equation}
\end{theorem} %\qed \\
Note, that while the above theorem is deterministic, its proof relies on probabilistic techniques combined with Kirszbraun's theorem to yield a so-called extension mapping~\cite{johnson1984extensions}.  In particular, it can be shown, that one of the many such embedding maps is simply a linear projection matrix $R$ with entries $r_{ij}$ that are i.i.d. random variables sampled from a normal distribution. In particular, the JL Theorem may be proved using the following lemma.
\begin{lemma}
Let $\mathcal{Z}$ be a set of $N$ points in $\mathbb{R}^n$ and let $G(\boldsymbol{z})$ be the random projection defined by
\begin{equation*}
    G (\boldsymbol{z}) = \frac{1}{\sqrt{\mu}} R \, \boldsymbol{z}, \quad \boldsymbol{z} \in \mathbb{R}^n,
\end{equation*}
where $R = [r_{ij}] \in \mathbb{R}^{\mu \times n}$ has components that are i.i.d.~random variables sampled from a normal distribution.
Then, $\forall \, \boldsymbol{z} \in \mathcal{Z}$
\begin{equation*}
(1-\epsilon) \lVert \boldsymbol{z} \rVert ^2  \leq \lVert G (\boldsymbol{z}) \rVert ^2 \leq (1+\epsilon) \lVert \boldsymbol{z}\rVert ^2 
\end{equation*}
is true with probability $p \geq 1-2\exp \left(-(\epsilon^2-\epsilon^3)\frac{\mu}{4} \right)$.
\end{lemma}
Similar proofs have been given for distributions different from the normal one (see, e.g. \cite{achlioptas2003database,dasgupta2003elementary,vempala2005random,wang2012geometric}).

The above is a feature mapping, which  may result in a dimensionality reduction ($\mu<n$) or, in analogy to the case of kernel-based manifold learning methods, a projection into a higher dimensional space ($\mu>n$). We also note that while the above linear random projection is but one of the choices for constructing a JL embedding, it has been experimentally demonstrated and/or theoretically proven that appropriately constructed nonlinear random embeddings may outperform simple linear random projections. For example, in \cite{giryes2016deep} it was shown that deep networks with random weights for each layer result in even better approximation accuracy than the simple linear random projection.

Here, for learning the right-hand-side of the set of PDEs, we used RPNNs (in the form of an Extreme Learning Machine\cite{huang2006extreme}) with a single hidden layer \cite{fabiani2021numerical}. For $H$ hidden units in the hidden layer, the output of the proposed RPNN can be written as:
\begin{equation}
\begin{aligned}
    \hat{u}^{(i)}_t(\boldsymbol{x}_q,t_s)&=\hat{F}^{(i)}(\boldsymbol{z}^{(i)}(\boldsymbol{x}_q,t_s);\mathcal{W}^{(i)},\boldsymbol{b}^{(i)},\boldsymbol{\omega}^{o(i)},b^{o(i)})=\\
    &=\sum_{j=1}^H\omega_j^{o(i)}\psi(\boldsymbol{\omega}_j^{(i)}\cdot \boldsymbol{z}^{(i)}(\boldsymbol{x}_q,t_s)+b^{(i)}_j)+b^{o(i)},
\end{aligned}
    \label{eq:RPNN}
\end{equation}
where $\boldsymbol{z}^{(i)}(\boldsymbol{x}_q,t_s) \in \mathbb{R}^{\gamma(i)}$ denotes the inputs computed at each point in space $\boldsymbol{x}_q,\, q=1,\dots,M$, and time $t_s, s=1,2,\dots,N$, $\psi(.)$ is the activation function $\boldsymbol{\omega}^{o(i)}=(\omega_1^{o(i)},\omega_2^{o(i)},\dots,...\omega_H^{o(i)})\in \mathbb{R}^{1 \times H}$ are the external weights connecting the hidden layer and the output and $b^{o(i)} \in \mathbb{R}$ is the bias of the output node, while the matrix $\mathcal{W}^{(i)} \in \mathbb{R}^{H\times \gamma(i)}$ with rows $\boldsymbol{\omega}_j^{(i)} \in \mathbb{R}^{\gamma^{(i)}}$ and $\boldsymbol{b}^{(i)}=(b^{(i)}_1,b^{(i)}_2,\dots,b^{(i)}_H) \in \mathbb{R}^H$ are the weights connecting the input and the hidden layer and the biases of the hidden layer, respectively. Now, in the proposed RPNN scheme, $\bm{\omega}^{(i)}_j$ and $b^{(i)}_j$ are random variables drawn from appropriate uniform distributions, the output bias $b^o$ is set to zero, and the output weights $\bm{\omega}^{o(i)} \in \mathbb{R}^{1\times H}$ are determined solving a linear least squares problem:
\begin{equation}
    \hat{\bm{u}}^{(i)}_t=\boldsymbol{\omega}^{o(i)}\mathcal{A}^{(i)}, \qquad i=1,\dots,m,
    \label{LSP}
\end{equation}
where $\hat{\bm{u}}^{(i)}_t\in \mathbb{R}^{MN}$ is the vector collecting all the outputs $\hat{u}^{(i)}_t(\bm{x}_q,t_s)$ of the RPNN for $q=1,\dots,M$ and $s=1,\dots,N$,  and the matrix $\mathcal{A}^{(i)} \in \mathbb{R}^{H\times MN}$ is the matrix which elements $\mathcal{A}_{j,k}$ are given by:
\begin{equation}
    \mathcal{A}_{j,k}^{(i)}=\psi(\boldsymbol{\omega}^{(i)}_j \cdot \boldsymbol{z}^{(i)}_k + b^{(i)}_j),
\end{equation}
where $\bm{z}^{(i)}_k=\bm{z}^{(i)}(\bm{x}_q,t_s)$ and $k=q+(s-1)M$.
Regarding the regression problem, generally $H<<MN$, the system Eq. (\ref{LSP}) is over-determined. As the resulting projection matrix $\mathcal{A}$ is not guaranteed to be full row-rank the solution can be computed with Singular Value Decomposition (SVD). Given the SVD decomposition of $\mathcal{A}$, the pseudo inverse $\mathcal{A}^{+}$ is:
\begin{equation}
    \mathcal{A}^{(i)}=U\Sigma V^T,  \qquad (\mathcal{A}^{(i)})^{+}=V \Sigma^{+} U^T,
\end{equation}
where $U \in \mathbb{R}^{H\times H}$ and $V \in \mathbb{R}^{MN \times MN}$ are the unitary matrices of left and right eigenvectors
respectively, and $\Sigma \in \mathbb{R}^{H\times MN}$ is the diagonal matrix of $H$ singular values $\sigma_j$. Finally, in order to regularize the problem, we can select just $\tilde{H}<H$ singular values $\tilde{\sigma}$ that are greater than a given tolerance, i.e., $\tilde{\sigma} \in \{\sigma_j\, | \,\sigma_j>tol, j=1,\dots,H\}$. Hence, the output weights $\boldsymbol{\omega}^{o(i)}$ are computed as:
\begin{equation}
    \boldsymbol{\omega}^{o(i)}=\hat{\bm{u}}^{(i)}_t\tilde{V} \tilde{\Sigma}^{+} \tilde{U}^T.
\end{equation}
where $\tilde{U} \in \mathbb{R}^{H\times\tilde{H}}$, $\tilde{V} \in \mathbb{R}^{MN\times\tilde{H}}$ and $\tilde{\Sigma} \in \mathbb{R}^{\tilde{H}\times\tilde{H}}$ are restricted to the $\tilde{\sigma}$s.

For the regression problem, we aim at learning the right-hand-side of the PDEs from spatio-temporal data with single-layer RPNNs with $H$ random basis functions:
\begin{equation}
    \psi_j^{(i)}(\boldsymbol{z}^{(i)})=\psi(\boldsymbol{\omega}_j^{(i)} \cdot \boldsymbol{z}^{(i)} + b^{(i)}_j).
\end{equation}
Then the approximated function $\hat{F}^{(i)}$ is just a linear combination of the random basis functions $\psi_j^{(i)}$. For our computations, we selected as activation function the logistic sigmoid $\psi:y\in\mathbb{R}\rightarrow \psi(y) \in\mathbb{R}$ given by:
\begin{equation}
    \psi(y)= \frac{1}{1+\text{exp}(-y)},
\end{equation}
where, $y$ is given by linear combination $y=\boldsymbol{\omega}_j^{(i)} \cdot \boldsymbol{z}^{(i)}+b^{(i)}_j$.

\subsubsection{Random Sampling Procedure}
For the construction of the appropriate set of random basis functions for the solution of the inverse problem (i.e. that of learning the effective PDEs from data), we suggest a different random sampling procedure, than the one usually implemented in RPNNs and in particular in Extreme Learning Machines \cite{fabiani2021numerical,calabro2021extreme,schiassi2021extreme,dwivedi2020physics,dong2021local,dong2022local} for the solution of the forward problem, i.e. that of the numerical solution of Partial Differential Equations. Since in the inverse problem, we aim at solving a high-dimensional over-determined system ($MN>>H$) is important to parsimoniously select the underlying basis functions $\psi_j^{(i)}$, i.e. to seek for appropriate internal weights $W^{(i)}$ and biases $\bm{b}^{(i)}$ that lead to non-trivial functions.  

In general, the weights $\omega_j^{(i)}$ and biases $b_j^{(i)}$ are uniformly random sampled from a subset of the input/feature space, e.g., $\omega_j^{(i)},\, b_j^{(i)} \sim \mathcal{U}([-1,1]^{\gamma(i)}$, where an high dimension $\gamma(i)$ of the input/feature space leads to the phenomenon of curse of dimensionality. Indeed, it is necessary to use many function ($H \propto 10^{\gamma(i)}$) to correctly ``explore'' the input space and give a good basis function.

Hence, our goal is to construct $\omega_j^{(i)}$ and $b_j^{(i)}$ with a simple data-driven manifold learning in order to have a basis of functions $\psi_j^{(i)}$ that well describe the manifold $\mathcal{M}^{(i)}$ where the data  $\bm{z}^{(i)}(\bm{x}_q,t_s) \in \mathcal{M}^{(i)}, \forall q,\forall s$ are embedded. It is well-known that the output of a neuron is given by a ridge function $f:\mathbb{R}^H \rightarrow \mathbb{R}$ such that $f(z_1,\dots,z_n)=g(\boldsymbol{a}^T\cdot \boldsymbol{z})$, where $g:\mathbb{R} \rightarrow \mathbb{R}$ and $\boldsymbol{a} \in \mathbb{R}^n $. The inflection point of the logistic sigmoid is at ($y=0$, $\psi(y)=1/2$). The points that satisfy the following relation \cite{fabiani2021numerical,calabro2021extreme}:
\begin{equation}
    y= \boldsymbol{\omega}^{(i)}_j \cdot \boldsymbol{z}^{(i)}(\bm{x}_q,t_s)+b_j^{(i)}=0
    \label{eq:dir_hyperplane}
\end{equation}
form an hyperplane $\mathcal{H}_j^{(i)}$ of $\mathbb{R}^{MN}$ (MN dimension of $\boldsymbol{z}$) defined by the direction of $\boldsymbol{\omega}_j^{(i)}$. Along $\mathcal{H}_j$, $\psi_j^{(i)}$ is constantly $1/2$. We call the points $c^{(i)}_j \in \mathcal{H}_j^{(i)}$ the \emph{centers} of the ridge function $\psi_j^{(i)}$.
Here the goal is to select $H$ centers $c^{(i)}_j$ that are on the manifold $\mathcal{M}^{(i)}$ (note that this is not achieved by random weights and biases) and find directions $w_j^{(i)}$ that make $\psi_j^{(i)}$ non-constant/non-trivial along $\mathcal{M}^{(i)}$ (note that for ridge functions there are many directions for which this does not happen).

Thus, here, being $H<<MN$ we suggest to uniformly random sample $H$ points $c^{(i)}_j$ from $\boldsymbol{z}(\bm{x}_q,t_s)$
to be the centers of the functions $\psi_j^{(i)}$: in this way the inflection points of $\psi_j^{(i)}$ are on the manifold $\mathcal{M}$.
Also, we independently randomly sample other $H$ points $\tilde{c}_j^{(i)}$ from the inputs $\boldsymbol{z}(\bm{x}_q,t_s)$. Then, we construct the hidden weights as:
\begin{equation}
    \boldsymbol{\omega}^{(i)}_j=\tilde{c}^{(i)}_j-c^{(i)}_j,
\end{equation}
in order to  set the direction $\boldsymbol{\omega}^{(i)}_j$ of the hyperplane $H_j^{(i)}$ parallel to the one connecting $\tilde{c}^{(i)}_j$ and $c^{(i)}_j$. By doing so, the ridge function will be constant on a direction orthogonal to the connection between two points in the manifold $\mathcal{M}^{(i)}$ and along this line will change in value, so it will be able to discriminate between the points lying on this direction. Thus, the biases $b_j^{(i)}$ are set as:
\begin{equation}
    b^{(i)}_j=-\boldsymbol{\omega}^{(i)}_j \cdot c^{(i)}_j.
\end{equation}
Eq. \eqref{eq:dir_hyperplane} ensures that $c^{(i)}_j \in \mathcal{H}_j^{(i)}$ is a center of the ridge function.

\section{Coarse-grained numerical bifurcation analysis from spatio-temporal data}
For assessing the performance of the proposed scheme, we selected the celebrated, well studied FitzHugh-Nagumo (FHN) model first introduced in \cite{fitzhugh1961impulses} to simplify the Hodgkin-Huxley model into a two-dimensional system of ODEs to describe the dynamics of the voltage  across a nerve cell. In particular, we consider the FHN equations which add a spatial diffusion term to describe the propagation of an action potential as a traveling wave. The bifurcation diagram of the one-dimensional set of PDEs is known to have a turning point and two supercritical Andronov-Hopf bifurcation points. In what follows, we describe the model along with the initial and boundary conditions and then we present the $D1Q3$ Lattice Boltzmann model.

\subsection{The Macroscale model: the FitzHugh-Nagumo Partial Differential Equations}
The evolution of activation $u:[x_0,x_{end}]\times[t_0,t_{end}]\rightarrow \mathbb{R}$ and inhibition $v:[x_0,x_{end}]\times[t_0,t_{end}]\rightarrow \mathbb{R}$ dynamics are described by the following two coupled nonlinear parabolic PDEs:
\begin{equation}
\begin{aligned}
&\frac{\partial u(x,t)}{\partial t}=D^{u}\frac{\partial^{2}u(x,t)}{\partial x^2}+u(x,t)-u(x,t)^3-v(x,t),\\ 
&\frac{\partial v(x,t)}{\partial t}=D^{v}\frac{\partial^{2}v(x,t)}{\partial x^2}+\varepsilon(u(x,t)-\alpha_1v(x,t)-\alpha_0),
\end{aligned}
\end{equation}
with homogeneous von Neumann Boundary conditions:
\begin{equation}
\begin{aligned}
    &\frac{du(x_{end},t)}{dx}=0, \quad \frac{dv(x_0,t)}{dx}=0.\\
    &\frac{du(x_{end},t)}{dx}=0, \quad \frac{dv(x_0,t)}{dx}=0.
\end{aligned}
\end{equation}
$\alpha_0$ and $\alpha_1$ are  parameters, $\varepsilon$ is the kinetic bifurcation parameter. 

For our simulations, we have set $x_0=0$, $x_{end}=20$, $\alpha_1=2,\alpha_0=-0.03, D^u=1, D^v=4$ and varied the bifurcation parameter $\varepsilon$ in the interval $[0.005, 0.955] \cite{Theodoropoulos2000}$.  
%used a uniform discretization of the spatial domain $[x_0,x_{end}]=[0,20]$ with a step $\Delta x=0.2$ and of the time domain on $[t_0,t_{end}]=[0,450]$ with a time step $\Delta t=0.01$.
%We integrate in time numerically using the Lattice Boltzmann Method
For our simulations, in order to explore the dynamic behaviour, we considered various initial conditions $u_0(x)=u(x,0)$ and $v_0(x)=v(x,0)$ selected randomly as follows:
\begin{equation}
\begin{aligned}
    &u_0(x)=w \tanh{\big(\alpha (x-c)\big)}+\beta\\
    &v_0(x)=0.12 \cdot u_0(x).\\
    & w \sim \mathcal{U}(0.8,1.2),
    & \alpha \sim \mathcal{U}(0.5,1)\\
    & c \sim \mathcal{U}(2,18),
    & \beta \sim \mathcal{U}(-0.4,0),
\end{aligned}
\label{eq:random_initial_cond}
\end{equation}
where $\mathcal{U}(a,b)$ denotes the uniform distribution in the interval $[a,b]$.

\subsection{The $D1Q3$ Lattice Boltzmann model}
%Lattice-Boltzman modeling (LBM) (\cite{Lee2020,Chen1998, Yeomans2002}) can be used as a mesoscopic numerical simulation approach for identifying spatio-temporal dynamics of finite-difference-type discretizations of the Boltzman-BGK equations (\cite{Bhatnagar1954}).
The Lattice Boltzmann model serves as our fine-scale simulator. The statistical description of the system at a mesoscopic level uses the concept of distribution function $f(\vec{r},\vec{c},t)$, i.e. $f(\vec{r},\vec{c},t)d\vec{r} d\vec{c} dt$ is the infinitesimal probability of having particles at location $\vec{r}$ with velocities $\vec{c}$ at a given time $t$, for reducing the high-number of equations and unknowns. Then, at this level, a system without an external force is governed by the Boltzmann Transport equation \cite{Bhatnagar1954}:
\begin{equation}
    \frac{\partial f}{\partial t }+\vec{c}\cdot \nabla f=\mathcal{R}(f),
    \label{eq:Boltzmann_transport}
\end{equation}
where the term $\mathcal{R}(f)$ describes the rate of collisions between particles. In 1954, Bhatnagar, Gross and Krook (BGK) \cite{bhathnagor1954model} introduced an approximation model for the collision operator:
\begin{equation}
    \mathcal{R}(f)=\frac{1}{\tau}(f^{eq}-f),
    \label{eq:BGK_operator}
\end{equation}
where $\tau$ is the so-called relaxing time coefficient and $f^{eq}$ denote the local equilibrium distribution function.

In the LBM, Eq.\eqref{eq:Boltzmann_transport}-\eqref{eq:BGK_operator} is collocated (assumed valid) along specific directions $\vec{c}_i$ on a lattice:
\begin{equation}
    \frac{\partial f_i}{\partial t }+\vec{c}_i\cdot \nabla f_i=\frac{1}{\tau}(f_i^{eq}-f_i)
    \label{eq:Boltzmann_direction}
\end{equation}
and then Eq.\eqref{eq:Boltzmann_direction} is discretized with a time step $\Delta t$ as follows:
\begin{equation}
    f_i(\vec{r}+\vec{c}_i\Delta t,t+\Delta t)=f_i(\vec{r},t)+\frac{\Delta t}{\tau}(f_i^{eq}-f_i).
    \label{eq:Boltzmann_discretized}
\end{equation}
One common interpretation of Eq.\eqref{eq:Boltzmann_discretized} is to think about the distribution functions as fictitious particles that stream and collide along specified linkages of the lattice. Lattices are usually denoted by the notation $DnQm$, where $n$ is the spatial dimension of the problem and $m$ refer to the number of connections of each node in the lattice. The node in the lattices coincide with the points of a spatial grid with a spatial step $\Delta x$.

Here, in order to estimate the coarse-scale observables $u$ and $v$ of the FHN dynamics, we considered the $D1Q3$ implementation, i.e. we used the one-dimensional lattice with three velocities $c_i$: particles can stream to the right ($c_1=\frac{\Delta x}{\Delta t}$), to the left ($c_{-1}=-\frac{\Delta x}{\Delta t}$) or staying still on the node ($c_0=0$). Also, we assume the coexistence of two different distribution functions for describing the distribution of the activator particles $f^u_i$ and the distribution of the inhibitor particles $f^v_i$, where the subscript $i$ refer to the associated direction. Therefore, one can figure that at each instant there are six fictitious particles on each node of the lattice: two resting on the node (with distribution $f^u_0$ and $f^v_0$), two moving on the left (with distribution $f^u_{-1}$ and $f^v_{-1}$) and two moving on the right (with distribution $f^u_{1}$ and $f^v_{1}$). The relation between the above distributions and the coarse-scale density $u$ and $v$ is given by the zeroth moment (across the velocity directions) of the overall distribution function:
\begin{equation}
\begin{aligned}
    &u(x_j,t_k)= \sum^1_{i=-1}f^u_i(x_j,t_k), \\
    &v(x_j,t_k)= \sum^1_{i=-1}f^v_i(x_j,t_k).
    \label{eq:coarse_lbm}
\end{aligned}
\end{equation}
The coexistence of multiple distributions renders necessary to introduce weights $\omega_i$ for the connections in the lattice that should satisfy the following properties:
\begin{itemize}
    \item [(a)] Normalization $\omega_0+\omega_1+\omega_{-1}=1$
    \item [(b)] Simmetry  $\omega_1-\omega_{-1}=0$
    \item [(c)] Isotropy:
    \begin{itemize}
        \item [(c.1)] $\omega_0c_0^2+\omega_1c_1^2+\omega_{-1}c_{-1}^2=c_s^2$
        \item [(c.2)]$\omega_0c_0^3+\omega_1c_1^3+\omega_{-1}c_{-1}^3=0$
        \item [(c.3)] $\omega_0c_0^4+\omega_1c_1^4+\omega_{-1}c_{-1}^4=3c_s^4$,
    \end{itemize}
\end{itemize}
where $c_s$ is the speed of sound in the lattice. Thus, the weights are equal to $\omega_{\pm 1}=1/6$ for the moving particles and $\omega_0=4/6$ for the resting particle. The resulting speed of sound in the lattice is $c_s=\frac{\sqrt{3}\Delta x}{3\Delta t}$.

As the BGK operator \eqref{eq:BGK_operator} suggests, one key step in applying LBM for solving reaction-advection-diffusion PDEs is to determine the local equilibrium distribution function $f^{eq}$ associated to a given model. For particles with macroscopic density $\rho$ that move in a medium macroscopic velocity $\vec{u}_m$, the Maxwell distribution is:
\begin{equation}
    \begin{aligned}
       f^{eq}(\vec{c})&=\frac{\rho}{(2\pi RT)^{d/2}}\text{exp}\biggl( -\frac{(\vec{c}-\vec{u}_m)^2}{2RT} \biggr)=\\
       &=\frac{\rho}{(2\pi RT)^{d/2}}\text{exp}\biggl(-\frac{\vec{c}\cdot \vec{c}}{2RT} \biggr)\text{exp}\biggl( -\frac{-2\vec{c}\cdot\vec{u}_m+\vec{u}_m\cdot\vec{u}_m}{2RT} \biggr),
    \end{aligned}
    \label{eq:eq_medium}
\end{equation}
where $d$ is the spatial dimension of the problem, $T$ is the temperature and $R$ is the universal gas constant. The exponential in Eq. \eqref{eq:eq_medium} can be expanded using Taylor series, ignoring terms of order $O(u^3)$ and higher, thus obtaining:
\begin{equation}
    f^{eq}(\vec{c})=\rho \omega(\vec{c})\biggl[ 1+ \frac{2\vec{c}\cdot \vec{u}_m-\vec{u}_m\cdot \vec{u}_m}{2 c_s^2}+\frac{(\vec{c}\cdot \vec{u}_m)^2}{2c_s^4}\biggr],
    \label{eq:eq_medium_exp_taylor}
\end{equation}
with $\omega(\vec{c})=(2\pi RT)^{-d/2}\text{exp}\biggl(-\dfrac{\vec{c}\cdot \vec{c}}{2RT} \biggr)$ and $RT=c_s^2$, with $c_s$ speed of the sound.

Now, since the FHN PDEs are only diffusive, i.e. there are no advection terms, the medium is stationary ($\vec{u}_m=0$) and the equilibrium distribution function, discretized on the lattice direction $c_i$, is simplified in:
\begin{equation}
\begin{aligned}
    &f^{u,eq}_i(x_j,t_k)= \omega_i u(x_j,t_k), \quad i=-1,0,1\\
    &f^{v,eq}_i(x_j,t_k)= \omega_i v(x_j,t_k).
    \label{eq:lbm}
\end{aligned}
\end{equation}

Now, in the FHN model, we need to consider also reaction terms $R^l_i$ and so finally, the time evolution of the microscopic simulator associated to the FHN on a given $D1Q3$ lattice is:
\begin{equation}
    f_i^l(x_{j+i},t_{k+1})=f_i^l(x_j,t_k)+\frac{\Delta t}{\tau^l}(f_i^{l,eq}(x_j,t_k)-f_i^{l}(x_j,t_k))+\Delta t R^l_i(x_j,t_k), l\in \{u,v\}
\end{equation}
where the superscript $l$ denotes the activator $u$ and the inhibitor $v$ and the reaction terms $R^l_i$ are directly derived by:
\begin{equation}
\begin{aligned}
    &R^u_i(x_j,t_k)= \omega_i(u(x_j,t_k)-u^3(x_j,t_k)-v(x_j,t_k)), \\
    &R^v_i(x_j,t_k)= \omega_i \, \varepsilon(u(x_j,t_k)-\alpha_1v(x_j,t_k)-\alpha_0).
\end{aligned}
\end{equation}

Finally, the relaxation coefficient $\dfrac{\Delta t}{\tau^l}$ is related to the macroscopic kinematic viscosity $D^l$ of the FHN model and in general depends on the speed of the sound $c_s$ associated to the lattice \cite{Qian1995}:
\begin{equation}
    \dfrac{\Delta t}{\tau^l}=\frac{2}{1+\frac{2}{c_s^2\Delta t}D^l}=\frac{2}{1+6D^l\frac{\Delta t}{\Delta x^2}}.
\end{equation}

\section{Algorithm flow chart}
Summarizing, the proposed three-tier algorithm for constructing bifurcation diagrams from data is provided in the form of a pseudo code in Algorithm \ref{alg:ml-bds}. The first two steps are related to the identification of the effective coarse scale observables and the learning of the right-hand-side of the effective PDEs. The third step is where the pseudo-arc-length continuation method is applied for the tracing of the solution branch through saddle-node bifurcations.
\begin{algorithm}[ht!]
\caption{Construct coarse bifurcation diagrams from spatio-temporal data from microscopic (here LB) simulations  \label{alg:ml-bds}}
\small
\smallskip
\begin{algorithmic}[0] %[0] for no numbering 
\footnotesize \Require Grid $N_{\varepsilon}$ of $\boldsymbol{\varepsilon}$ \Comment{\footnotesize set grid for the values of the bifurcation parameter $\varepsilon$}
\footnotesize \Require  $\boldsymbol{x}$ and $\boldsymbol{t}$ be the space and time grid
\State\textbf{1. \footnotesize Use Diffusions Maps to identify a parsimonious set of coarse scale observables from data (here produced by Lattice Boltzmann simulations)}
\State $L_t^{u}\gets$ $0$, $L_t^{v}\gets$ $0$ 
\For{$\varepsilon=1,\dots,N_{\varepsilon}$}
\State Select $w\sim \mathcal{U}(0.8,1.2)$, $c\sim \mathcal{U}(2,18)$,\\
$\alpha\sim\mathcal{U}(0.5,1)$, $\beta\sim\mathcal{U}(-0.4,0)$.
\State $u(\boldsymbol{x},0,\varepsilon)\gets wtanh(\alpha(\boldsymbol{x}-c))+\beta$
\State $v(\boldsymbol{x},0,\varepsilon)\gets0.12\cdot u_{i}(\boldsymbol{x},0,\varepsilon)$ \Comment{see eq. \ref{eq:random_initial_cond}}
\State $u(\boldsymbol{x},\boldsymbol{t},\varepsilon),v(\boldsymbol{x},\boldsymbol{t},\varepsilon)\gets LBM(u(\boldsymbol{x},0,\varepsilon),v(\boldsymbol{x},0,\varepsilon),\boldsymbol{t})$ \Comment{Lattice Boltzman simulator, see eq. \ref{eq:coarse_lbm}}
\footnotesize \State  Compute $u_{t},v_{t},u_{x},v_{x},u_{xx},v_{xx}$ \Comment{\footnotesize here, using central finite differences.}
\footnotesize  \State Compute the first $\mu$ Diffusion Maps (DM) eigenvectors:\\ 
$[\phi_1^u,\dots,\phi_{\mu}^u]\gets$ DM($u,v,u_{t},u_{x},v_{x},u_{xx},v_{xx}$)\\
$[\phi_1^v,\dots,\phi_{\mu}^v]\gets$ DM($u,v,v_t,u_x,v_x,u_{xx}^,v_{xx}$) \Comment{see eq. \ref{embedding}}
\State $\boldsymbol{z}\gets[u,v,u_x.v_x.u_{xx},v_{xx}]$
\For{every $q \subset z$}
\For{$k=1,\dots,\mu$}
\State $\hat{\phi^u_k}\gets GP(q)$\\
$\hat{\phi^v_k}\gets GP(q)$\Comment{GP: Gaussian process regressor}
\State $L_{\boldsymbol{\phi}^u_k}\gets\frac{1}{N}\sum^N_{i=1}(\phi^u_{k,i}-\hat{\phi^u_{k,i}})^2$
\State $L_{\boldsymbol{\phi}^v_k}\gets\frac{1}{N}\sum^N_{i=1}(\phi^v_{k,i}-\hat{\phi^v_{k,i}})^2$\Comment{see eq. \ref{reglosserror}}
\EndFor
\State $L_t^{u}(q)\gets L_t^{u}(q)+(\sum_{j=1}^{\mu}L^2_{\phi^u_j})^{\frac{1}{2}}$
\State $L_t^{v}(q)\gets L_t^{v}(q) +(\sum_{j=1}^{\mu}L^2_{\phi^v_j})^{\frac{1}{2}}$ \Comment{see eq. \ref{totalreglosserror} and \ref{reglosserror}}
\EndFor
\EndFor
\State $\boldsymbol{z}_u\gets\{q^*: L_t^{u}(q^*)=min(L_t^{u})\}$\\
$\boldsymbol{z}_v\gets\{q^*: L_t^{u}(q^*)=min(L_t^{v})\}$
\Comment{\footnotesize extract the effective features}\\
\textbf{2. \footnotesize Based on the extracted set of coarse variables from Step 2, learn the right-hand-sides of the coarse scale PDEs}
\State Train the FNNs/RPNNs:\\ $\hat{F}^u_r\equiv\hat{u}_{t,r}\gets$FFN/RPNN$(\boldsymbol{z}_u,\varepsilon)$\\ $\hat{F}^v_r\equiv\hat{v}_{t,r}\gets$FNN/RPNN$(\boldsymbol{z}_v,\varepsilon)$\\
\textbf{3. \footnotesize Wrap around the machine learning models the numerical bifurcation analysis toolkit (here the pseudo arc-length continuation method) to systematically study the emergent dynamics}
\end{algorithmic}
\end{algorithm}

\section{Numerical Results}
\subsection{Numerical bifurcation analysis of the FHN PDEs}
\begin{figure}
    \centering
    \subfigure[]{
    \includegraphics[width=0.47 \textwidth]{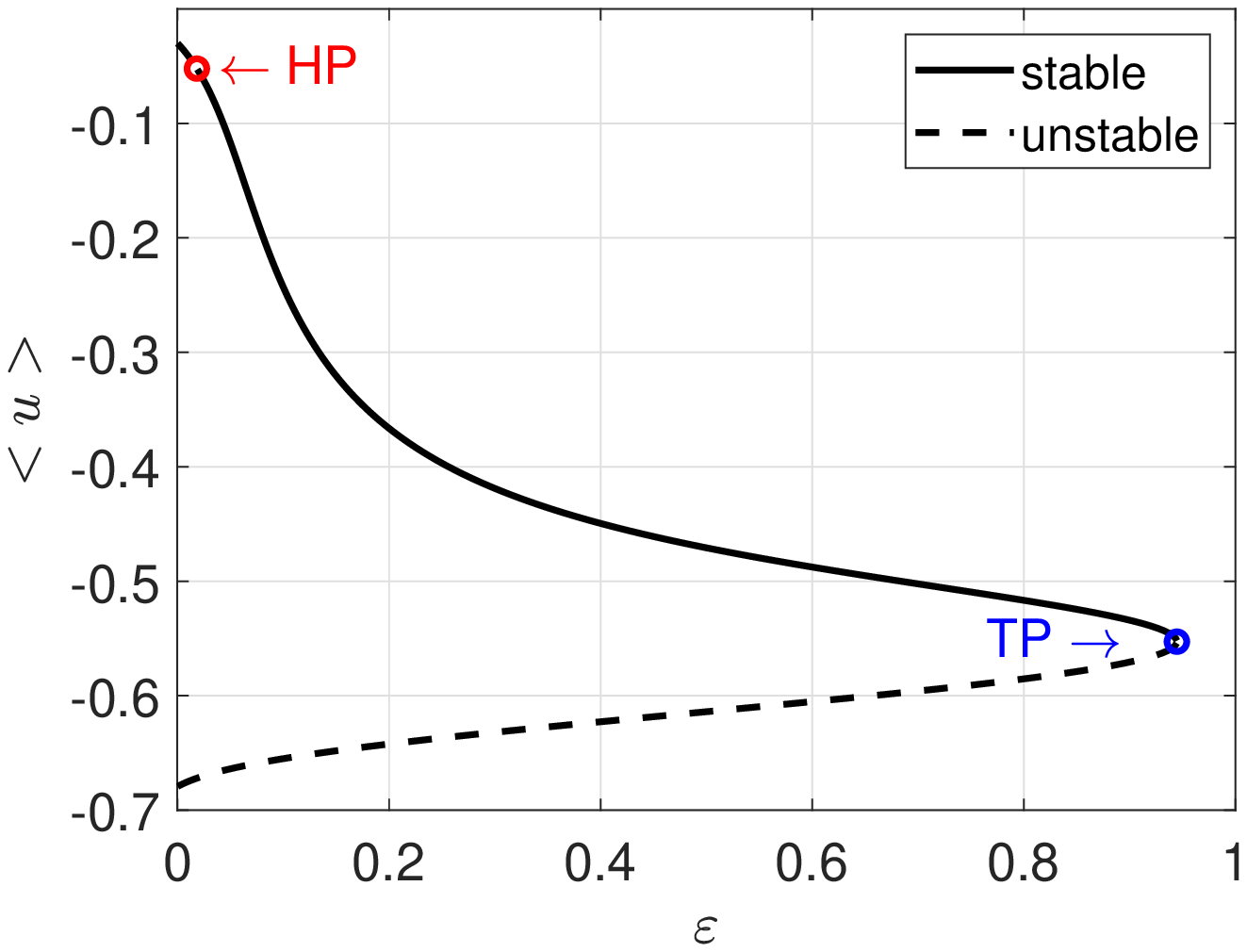}
    }
    \subfigure[]{
    \includegraphics[width=0.47 \textwidth]{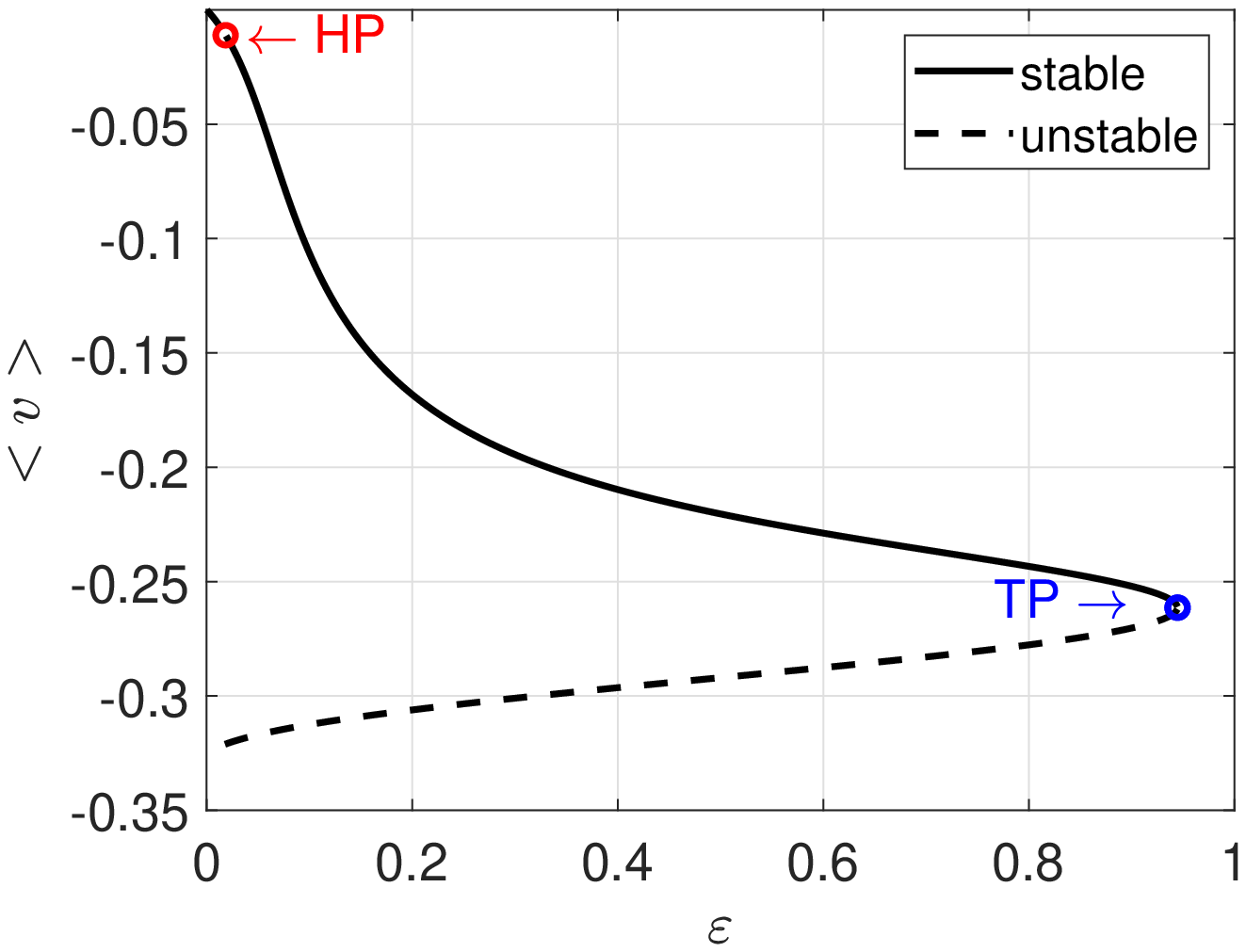}
    }
    \caption{Reference bifurcation diagram of the FHN PDEs with respect to $\varepsilon$ as computed with FD and $N=200$ points. (a) Mean values $<u>$ for stable and unstable branches, (b) Mean values $<v>$ for stable and unstable branches. Andronov-Hopf Point: $HP_{\varepsilon}$=0.01827931. Turning Point: $TP_{\varepsilon}$=0.94457768.}
    \label{fig:FD_bif_diag}
\end{figure}

For comparison purposes, we first constructed the bifurcation diagram of the FHN PDEs using central finite differences. The discretization of the one-dimensional PDEs in $M$ points with second-order central finite differences in the unit interval $0\leq x \leq 20$ leads to the following system of $2(M-2)$ non-linear algebraic equations $\forall x_j=(j-1)h, j=2,\dots M-1 $, $h=\frac{1}{M-1}$:
 \begin{equation*}
     \begin{aligned}
     F_j^u(u,v)=\frac{D^u}{h^2}(u_{j+1}-2u_{j}+u_{j-1})+u_j-u^3_j-v_j=0\\
     F_j^v(u,v)=\frac{D^v}{h^2}(v_{j+1}-2v_{j}+v_{j-1})+\varepsilon(u_j-\alpha_1v_j-\alpha_0)=0.
      \end{aligned}
 \end{equation*}
At the boundaries, we imposed homogeneous von Neumann boundary conditions.
The  above $2(M-2)$ set of non-linear algebraic equations is solved iteratively using Newton's method. The non-null elements of the Jacobian matrix are given by:
\begin{equation*}
    \begin{aligned}
    \frac{\partial F_j^u}{\partial u_{j-1}}=\frac{D^u}{h^2} ;
    \frac{\partial F_j^u}{\partial u_j}=-D^u\frac{2}{h^2}-3u_j^2 ; \frac{\partial F_j^u}{\partial u_{j+1}}=\frac{D^u}{h^2};
    \frac{\partial F_j^u}{\partial v_j}=-1\\
    \frac{\partial F_j^v}{\partial v_{j-1}}=\frac{D^v}{h^2} ; 
    \frac{\partial F_j^v}{\partial v_j}=-D^v\frac{2}{h^2}-\varepsilon\alpha_1v_j ; 
    \frac{\partial F_j^v}{\partial v_{j+1}}=\frac{D^u}{h^2};
    \frac{\partial F_j^v}{\partial u_j}=\varepsilon.
    \end{aligned}
\end{equation*}
To trace the solution branch along the critical points, we used the pseudo arc-length-continuation method (\cite{Chan1982, Glowinski1985, govaerts2000numerical}). 
This involves the parametrization of $u(x)$, $v(x)$ and $\varepsilon(x)$ by the arc-length $s$ on the solution branch. The solution is sought in terms of $\tilde{u}(x,s)$, $\tilde{v}(x,s)$ and $\tilde{\varepsilon}(s)$ in an iterative manner, by solving until convergence the following augmented system:
\begin{equation}
\begin{bmatrix}
\nabla_u \boldsymbol{F^u} &
\nabla_v \boldsymbol{F^u} &
\nabla_{\varepsilon}\boldsymbol{F^u}\\
\nabla_u \boldsymbol{F^v} &
\nabla_v \boldsymbol{F^v} &
\nabla_{\varepsilon}\boldsymbol{F^v}\\
\nabla_{u} {\boldsymbol{N}} & \nabla_v \boldsymbol{N} & \nabla_{\varepsilon} {\boldsymbol{N}}
\end{bmatrix} \begin{bmatrix}d u^{(n)}(x,s)\\d v^{(n)}(x,s)\\d\varepsilon^{(n)}(s) \end{bmatrix} =-\begin{bmatrix}\boldsymbol{F^u}(u^{(n)}(x,s),v^{(n)}(x,s),\varepsilon^{(n)}(s))\\\boldsymbol{F^v}(u^{(n)}(x,s),v^{(n)}(x,s),\varepsilon^{(n)}(s))\\\boldsymbol{N}(u^{(n)}(x,s),v^{(n)}(x,s),\varepsilon^{(n)}(s))\end{bmatrix},
\label{augmentedarclength}
\end{equation}
where 
\begin{equation*}
\begin{aligned}
\nabla_{\varepsilon }\boldsymbol{F^u}=\begin{bmatrix}
\frac{\partial F^u_1}{\partial \varepsilon} & \frac{\partial F^u_2}{\partial \varepsilon} & \dots & \frac{F^u_N}{\partial \varepsilon}
\end{bmatrix}^T, \nabla_{\varepsilon }\boldsymbol{F^v}=\begin{bmatrix}
\frac{\partial F^v_1}{\partial \varepsilon} & \frac{\partial F^v_2}{\partial \varepsilon} & \dots & \frac{F^v_M}{\partial \varepsilon}
\end{bmatrix}^T,
\end{aligned}
\end{equation*}
and
\begin{equation*}
    \begin{split}
\boldsymbol{N}(u^{(n)}(x,s),v^{(n)}(x,s),\varepsilon^{(n)}(s))= & \\ (u^{(n)}(x,s)-&\tilde{u}(x,s)_{-2})^T\cdot \frac{(\tilde{u}(x)_{-2}-\tilde{u}(x)_{-1})}{ds}+\\
(v^{(n)}(x,s)-&\tilde{v}(x,s)_{-2})^T\cdot \frac{(\tilde{v}(x)_{-2}-\tilde{v}(x)_{-1})}{ds}+\\
(\varepsilon^{(n)}(s)-\tilde{\varepsilon}_{-2})&\cdot \frac{(\tilde{\varepsilon}_{-2}-\tilde{\varepsilon}_{-1})}{ds}-ds,
    \end{split}
\end{equation*}
where ($\tilde{u}(x)_{-2}$,$\tilde{v}(x)_{-2}$) and ($\tilde{u}(x)_{-1}$,$\tilde{v}(x)_{-1}$) are two already found consequent solutions for $\tilde{\varepsilon}_{-2}$ and $\tilde{\varepsilon}_{-1}$, respectively and $ds$ is the arc-length step for which a new solution around the previous solution $(\tilde{u}(x)_{-2},\tilde{v}(x)_{-2},\tilde{\varepsilon}_{-2})$ along the arc-length of the solution branch is being sought. 
The corresponding reference bifurcation diagram is shown in Figure \ref{fig:FD_bif_diag}. In this range of values, there is an Andronov-Hopf bifurcation at $\varepsilon\approx 0.018497$ and a fold point at $\varepsilon\approx 0.95874$. 

\begin{figure}
    \centering
    \subfigure[]{
    \includegraphics[width=0.47 \textwidth]{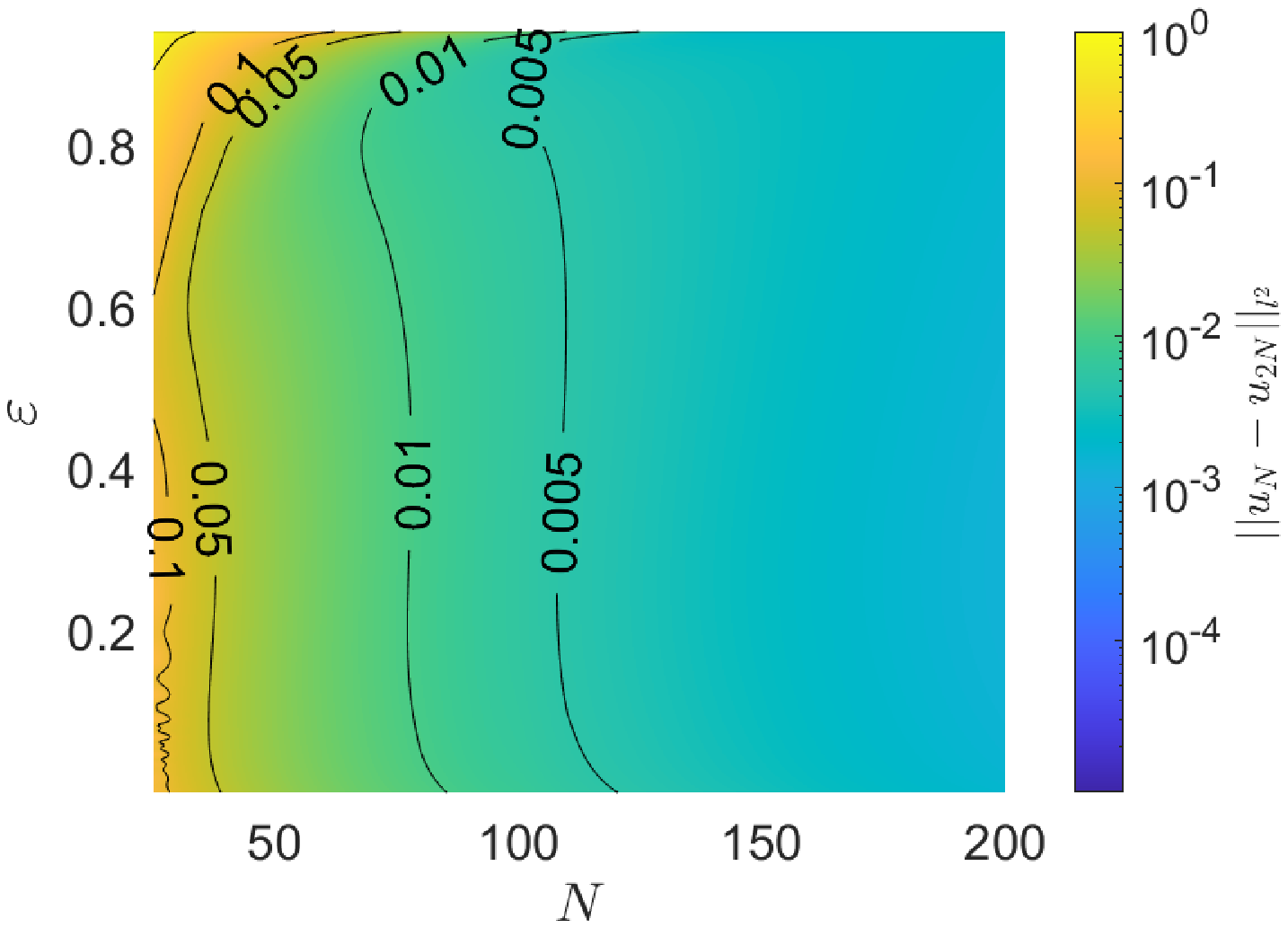}
    }
    \subfigure[]{
    \includegraphics[width=0.47 \textwidth]{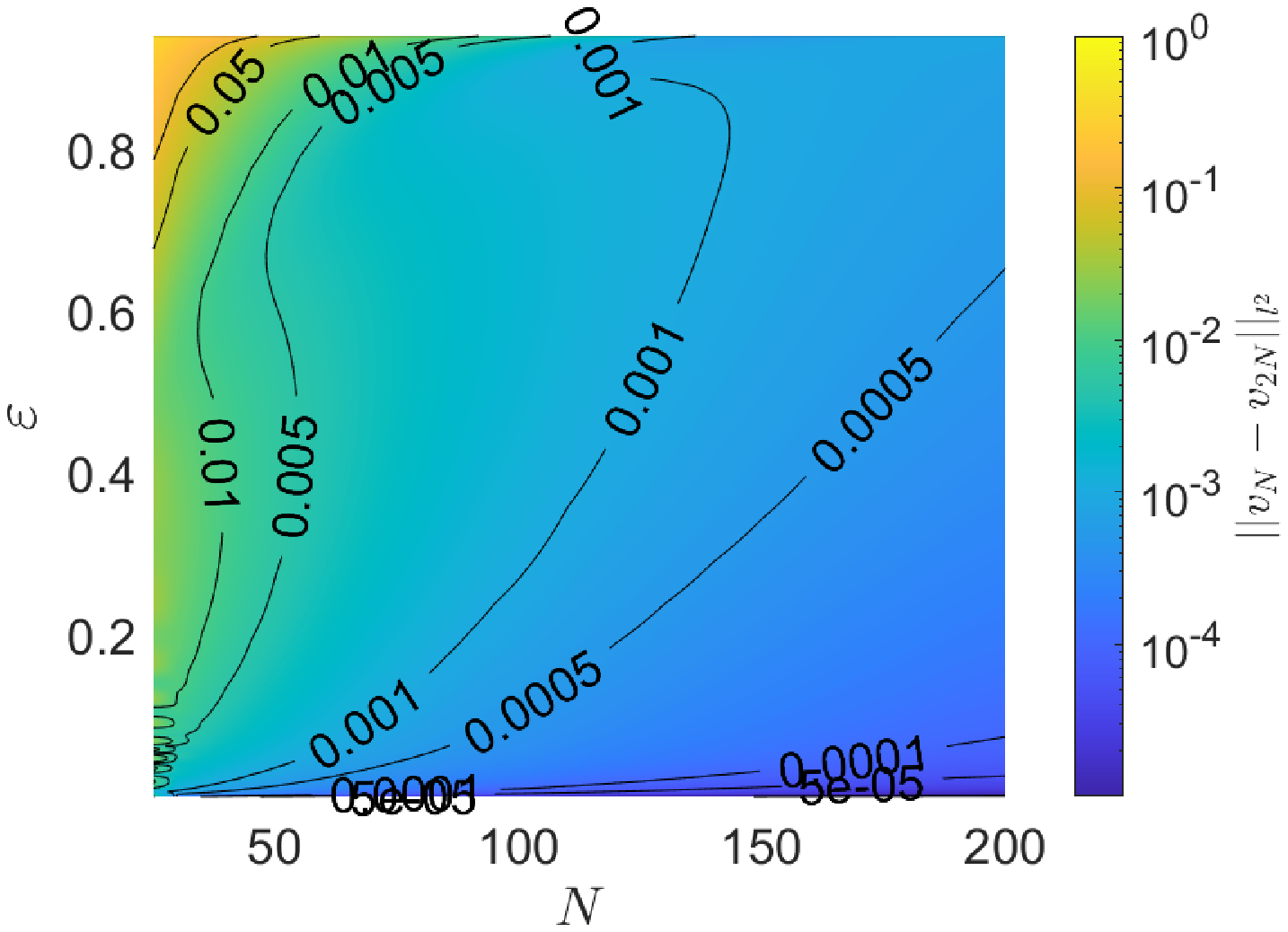}
    }
    \subfigure[]{
    \includegraphics[width=0.47 \textwidth]{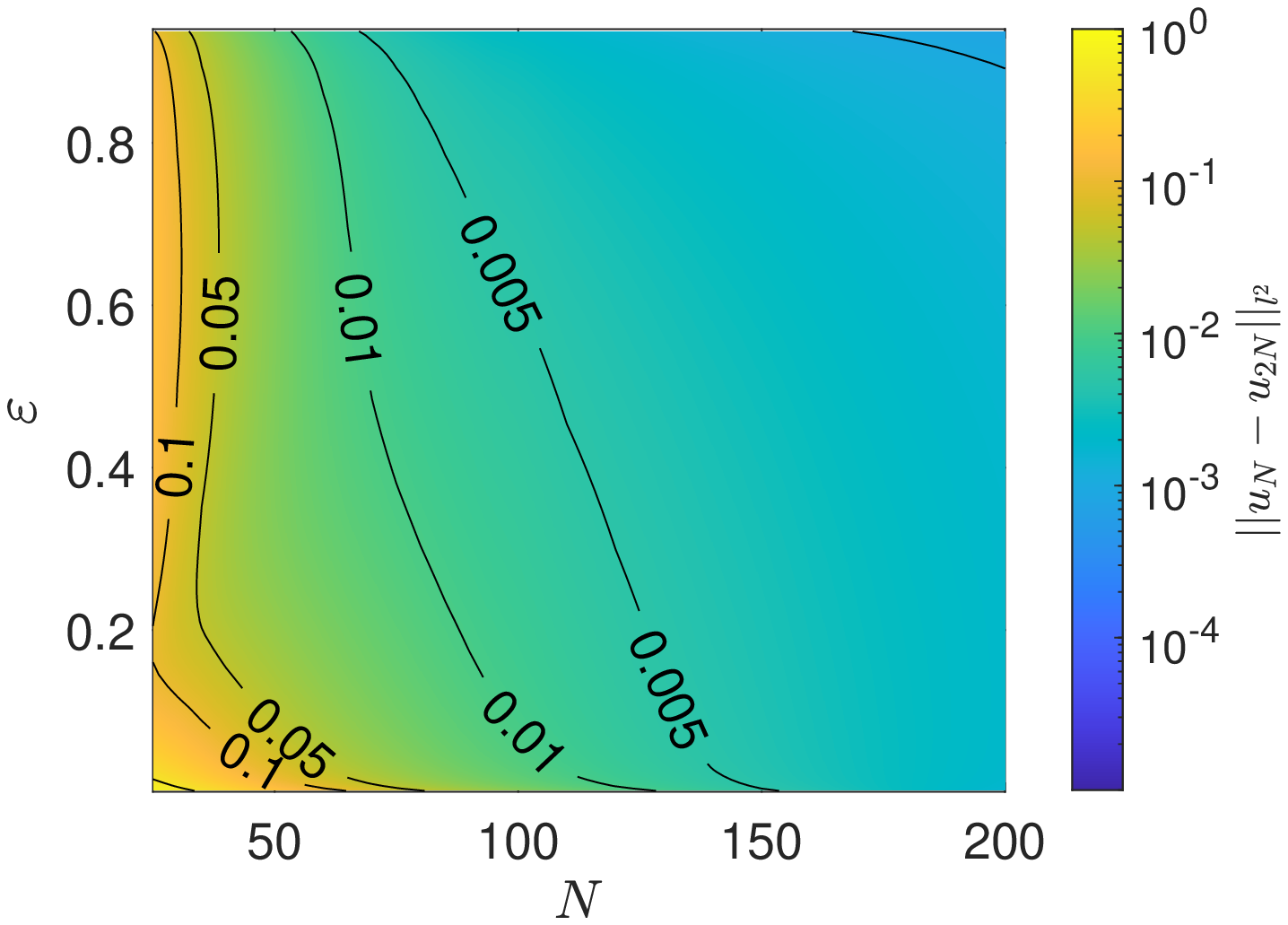}
    }
    \subfigure[]{
    \includegraphics[width=0.47 \textwidth]{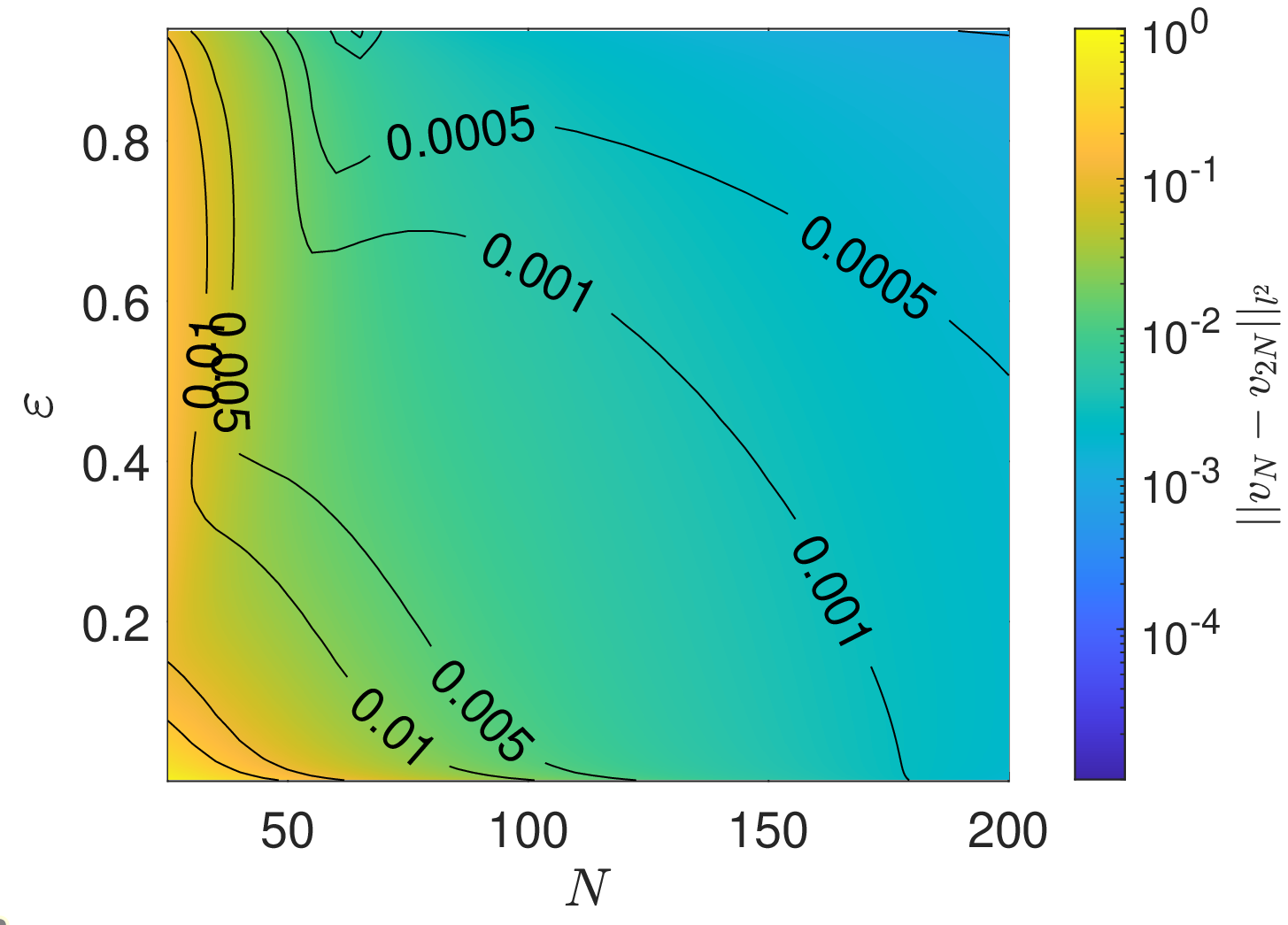}
    }
    \caption{Contour plot of the $l^2$ norms of the convergence of the solutions as computed with finite differences with respect to the size of the grid $N$ computed as $||u_{N}-u_{2N}||_{l^2}$, $||v_{N}-v_{2N}||_{l^2}$. The convergence error was computed on 1001 grid points, using linear piecewise interpolation. (a) upper branch for $u$, (b) upper branch for $v$, (c) lower branch for $u$, (d) lower branch for $v$.}
    \label{fig:CP_FD}
\end{figure}

\subsection{Numerical bifurcation analysis from microscopic simulations}
\begin{figure}[ht]
    \centering
    \subfigure[]{
    \includegraphics[width=0.47 \textwidth]{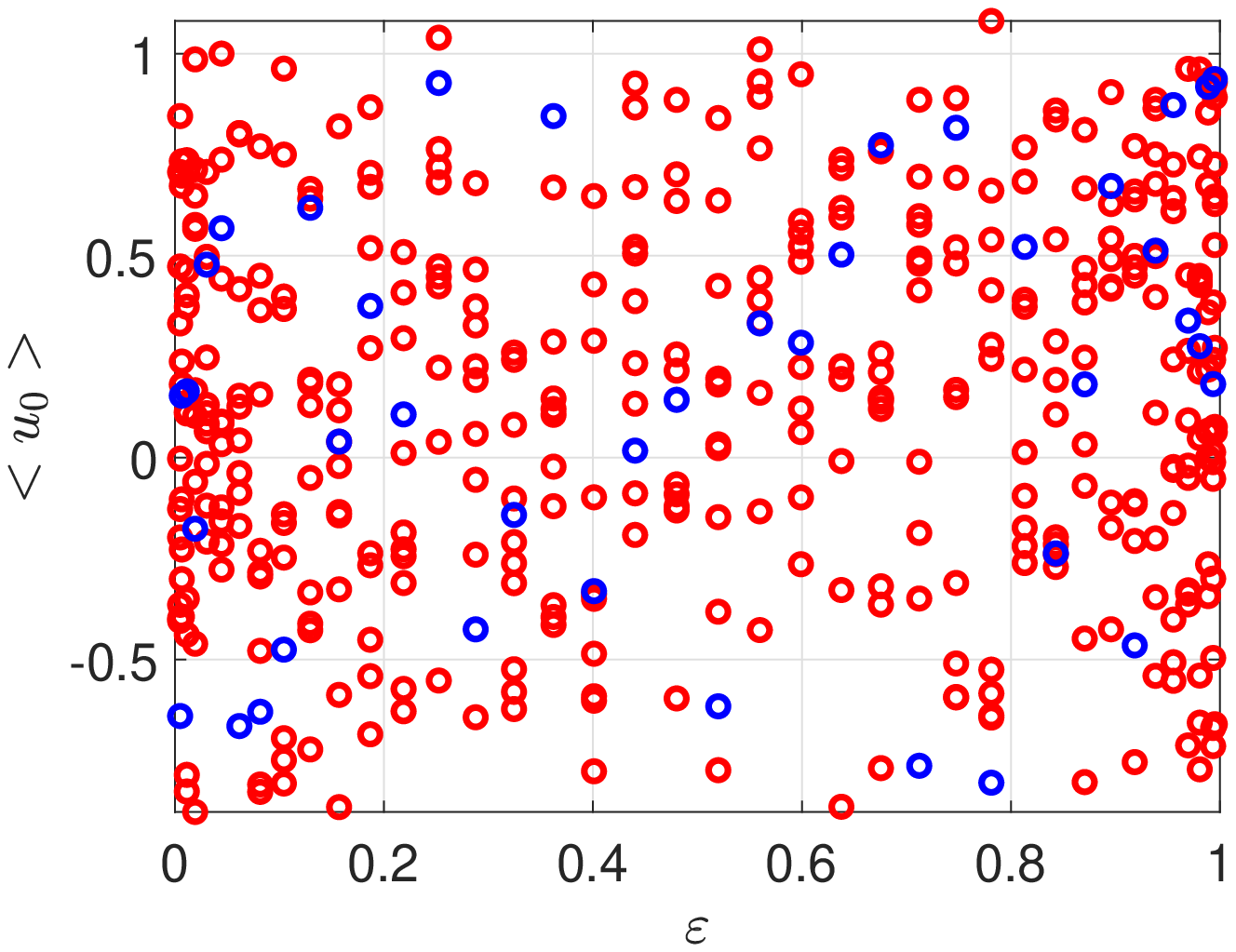}
    }
    \subfigure[]{
    \includegraphics[width=0.47 \textwidth]{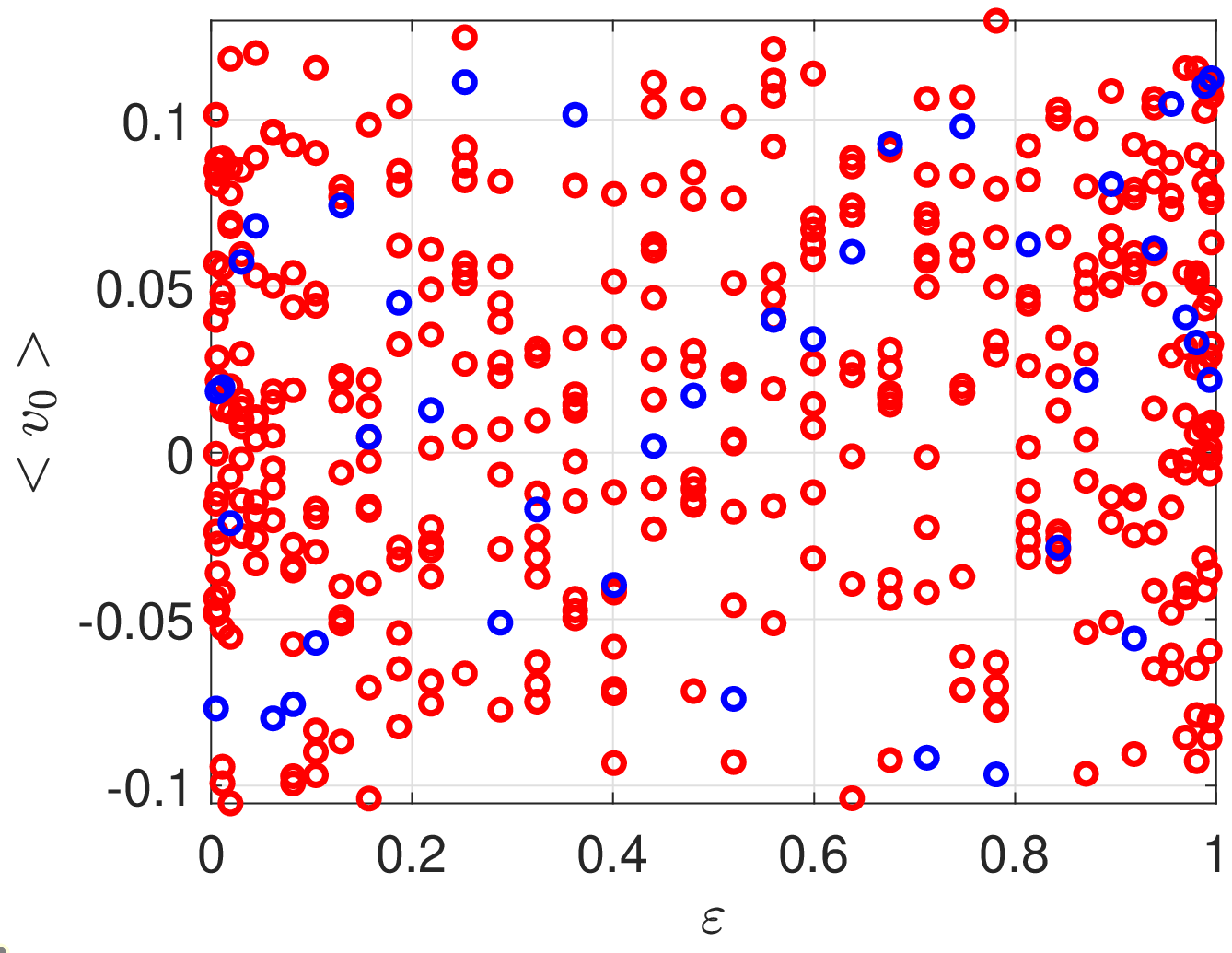}
    }
    \caption{Coarse initial conditions for (a) $u$ and (b) $v$ for the training. Every dot denotes a point whose $\varepsilon$ and mean $u$ (or $v$) were used for input data. Red dots are training points, blue points are test points. The grid is spanned with Chebychev-Gauss-Lobatto points for epsilons in the interval $[0.005,0.995]$ and the initial condition are randomly selected as in Eq. \eqref{eq:random_initial_cond}}
    \label{fig:initial_cond}
\end{figure}

We collected transients of $u(x,t)$ and $v(x,t)$ with a sampling rate of 1s, from 10 different random sampled initial conditions for 40 different values for the bifurcation parameter $\varepsilon$. In particular, we created a grid of 40 different $\varepsilon$ in [0005, 0.955] using Gauss-Chebychev-Lobatto points, while the 10 initial conditions are sampled according to Eq.(\ref{eq:random_initial_cond}).
Figure \ref{fig:initial_cond} depicts the total of 400 training initial conditions.
Thus, we end up with a dataset consisting of 40 (values of $\varepsilon$)$\times$10 (initial conditions)$\times$448 (time points ignoring the first 2s of the transient)$\times$40 (space points) $\backsimeq 7.168.000$ data points.\\
For learning the coarse-grained dynamics and construct the corresponding bifurcation diagram, we trained two FNNs and two single-layer RPNNs (one for each one of the variables $u$ and $v$). The FNNs were constructed using two hidden layers with 12 units in each layer. Hidden units were employed with the hyperbolic tangent sigmoid activation function, while the regularization parameter was tuned and set $\lambda=0.01$. 
For the training of the FNNs, we used the Deep Learning toolbox of MATLAB 2021a on an Intel Core i5-8265U with up to 3.9 GHz frequency with a memory of 8 GB.

\subsubsection{Numerical bifurcation analysis without feature selection}
Table \ref{errors_networks} summarises the performance of the two schemes on the training and on the test data sets.
\begin{table}[]
\begin{adjustbox}{width=1\textwidth}
\begin{tabular}{ccccc|cccc}
\hline
\multicolumn{1}{l}{}                                            & \multicolumn{4}{c|}{{test set}}                                                  & \multicolumn{4}{c}{{training set}}                                                  \\ \hline
\multicolumn{1}{l}{}                                            & {MSE ($u$)} & {$l^{\infty}$ ($u$)} & {MSE ($v$)} & {$l^{\infty}$ ($v$)} & {MSE ($u$)} & {$l^{\infty}$ ($u$)} & {MSE ($v$)} & {$l^{\infty}$ ($v$)} \\ \hline
{FNN}                                                    & 7.90e-09         & 2.26e-02                & 1.56e-09         & 6.63e-03                & 1.31e-09         & 7.00e-03                & 2.78e-10         & 2.58e-03                \\
{\begin{tabular}[c]{@{}c@{}} FNN(FS)\end{tabular}}  & 5.39e-08         & 2.93e-02                & 1.16e-08         & 7.65e-03                & 1.90e-08         & 2.64e-02                & 1.37e-09         & 4.70e-03                \\
{RPNN}                                                   & 2.91e-08         & 2.98e-02                & 4.50e-10         & 2.22e-03                & 6.90e-09         & 2.37e-02                & 7.06e-11         & 9.40e-04                \\
{\begin{tabular}[c]{@{}c@{}}RPNN(FS)\end{tabular}} & 7.10e-08         & 3.07e-02                & 1.73e-08         & 1.60e-02                & 2.16e-08         & 2.85e-02                & 6.30e-10         & 3.84e-03                \\ \hline
\end{tabular}
\end{adjustbox}
\caption{Mean-square error (MSE) and $l^{\infty}$ errors between the predicted $\hat{u}_t$ and $\hat{v}_t$ from the FNNs and RPNNs and the actual time derivatives $u_t$ and $v_t$ without and with feature selection (FS).}
\label{errors_networks}
\end{table}
As it is shown, for any practical purposes, both schemes resulted to equivalent numerical accuracy for all metrics. For the FNNs, the training phase (using the deep-learning toolbox in Matlab R2020b) required $\sim 1000$ epochs and $\sim 4$ hours with the minimum tolerance set to $1e-07$. 

\begin{figure}
    \centering
    \subfigure[]{
    \includegraphics[width=0.47 \textwidth]{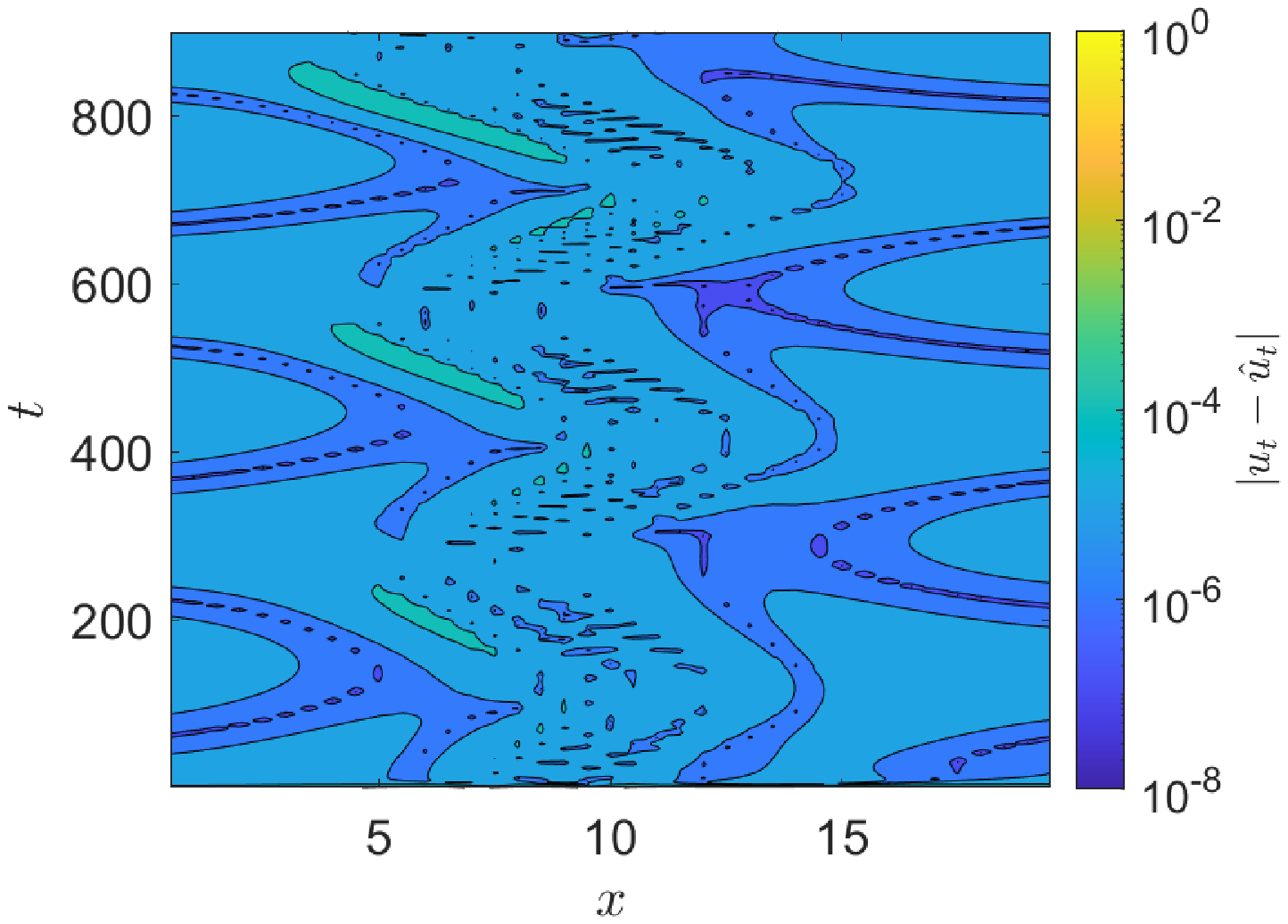}
    }
    \subfigure[]{
    \includegraphics[width=0.47 \textwidth]{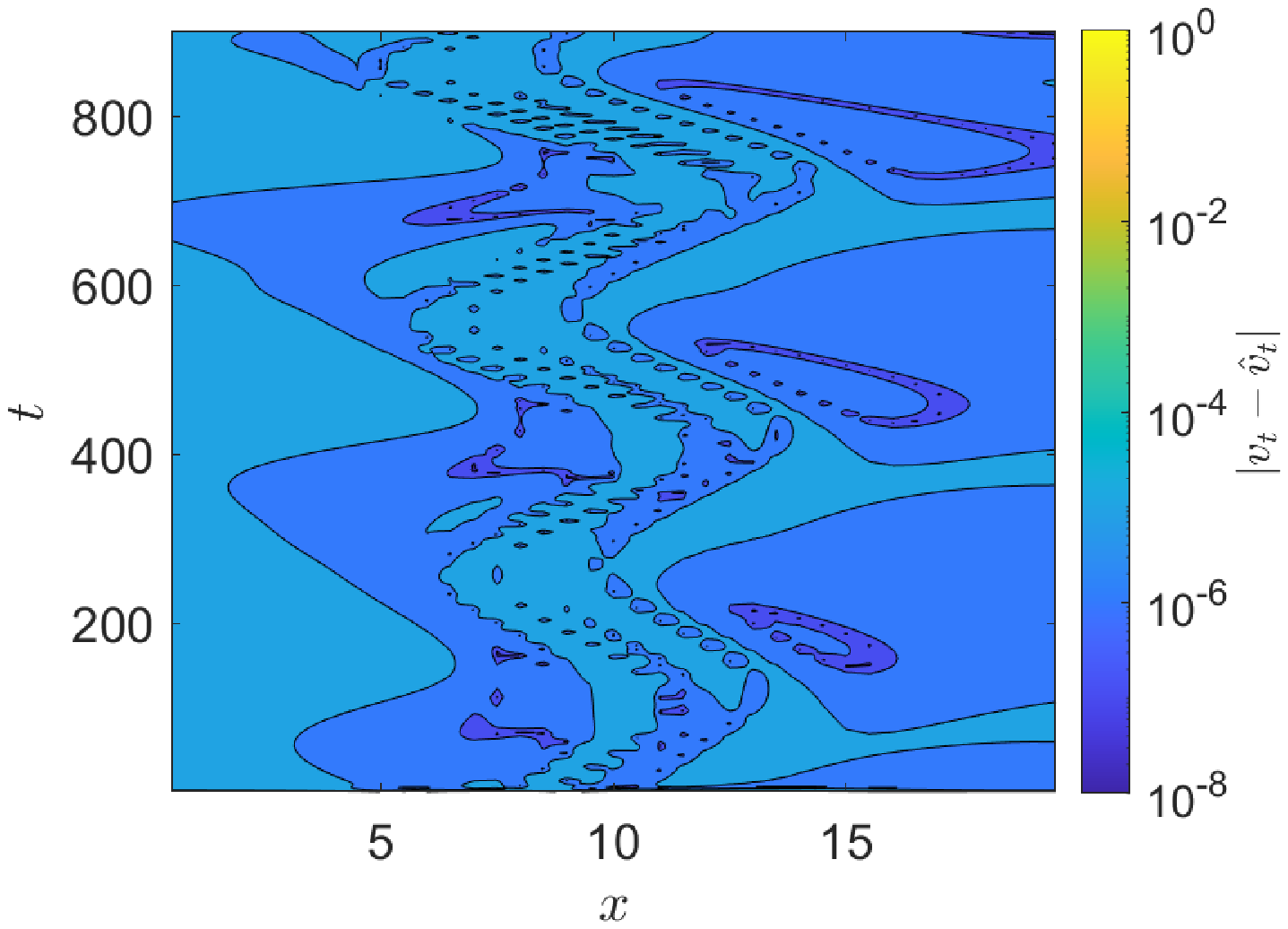}
    }
    \subfigure[]{
    \includegraphics[width=0.47 \textwidth]{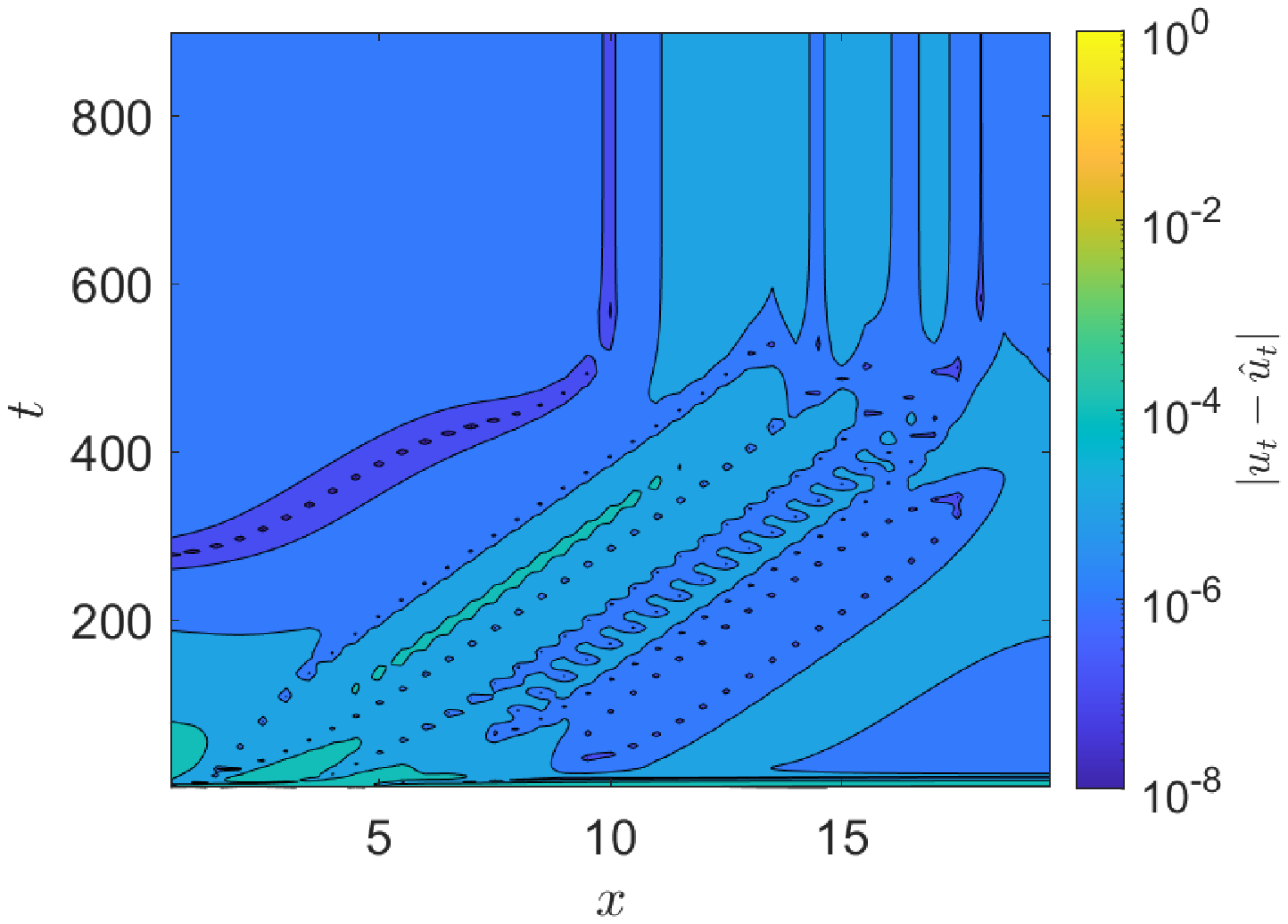}
    }
    \subfigure[]{
    \includegraphics[width=0.47 \textwidth]{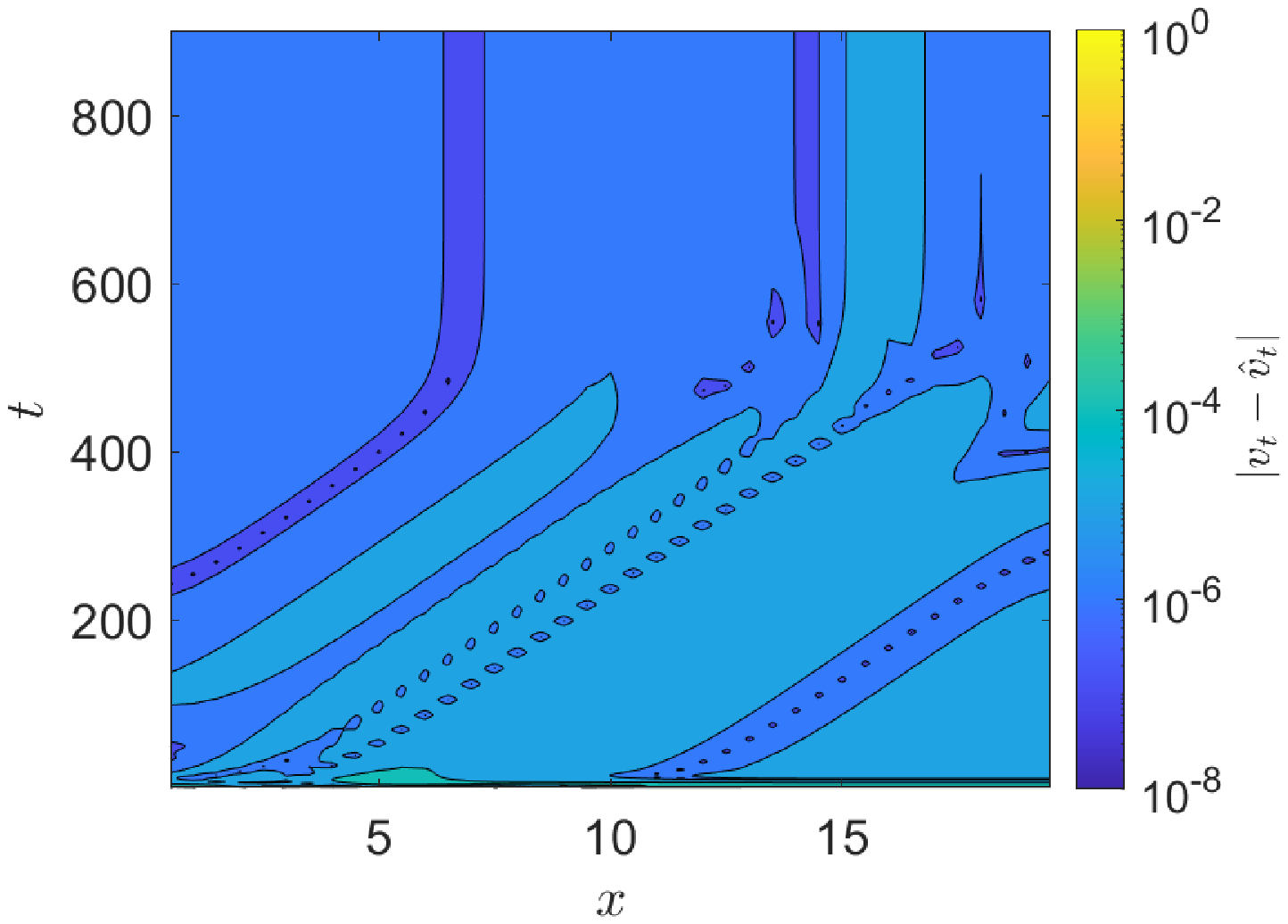}
    }
    \subfigure[]{
    \includegraphics[width=0.47 \textwidth]{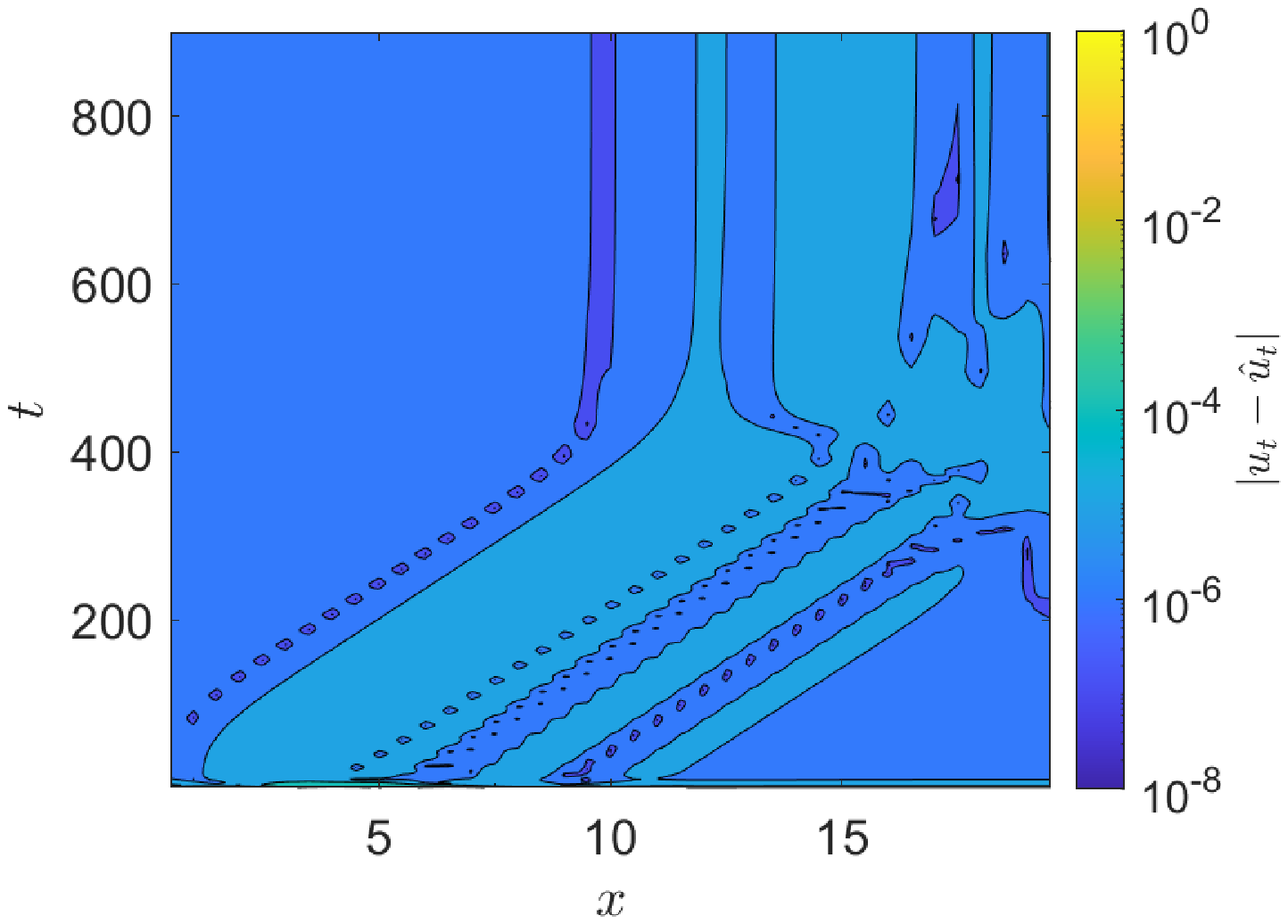}
    }
    \subfigure[]{
    \includegraphics[width=0.47 \textwidth]{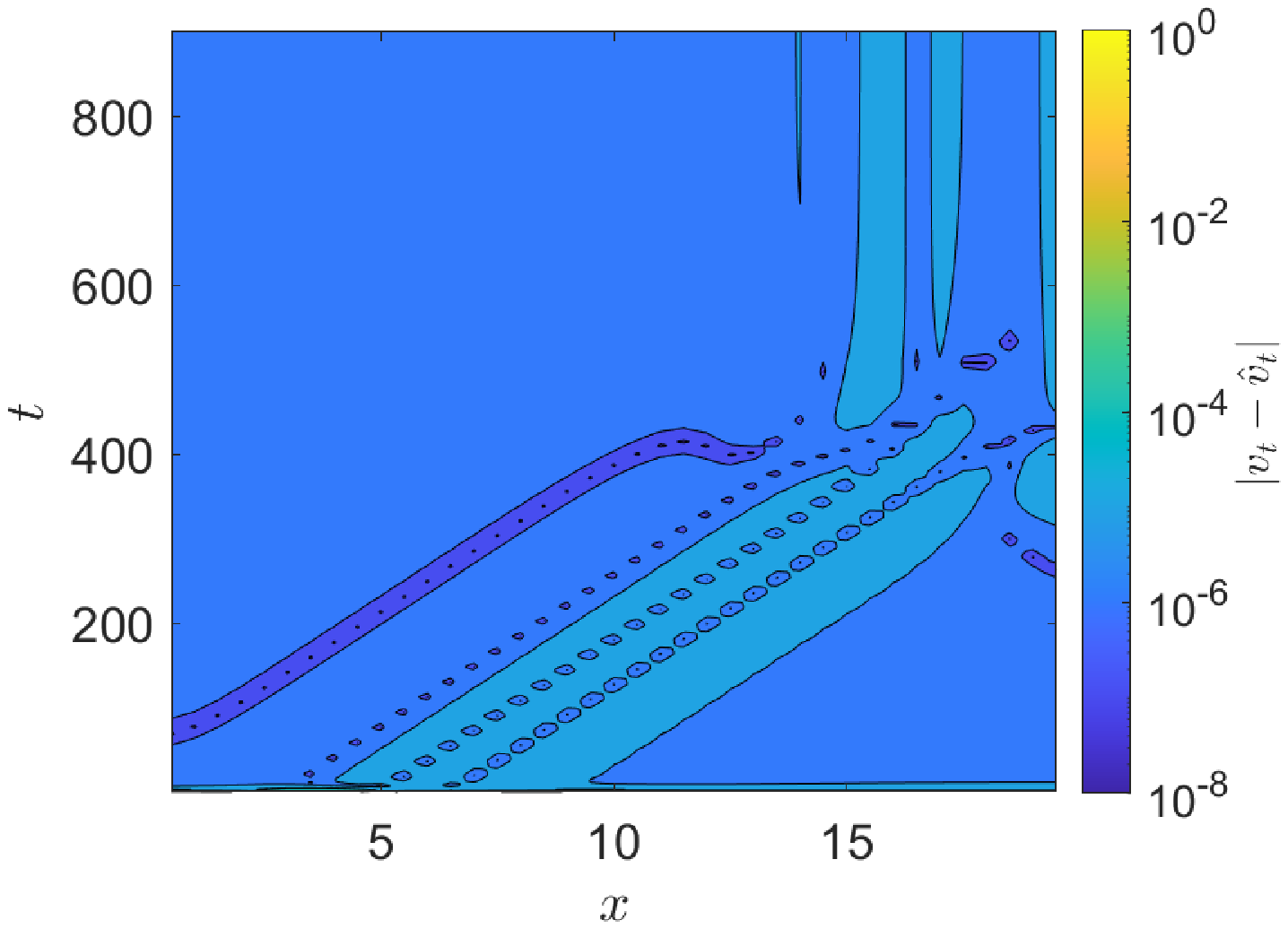}
    }
    \caption{Approximation accuracy in the test data without feature selection as obtained with FNNs. Contour plot of the absolute values of differences in space and time, of $|u_t(x,t)-\hat{u}_t(x,t)|$((a), (c), and (e)) and of $|v_t(x,t)-\hat{v}_t(x,t)|$ ((b), (d), and (f)) for characteristic values of $\varepsilon$: (a) and (b) $\varepsilon=0.0114$ near the Andronov-Hopf point, (c), (d) $\varepsilon=0.4$, (e) and (f) $\varepsilon=0.9383$ near the turning point.}
    \label{FNNsdifwithout}
\end{figure}

\begin{figure}
    \centering
    \subfigure[]{
    \includegraphics[width=0.47 \textwidth]{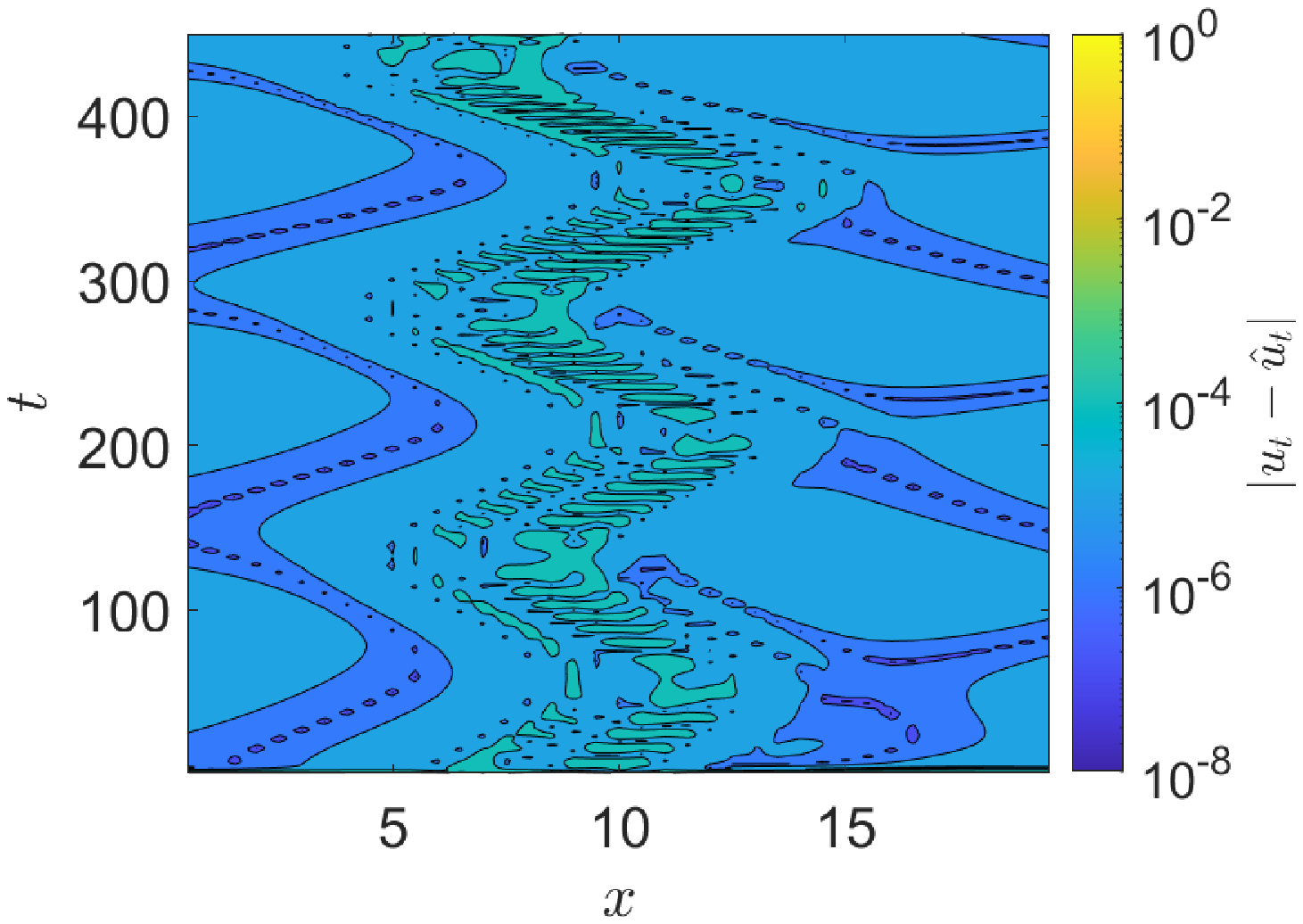}
    }
    \subfigure[]{
    \includegraphics[width=0.47 \textwidth]{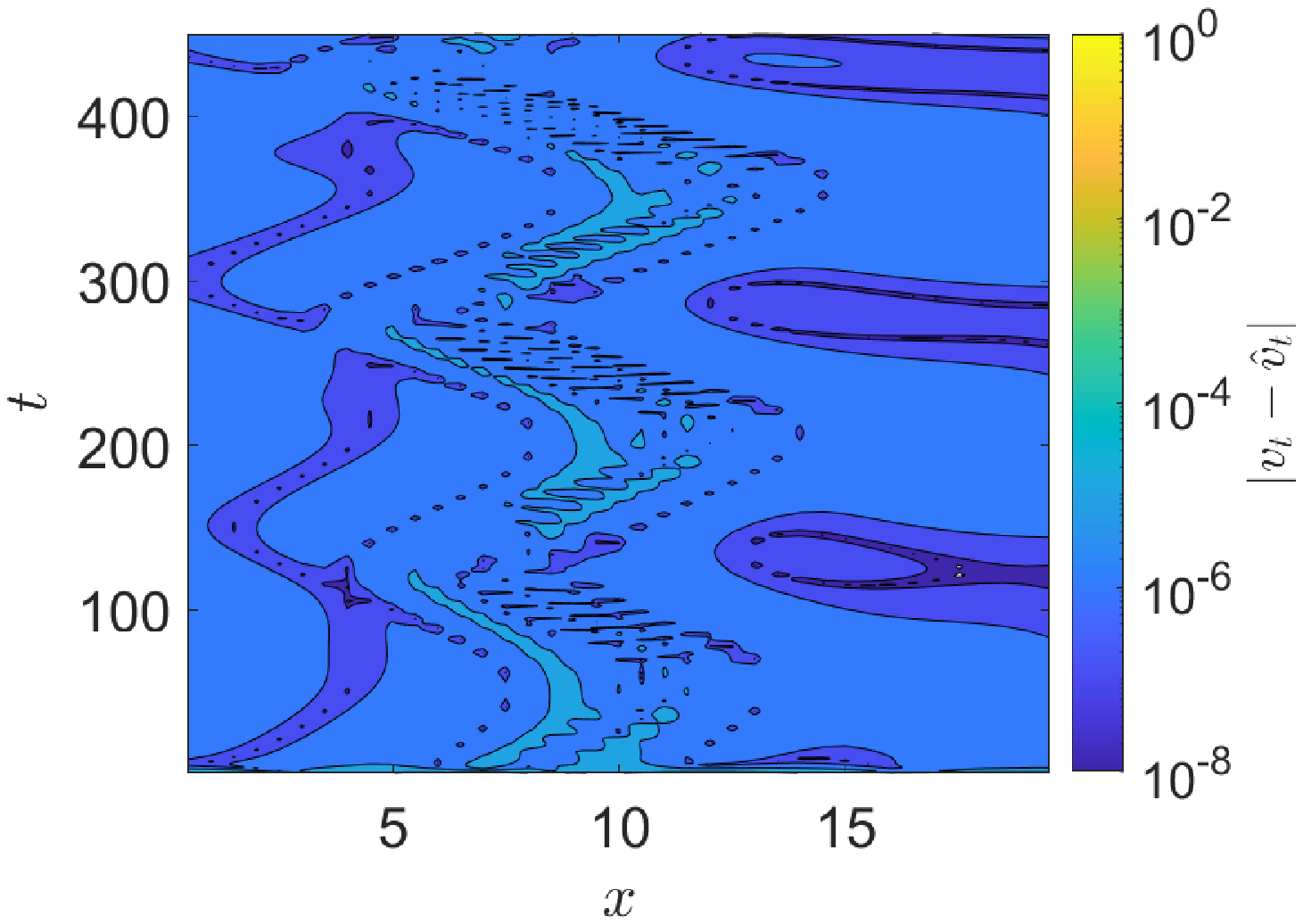}
    }
    \subfigure[]{
    \includegraphics[width=0.47 \textwidth]{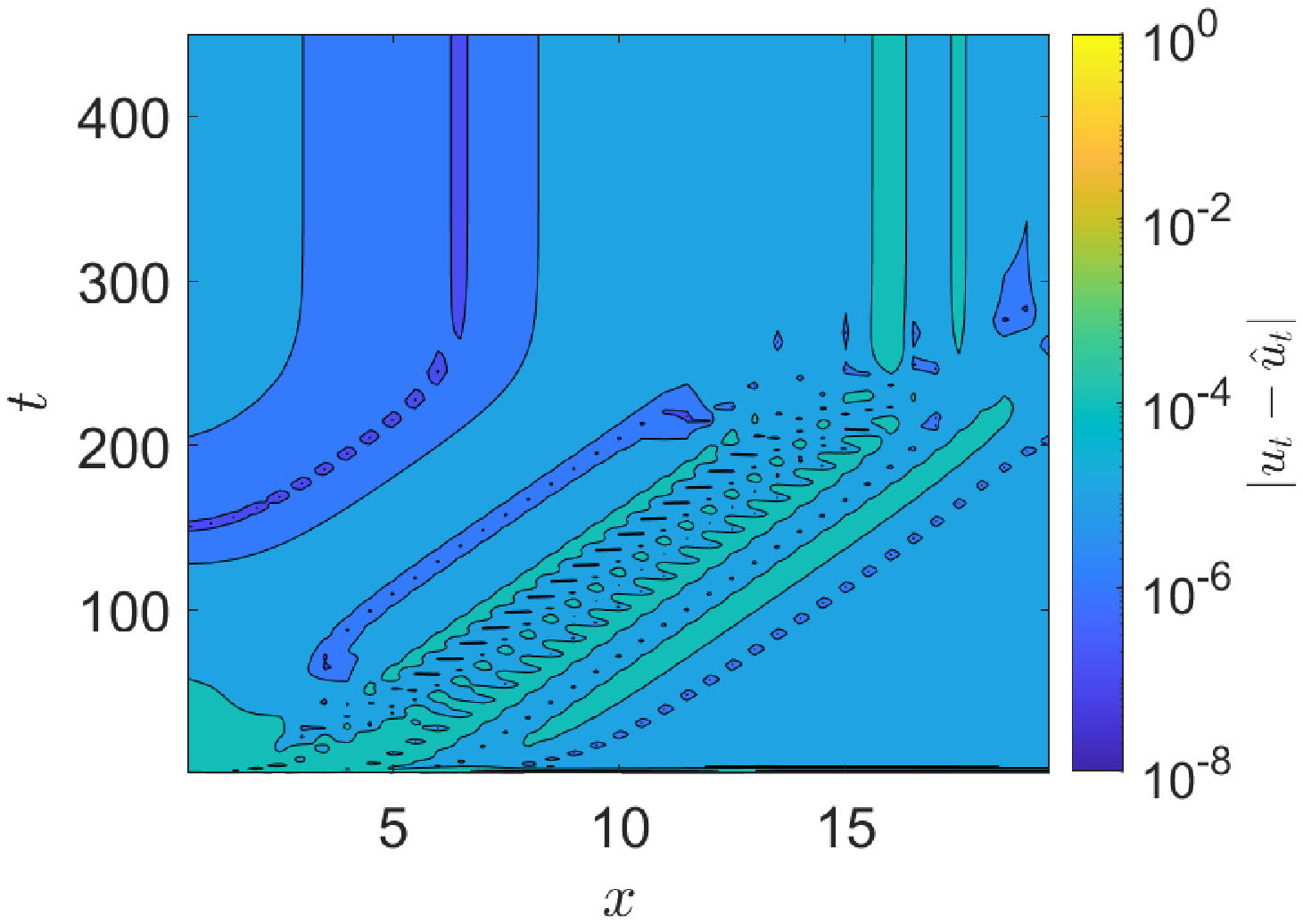}
    }
    \subfigure[]{
    \includegraphics[width=0.47 \textwidth]{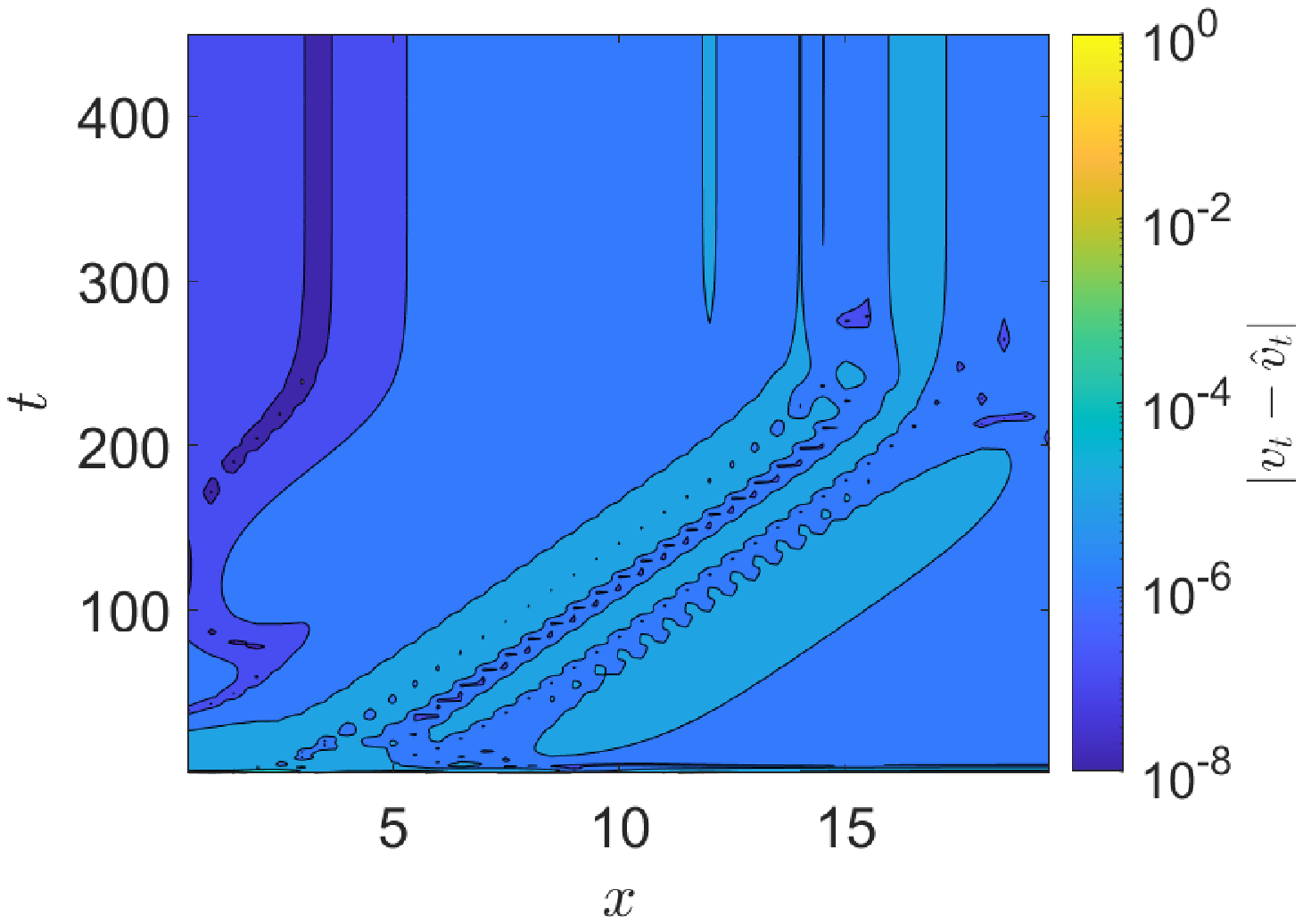}
    }
    \subfigure[]{
    \includegraphics[width=0.47 \textwidth]{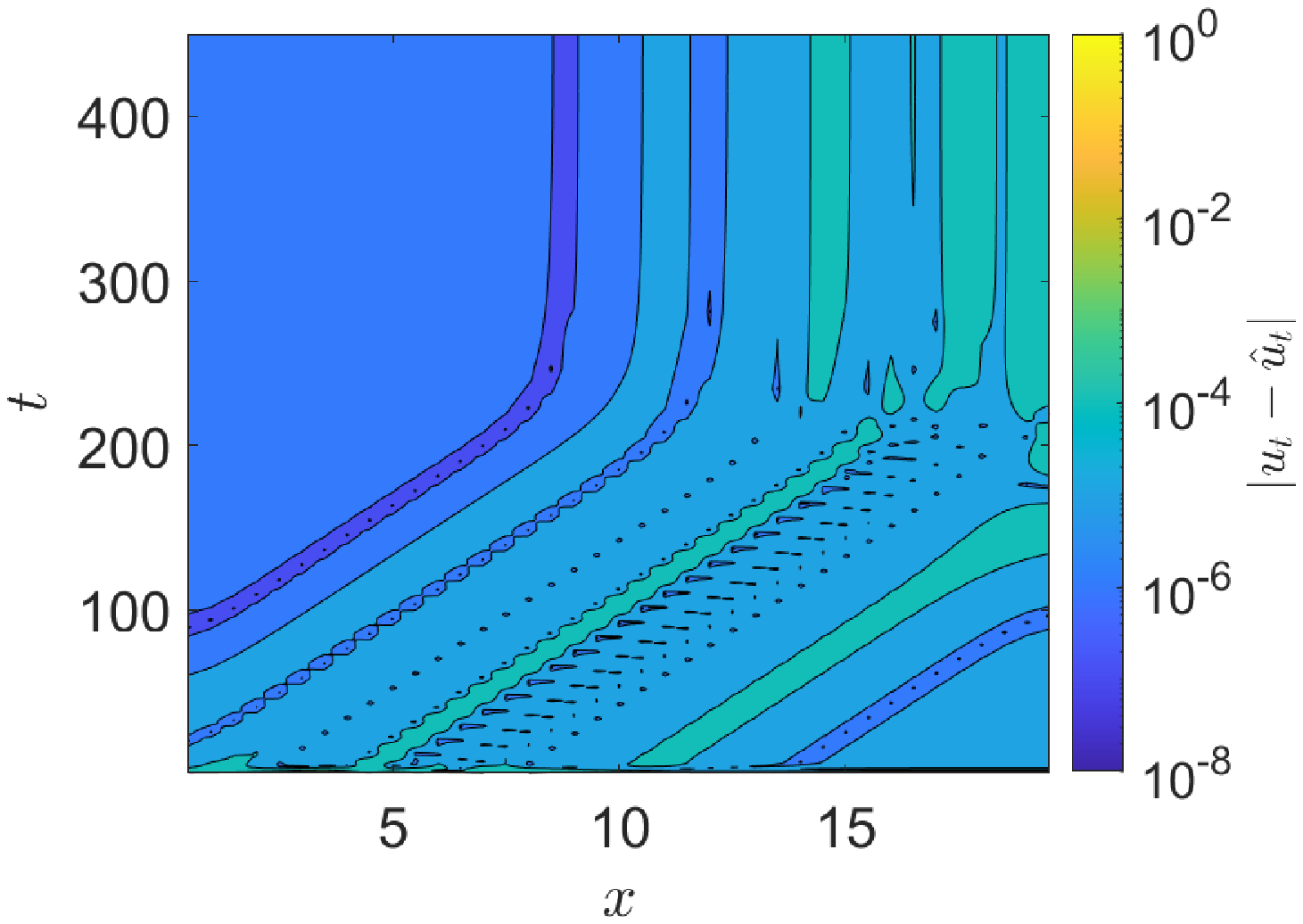}
    }
    \subfigure[]{
    \includegraphics[width=0.47 \textwidth]{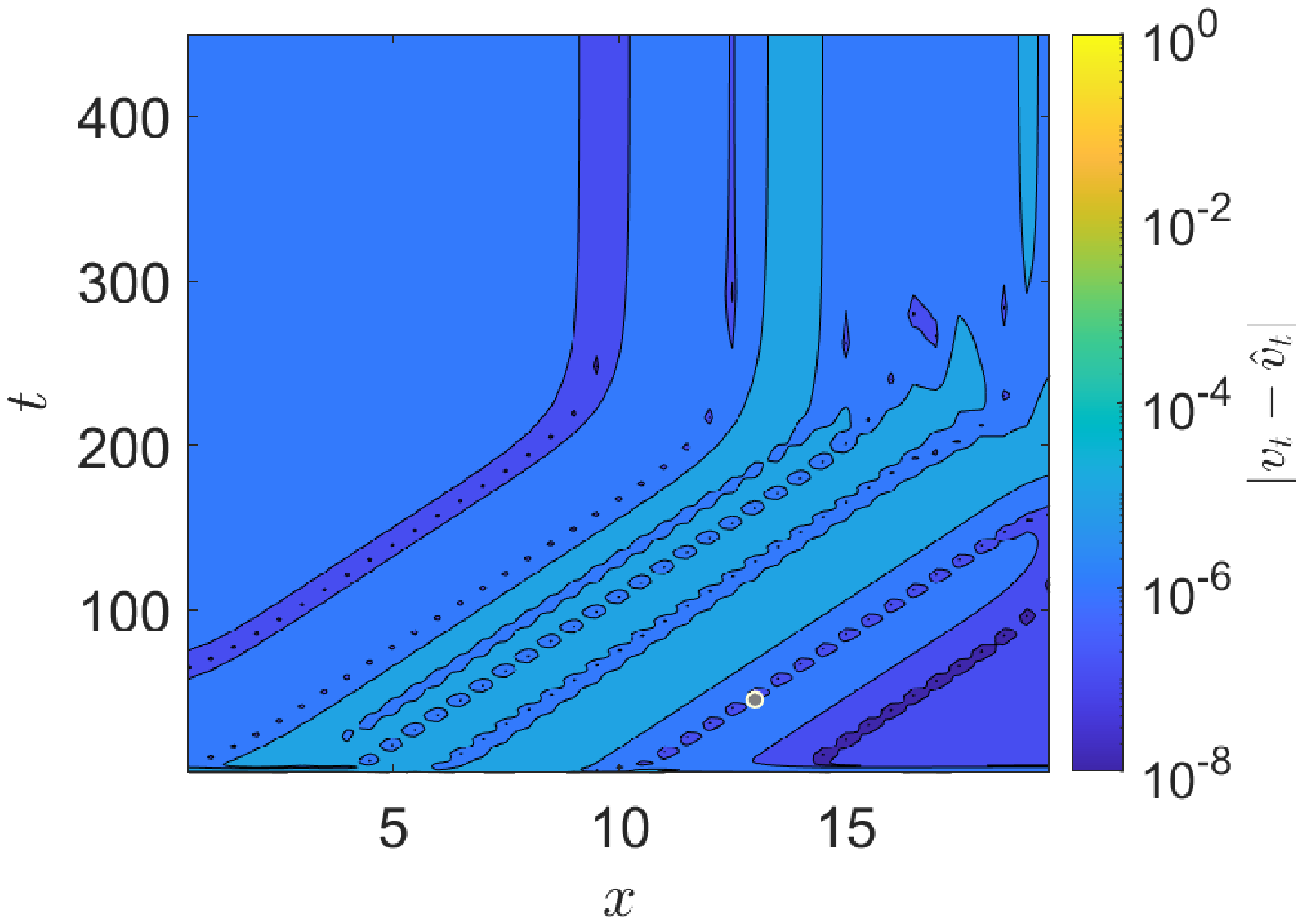}
    }
    \caption{Approximation accuracy in the test data without feature selection as obtained with RPNNs. Contour plot of the absolute values of differences in space and time, of $|u_t(x,t)-\hat{u}_t(x,t)|$ ((a), (c), and (e)) and of $|v_t(x,t)-\hat{v}_t(x,t)|$ ((b), (d), and (f)) for characteristic values of $\varepsilon$: (a) and (b) $\varepsilon=0.0114$ near the Andronov-Hopf point, (c), (d) $\varepsilon=0.4$, (e) and (f) $\varepsilon=0.9383$ near the turning point.}
    \label{RPNNsdifwithout}
\end{figure}

Differences between the predicted $\hat{u}_t(x,t)$ and $\hat{v}_t(x,t)$ and the actual values of the time derivatives $u_t(x,t)$ and $v_t(x,t)$ for three different values of $\varepsilon$ are shown in Figure \ref{FNNsdifwithout} when using FNNs and in Figure \ref{RPNNsdifwithout} when using RPNNs. 

\begin{figure}
    \centering
    \subfigure[]{
    \includegraphics[width=0.47  \textwidth]{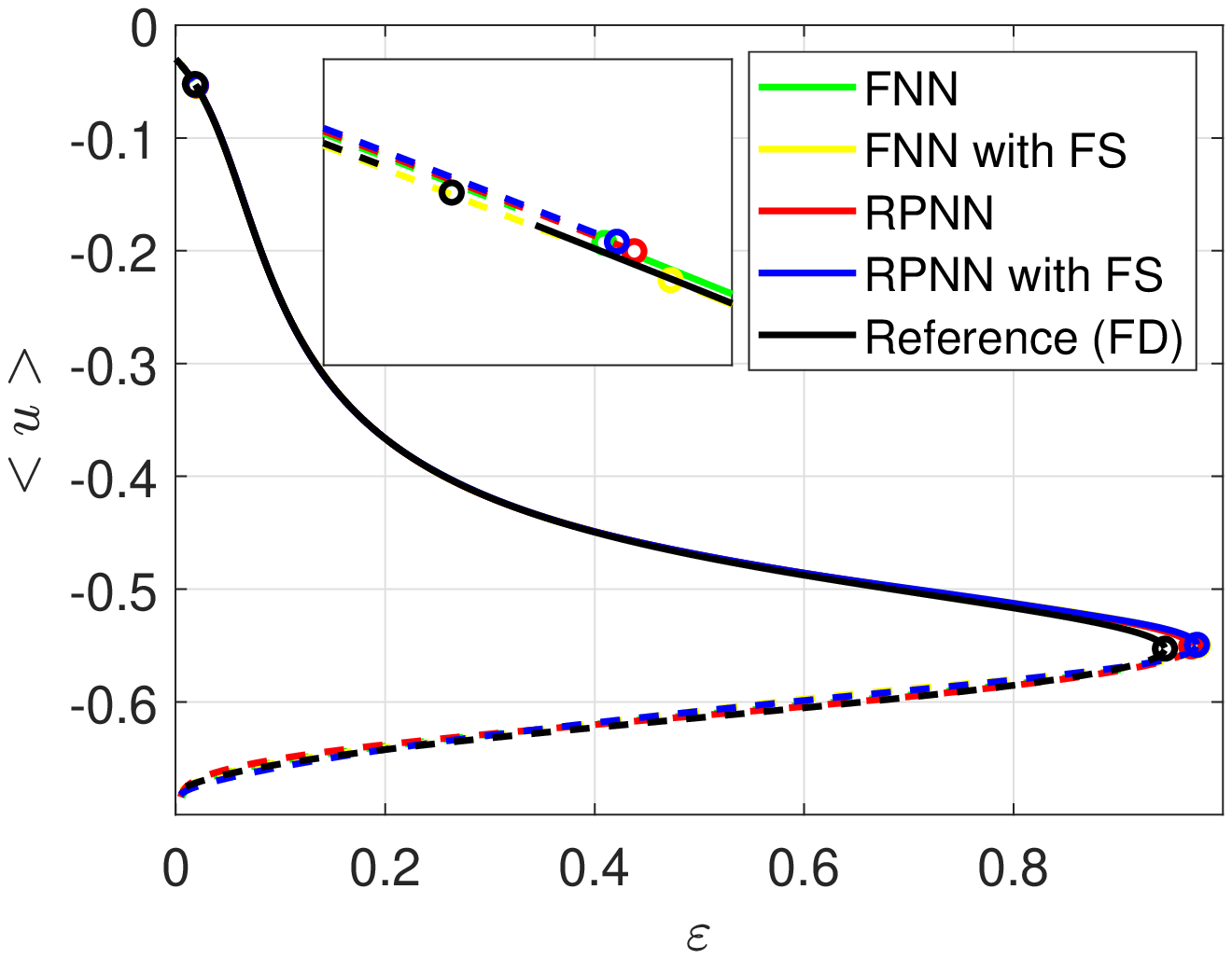}
    }
    \subfigure[]{
    \includegraphics[width=0.47\textwidth]{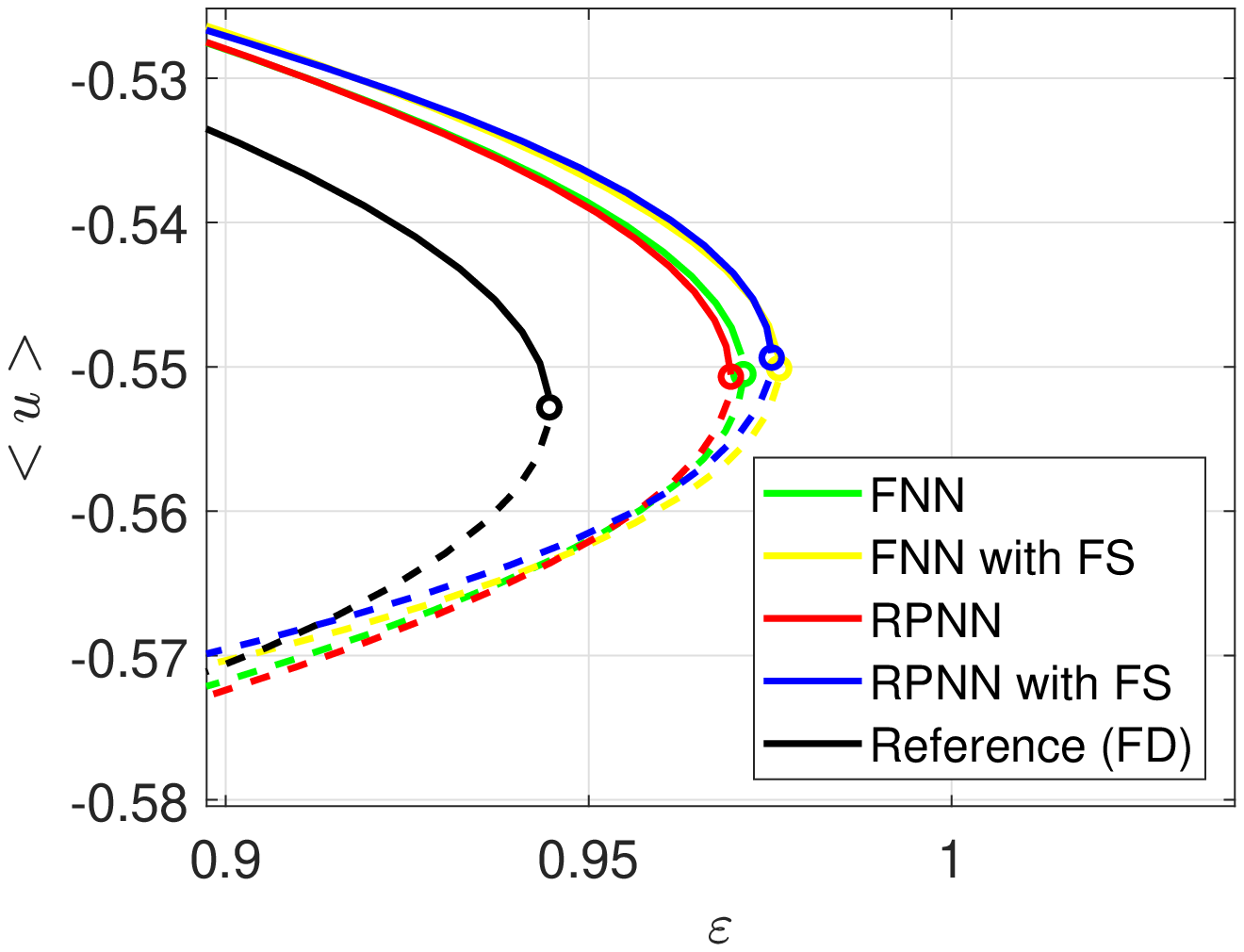}
    }
    \subfigure[]{
    \includegraphics[width=0.47\textwidth]{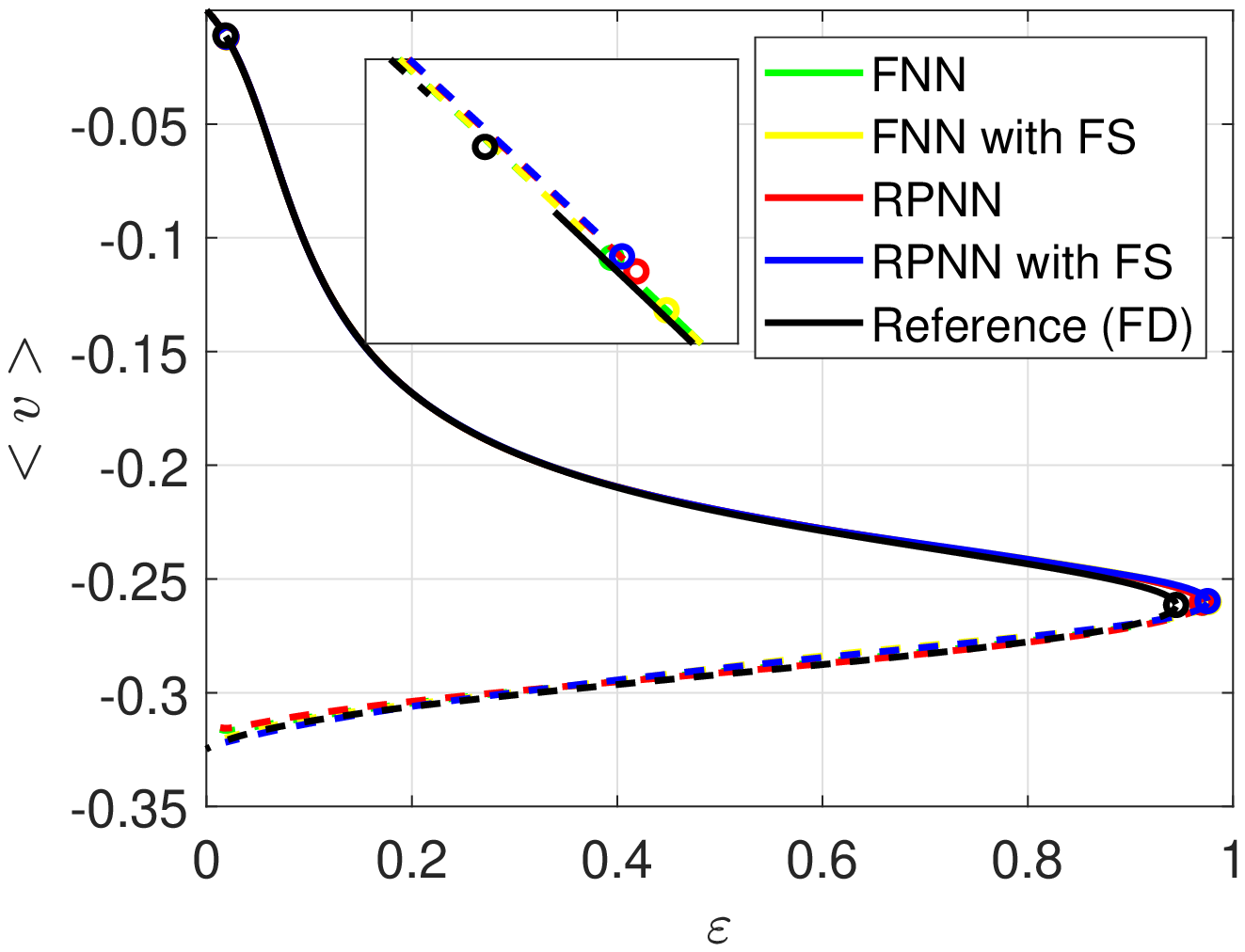}
    }
    \subfigure[]{
    \includegraphics[width=0.47\textwidth]{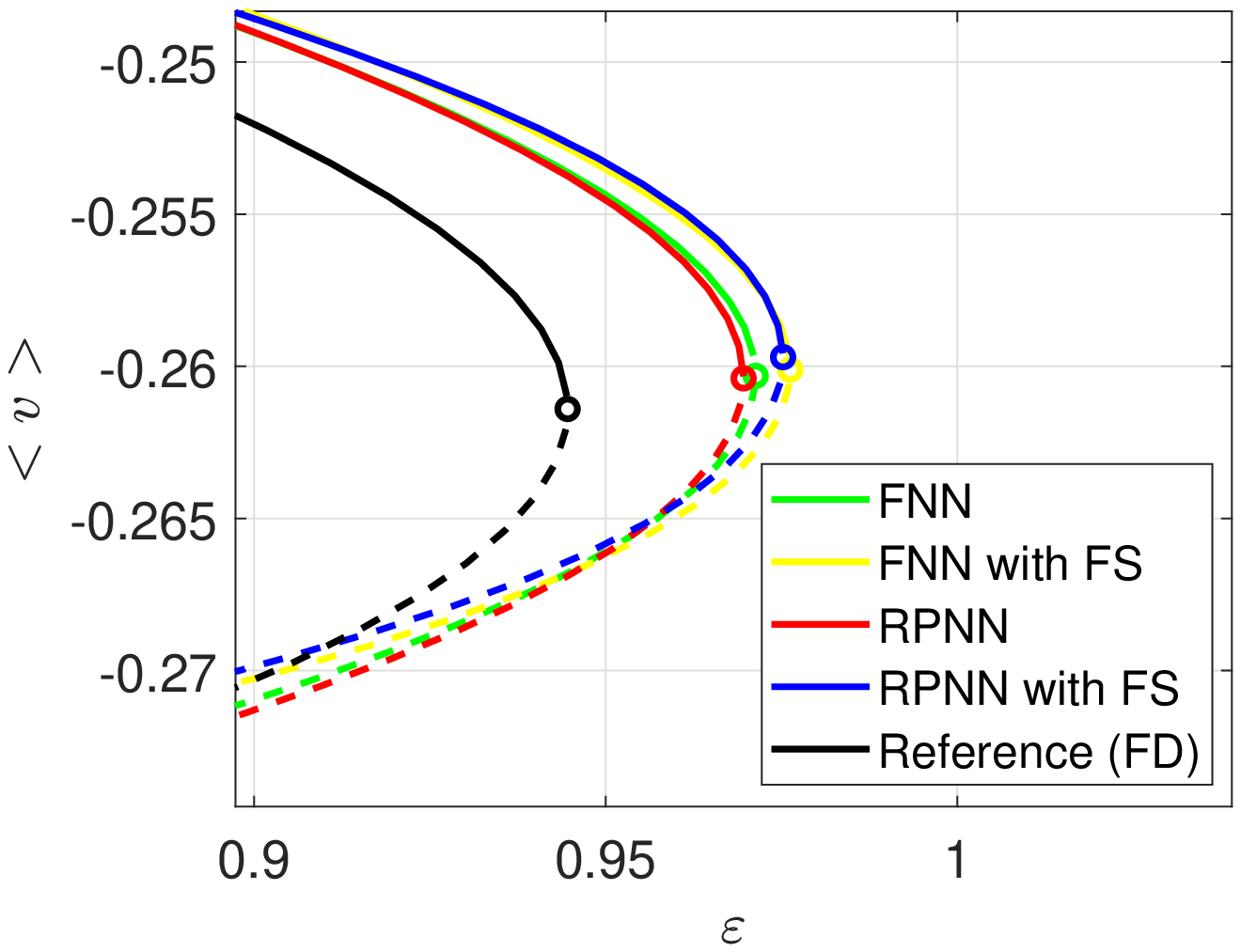}
    }
    \caption{Reconstructed bifurcation diagrams from the Lattice Boltzmann simulations of the FHN dynamics with respect to $\varepsilon$ with FNNs and RPNNs with and without feature selection. (a) Mean values $<u>$ for stable and unstable branches; the inset zooms near the Andronov-Hopf bifurcation point (b) zoom near the turning Point for $<u>$, (c) Mean values $<v>$ for stable and unstable branches; the inset zooms near the Andronov-Hopf bifurcation point, (d) zoom near the turning Point for $<v>$.}
    \label{bifdiag_final}
\end{figure}

Instead, for the proposed RPNN scheme, the training phase, i.e. the solution of the least-squares problem with regularization,  required around $8$ minutes, thus resulting to a training phase of at least 20 times faster than that of the FNNs.\\
After training, we used the FNNs and RPNNs to compute with finite differences the quantities required for performing the bifurcation analysis (see Eq.(\ref{augmentedarclength})), i.e.:
\begin{equation*}
\begin{aligned}
    \frac{\partial \hat{F}^u}{\partial u_j}=\frac{\hat{F}^u(u_j,v_j,\varepsilon)-\hat{F}^u(u_j+\delta,v_j,\varepsilon)}{2\delta}; 
    \frac{\partial \hat{F}^u}{\partial v_j}=\frac{\hat{F}^u(u_j,v_j,\varepsilon)-\hat{F}^u(u_j,v_j+\delta,\varepsilon)}{2\delta}\\
    \frac{\partial \hat{F}^v}{\partial u_j}=\frac{\hat{F}^v(u_j,v_j,\varepsilon)-\hat{F}^v(u_j+\delta,v_j,\varepsilon)}{2\delta}; 
    \frac{\partial \hat{F}^v}{\partial v_j}=\frac{\hat{F}^v(u_j,v_j,\varepsilon)-\hat{F}^v(u_j,v_j+\delta,\varepsilon)}{2\delta}\\
    \frac{\partial \hat{F}^u}{\partial \varepsilon}=\frac{\hat{F}^u(u_j,v_j,\varepsilon)-\hat{F}^u(u_j,v_j,\varepsilon+\delta)}{2\delta};
    \frac{\partial \hat{F}^v}{\partial \varepsilon}=\frac{\hat{F}^v(u_j,v_j,\varepsilon)-\hat{F}^v(u_j,v_j,\varepsilon+\delta)}{2\delta},
\end{aligned}
\label{numJac}
\end{equation*}
with $\delta=1e-06$.
The reconstructed bifurcation diagrams are shown in Figure \ref{bifdiag_final}.
Using the FNNs, we estimated the Andronov-Hopf point at $\varepsilon\approx 0.0191$ and the turning point at $\varepsilon\approx 0.9713$; using the RPNNs, we estimated the Andronov-Hopf point at $\varepsilon\approx 0.0193$ and the turning point at $\varepsilon\approx 0.9696$. We approximated the same points using the finite differences scheme in the previous section at $\varepsilon\approx 0.0183$ for the Andronov-Hopf point and at $\varepsilon\approx 0.9446$ for the turning point. Hence, compared to the FNNs, the RPNNs approximated slightly better the reference turning point.

\subsection{Numerical bifurcation analysis with feature selection}
\begin{figure}
     \centering
     \subfigure[]{
         \centering
         \includegraphics[width=0.4\textwidth]{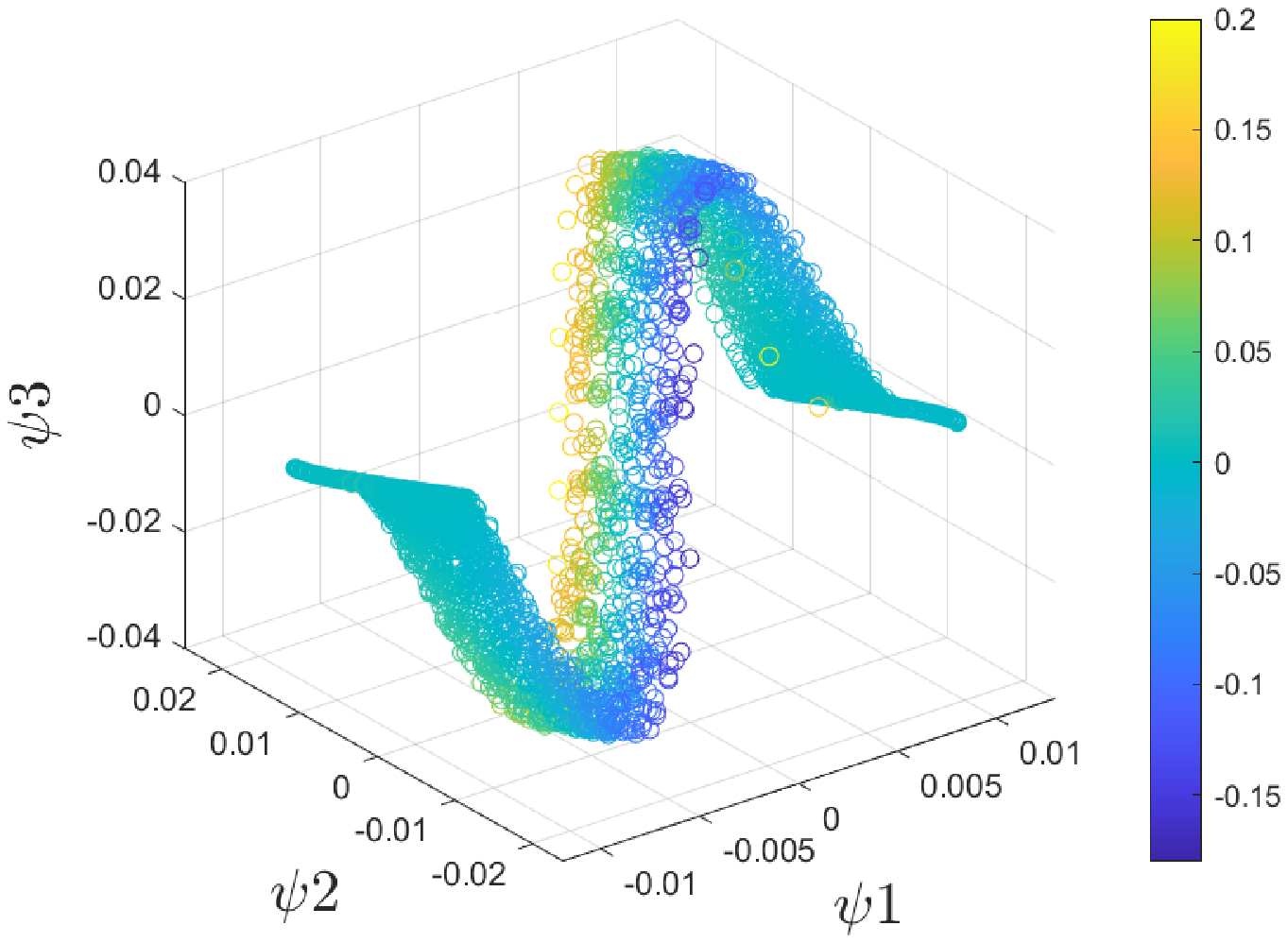}
         
         \label{udmap0.0114}
     }
     \hfill
     \centering
     \subfigure[]{
         \centering
         \includegraphics[width=0.4\textwidth]{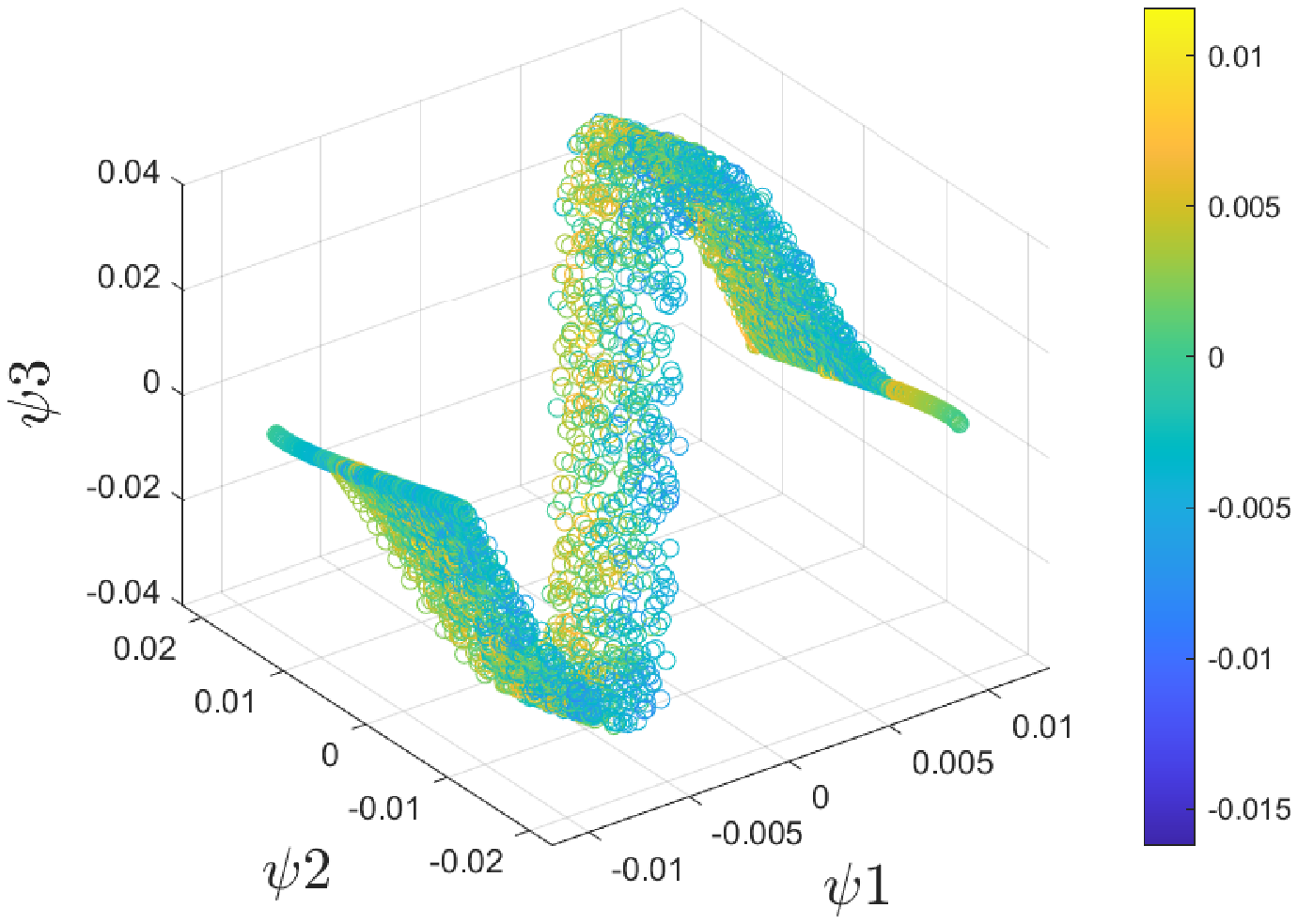}
         
         \label{vdmap0.0114}
     }
     \hfill
     \centering
     \subfigure[]{
         \centering
         \includegraphics[width=0.4\textwidth]{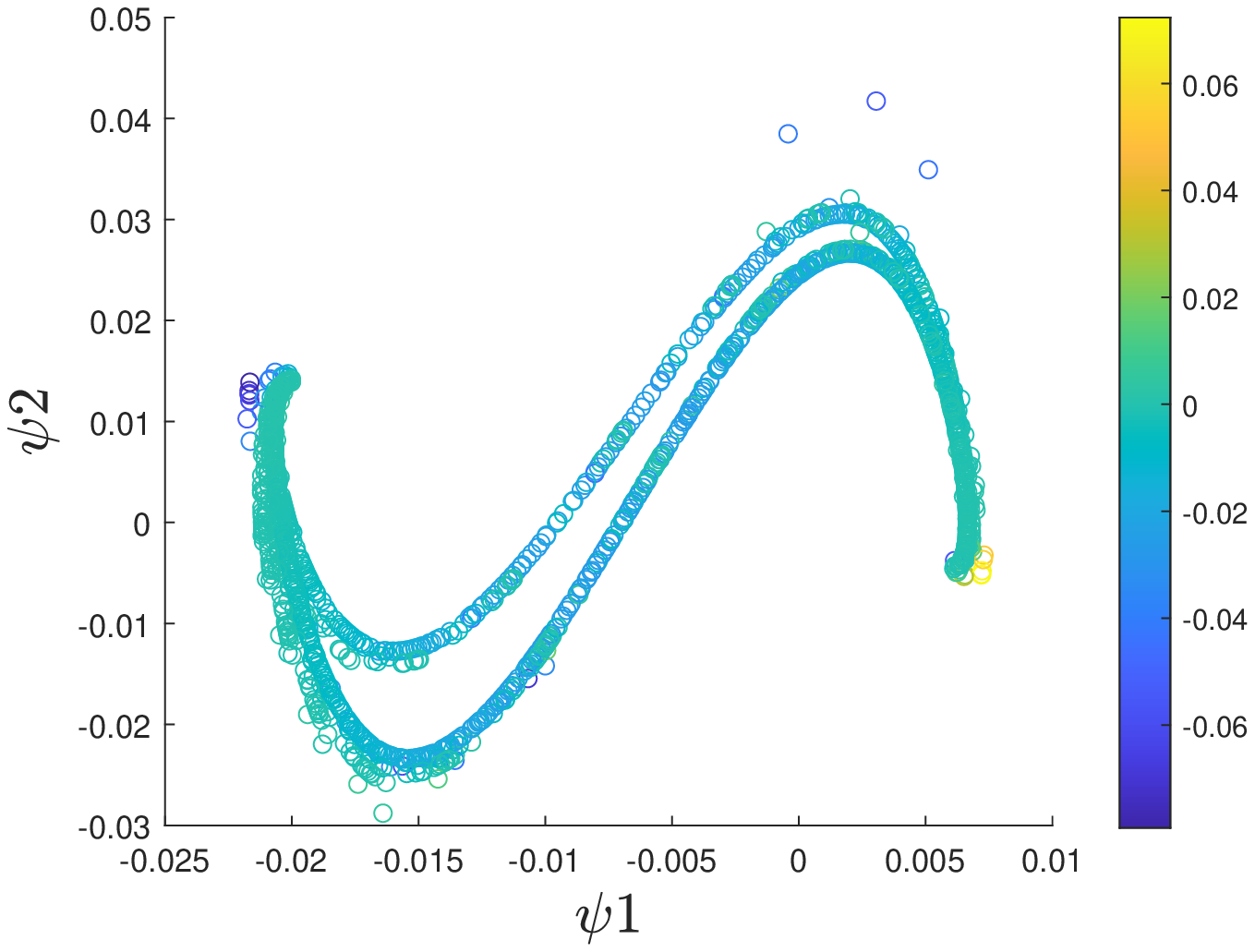}
         
         \label{udmap0.4010}
     }
     \hfill
     \centering
     \subfigure[]{
         \centering
         \includegraphics[width=0.4\textwidth]{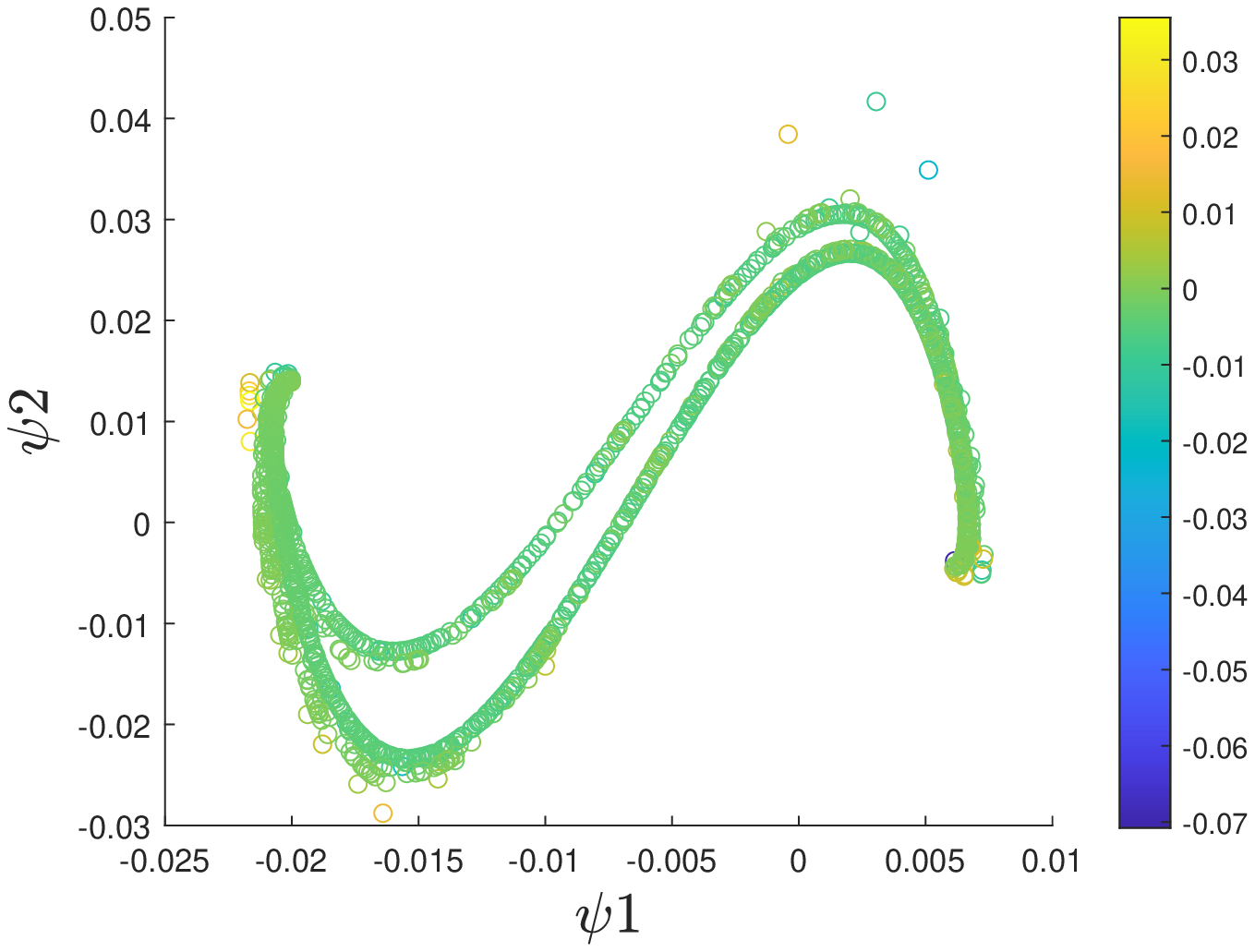}
         
         \label{vdmap0.4010}
     }
     \hfill
     \centering
     \subfigure[]{
         \centering
         \includegraphics[width=0.4\textwidth]{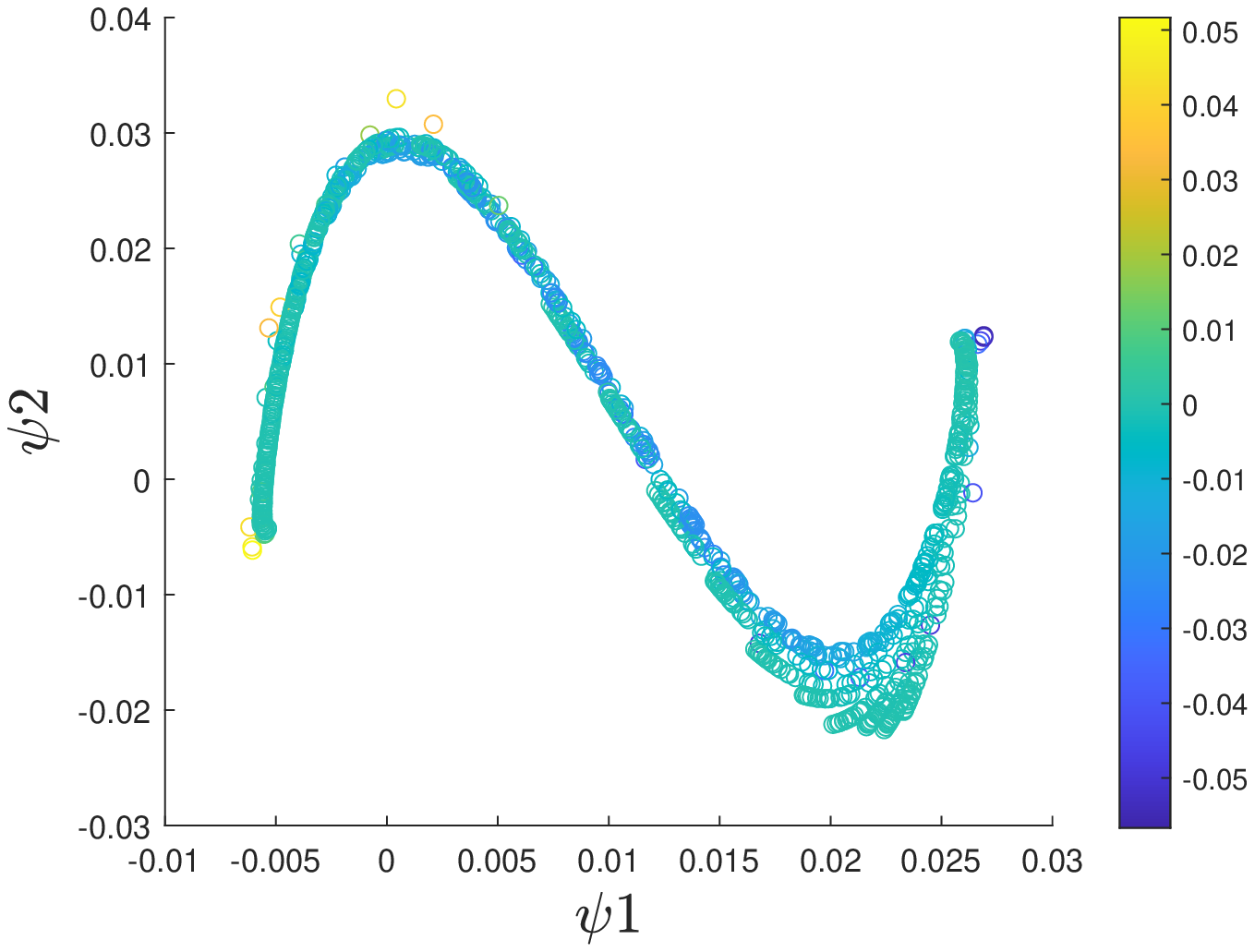}
         
         \label{udmap0.9554}
     }
     \hfill
     \centering
     \subfigure[]{
         \centering
         \includegraphics[width=0.4\textwidth]{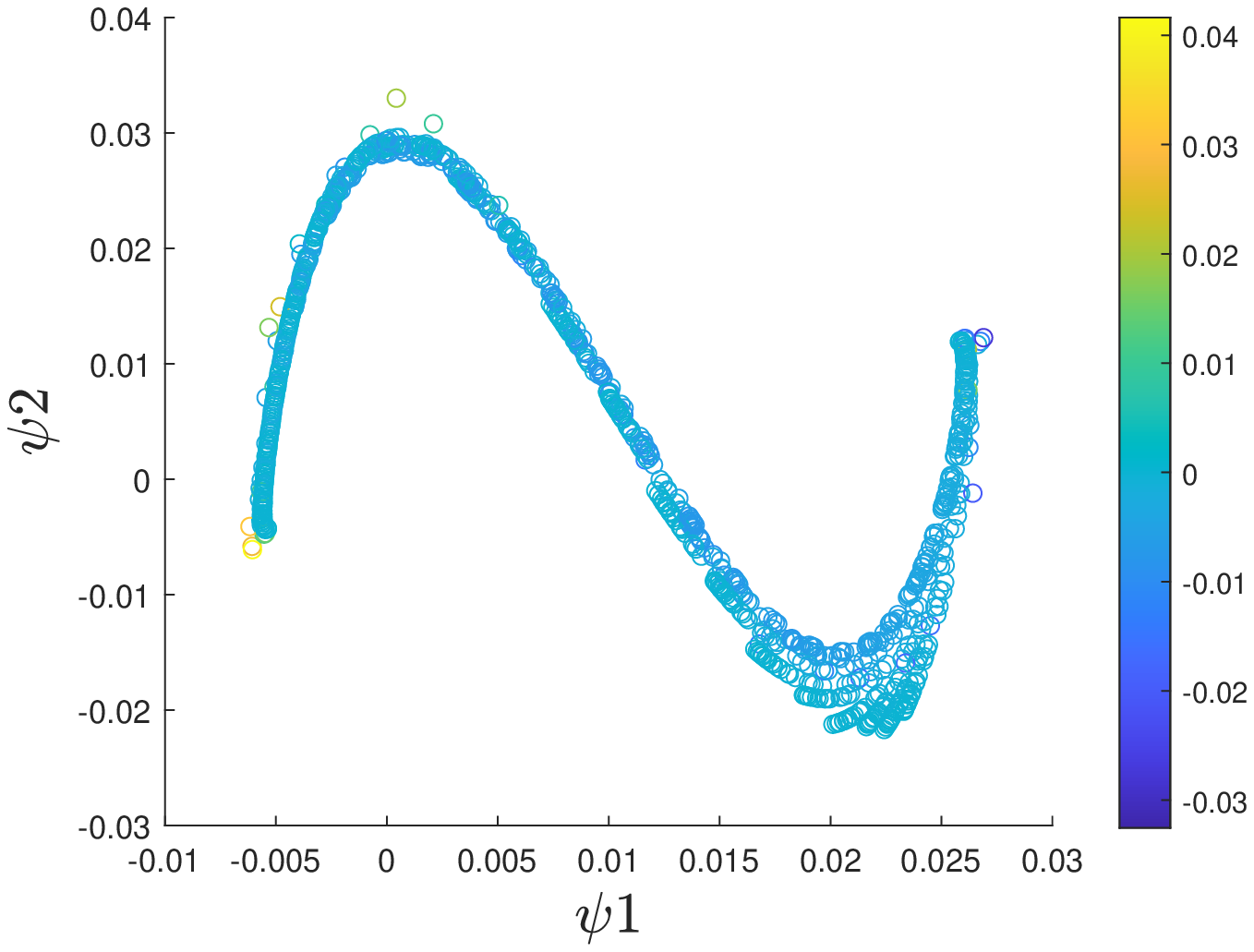}
         
         \label{vdmap0.9554}
     }
        \caption{(a) and (b): The three parsimonious Diffusion Maps coordinates for $\varepsilon=0.01114$ near the Andronov-Hopf point, respectively. (c) and (d): the two parsimonious Diffusion maps coordinates for $\varepsilon=0.4010$. (e) and (f): the two parsimonious Diffusion maps coordinates for $\varepsilon=0.9383$ near the turning point. Colors represent $u_t$ ((a), (c), (e)) and $v_t$ ((b), (d), (f)).}
        \label{difmaps}
\end{figure}

We used Diffusion Maps (setting the width parameter of the Gaussian kernel  to $\sigma=10$) to identify the three parsimonious leading eigenvectors as described in section \ref{sec:dmaps}. For our computations, we used the datafold package in python \cite{Lehmberg2020}. We denote them as $\phi_1,\phi_2,\phi_3$. The three parsimonious Diffusion Maps coordinates for different values of the parameter $\varepsilon$ are shown in Figure \ref{difmaps}. For $\varepsilon=0.114$ that is close to the Andronov-Hopf point, the embedded space is a two dimensional ``carpet'' in the three dimensional space. The oscillatory behaviour leads to different values of the time derivative which can be effectively parametrized as shown by the coloring of the manifold (Figures \ref{udmap0.0114}, \ref{vdmap0.0114}). For $\varepsilon=0.4010$ and $\varepsilon=0.9383$, the embedded space is a one dimensional line, since time derivatives converges rapidly to zero (Figures \ref{udmap0.4010},\ref{udmap0.9554},\ref{vdmap0.4010} and \ref{vdmap0.9554}).   
Based on the feature selection methodology, the ``good'' subsets of the input data domain are presented in Table \ref{reglos}. As expected, the best candidate features are the $(u,v,u_{xx})$ for $u_t$ and $(u,v,v_{xx})$ for $v_t$, which are the only features that indeed appear in the closed form of the FHN PDEs.

\begin{figure}
    \centering
    \subfigure[]{
    \includegraphics[width=0.47 \textwidth]{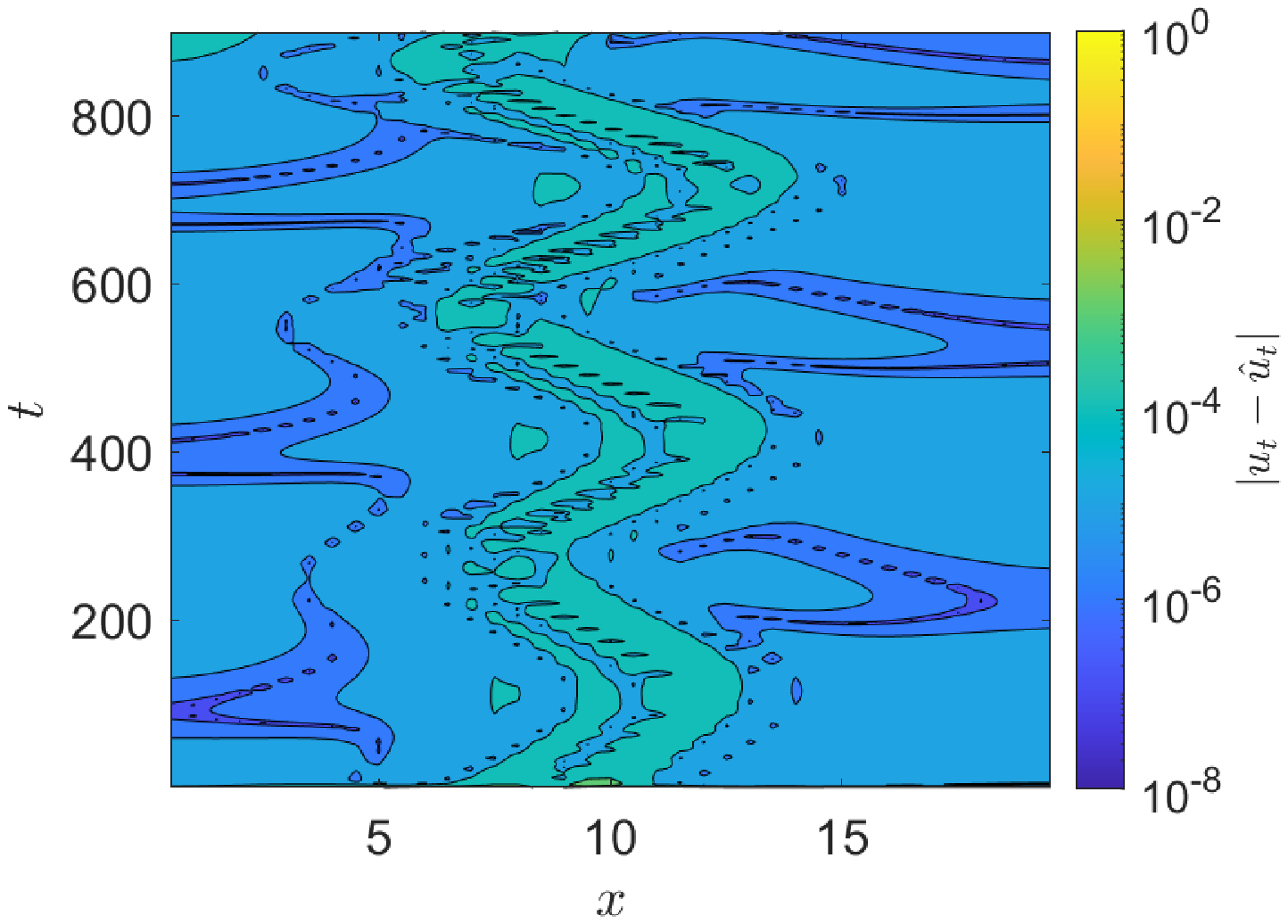}
    }
    \subfigure[]{
    \includegraphics[width=0.47 \textwidth]{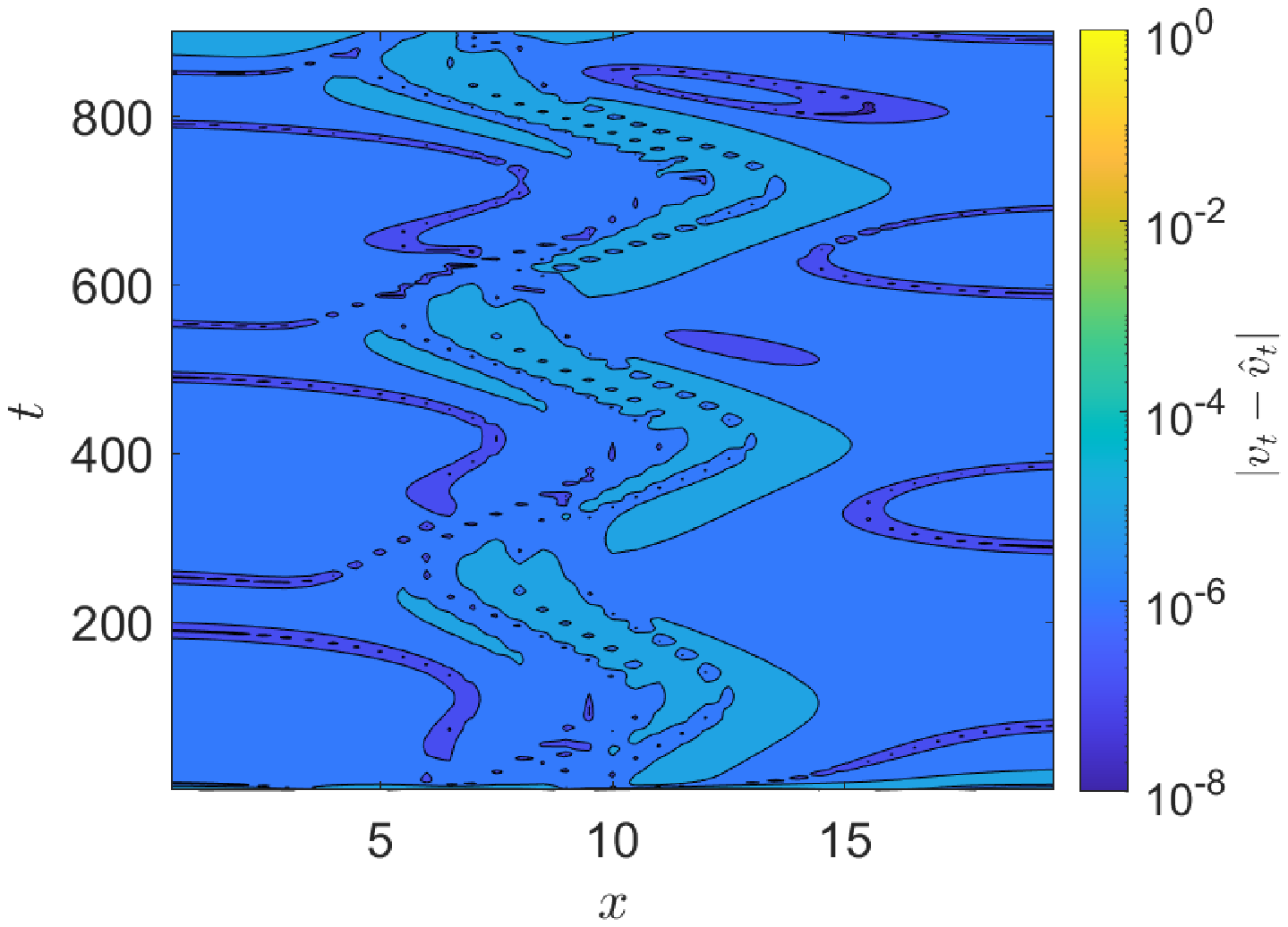}
    }
    \subfigure[]{
    \includegraphics[width=0.47 \textwidth]{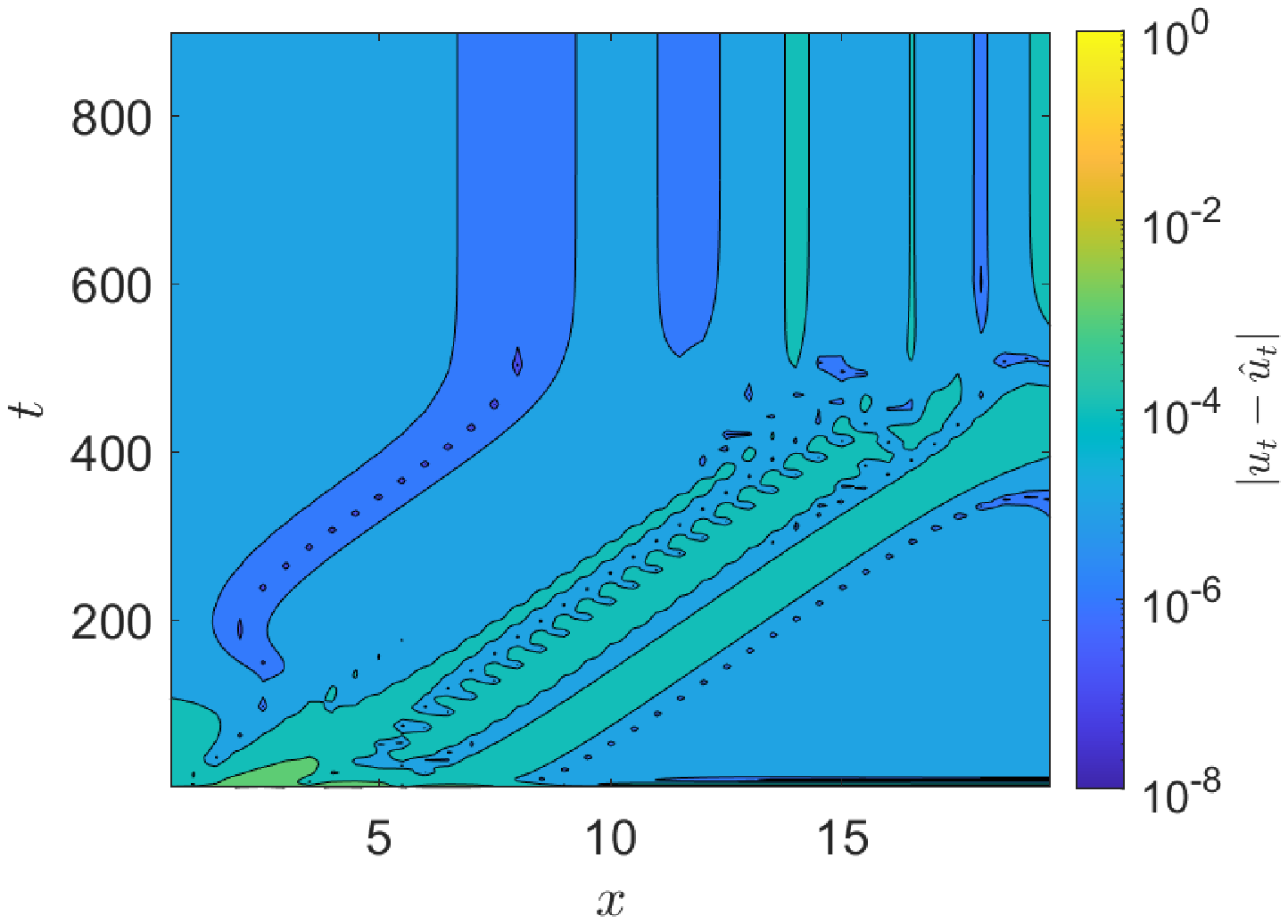}
    }
    \subfigure[]{
    \includegraphics[width=0.47 \textwidth]{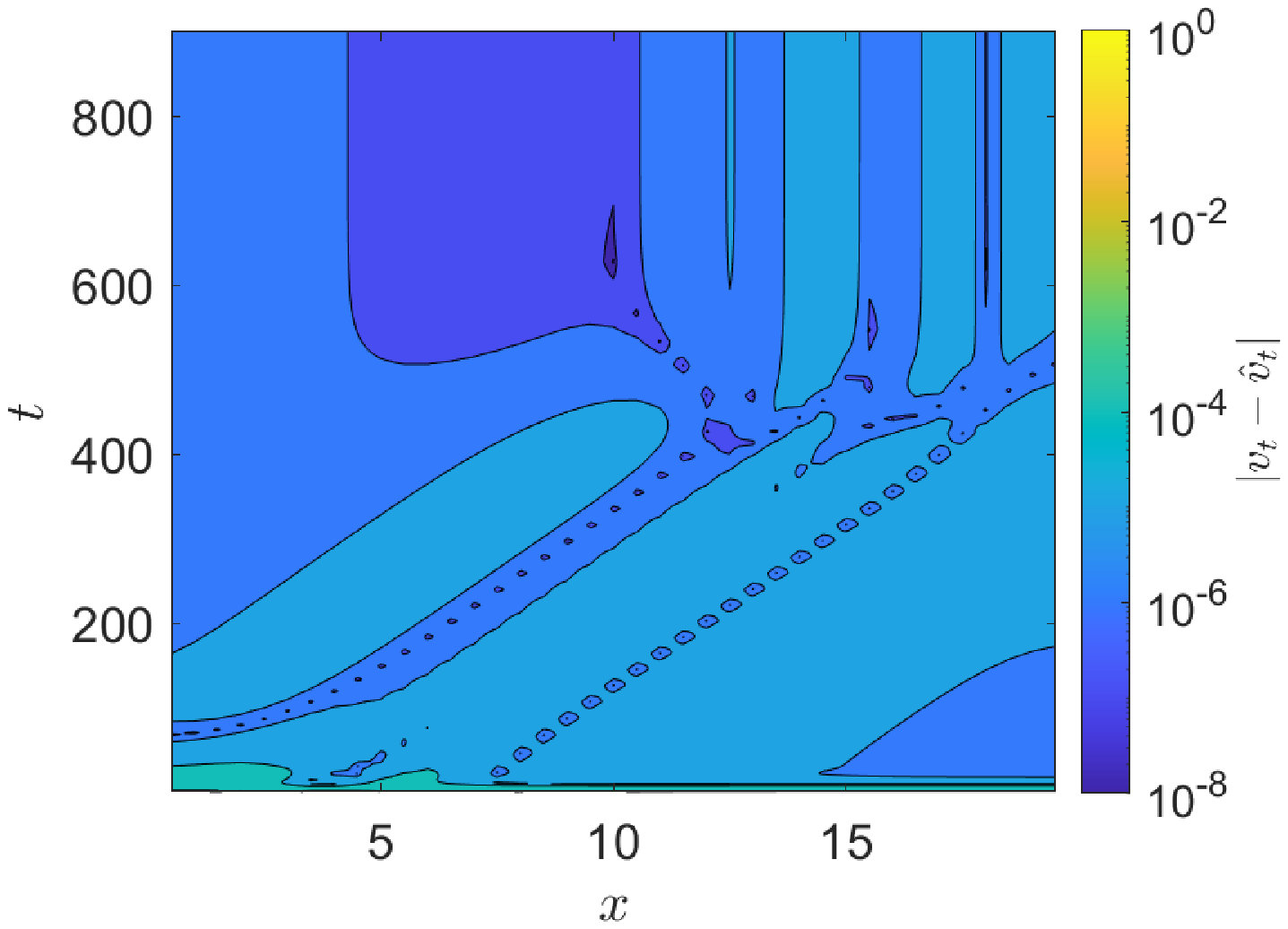}
    }
    \subfigure[]{
    \includegraphics[width=0.47 \textwidth]{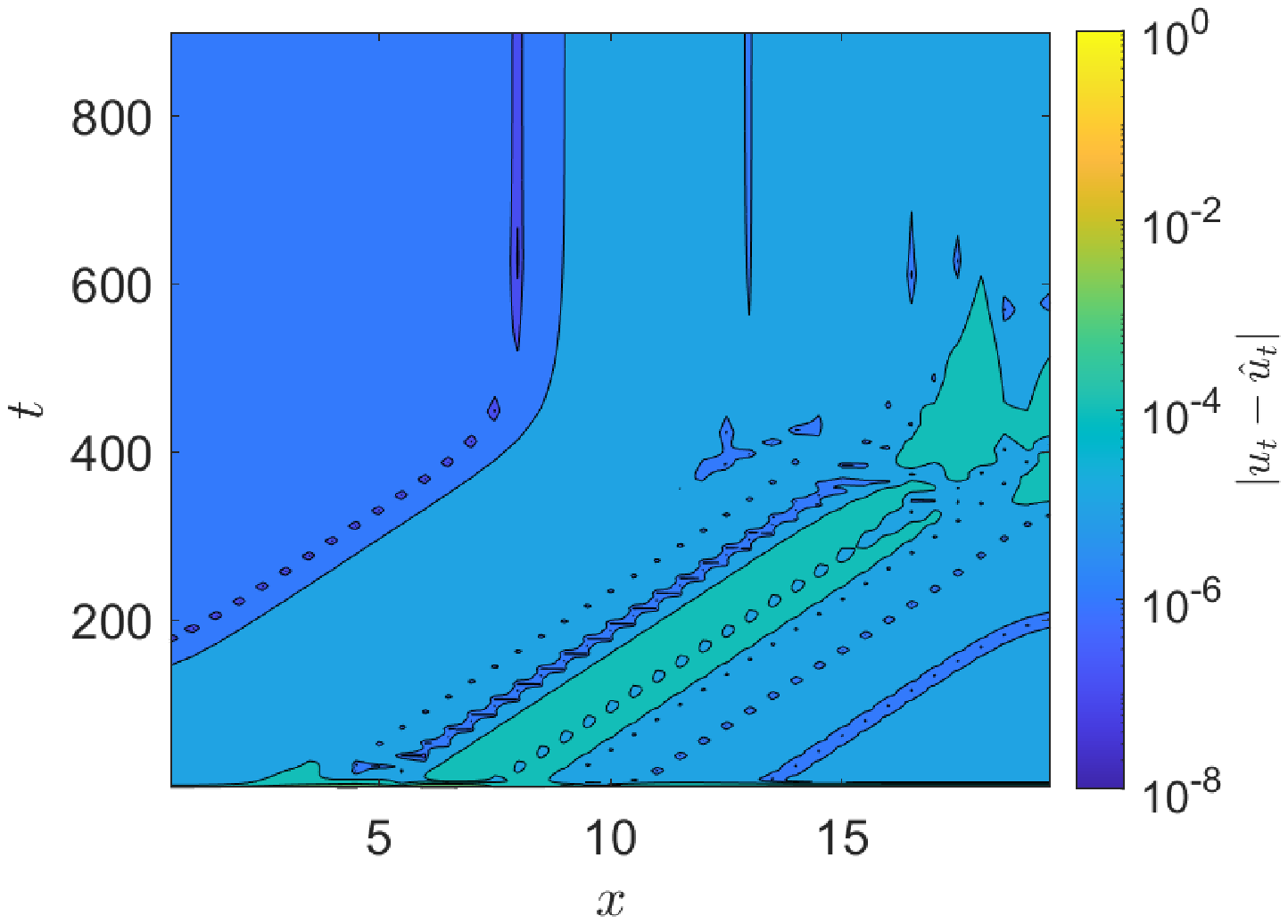}
    }
    \subfigure[]{
    \includegraphics[width=0.47 \textwidth]{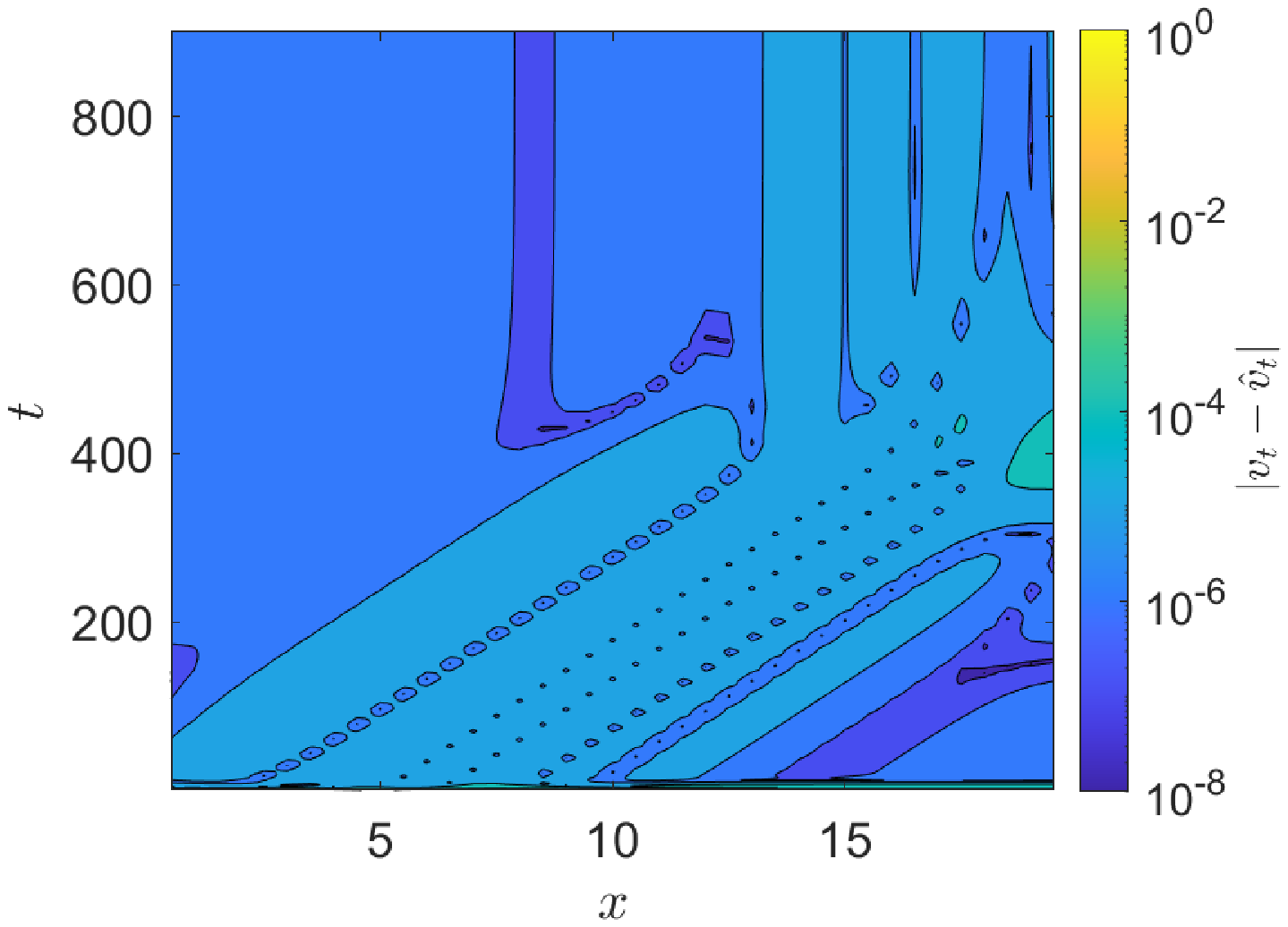}
    }
    \caption{Approximation accuracy in the test data with feature selection as obtained with the FNNs. Contour plot of the absolute values of differences in space and time, of
    $|u_t(x,t)-\hat{u}_t(x,t)|$ ((a), (c), and (e)) and of $|v_t(x,t)-\hat{v}_t(x,t)|$ ((b), (d), and (f)) for characteristic values of $\varepsilon$: (a) and (b) $\varepsilon=0.0114$ near the Andronov-Hopf point, (c), (d) $\varepsilon=0.4$, (e) and (f) $\varepsilon=0.9383$ near the turning point.}
    \label{difr}
\end{figure}

\begin{figure}
    \centering
    \subfigure[]{
    \includegraphics[width=0.47 \textwidth]{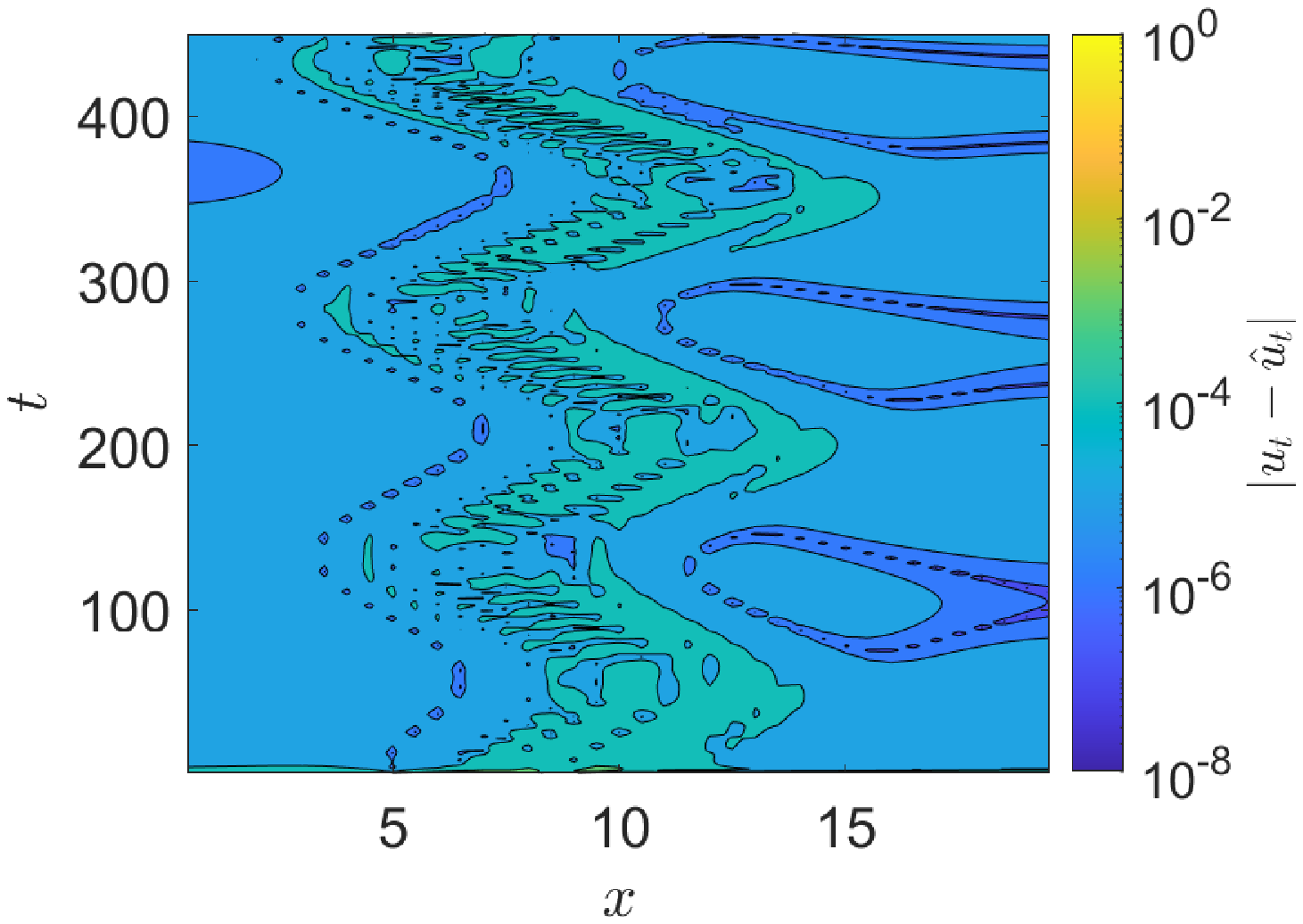}
    }
    \subfigure[]{
    \includegraphics[width=0.47 \textwidth]{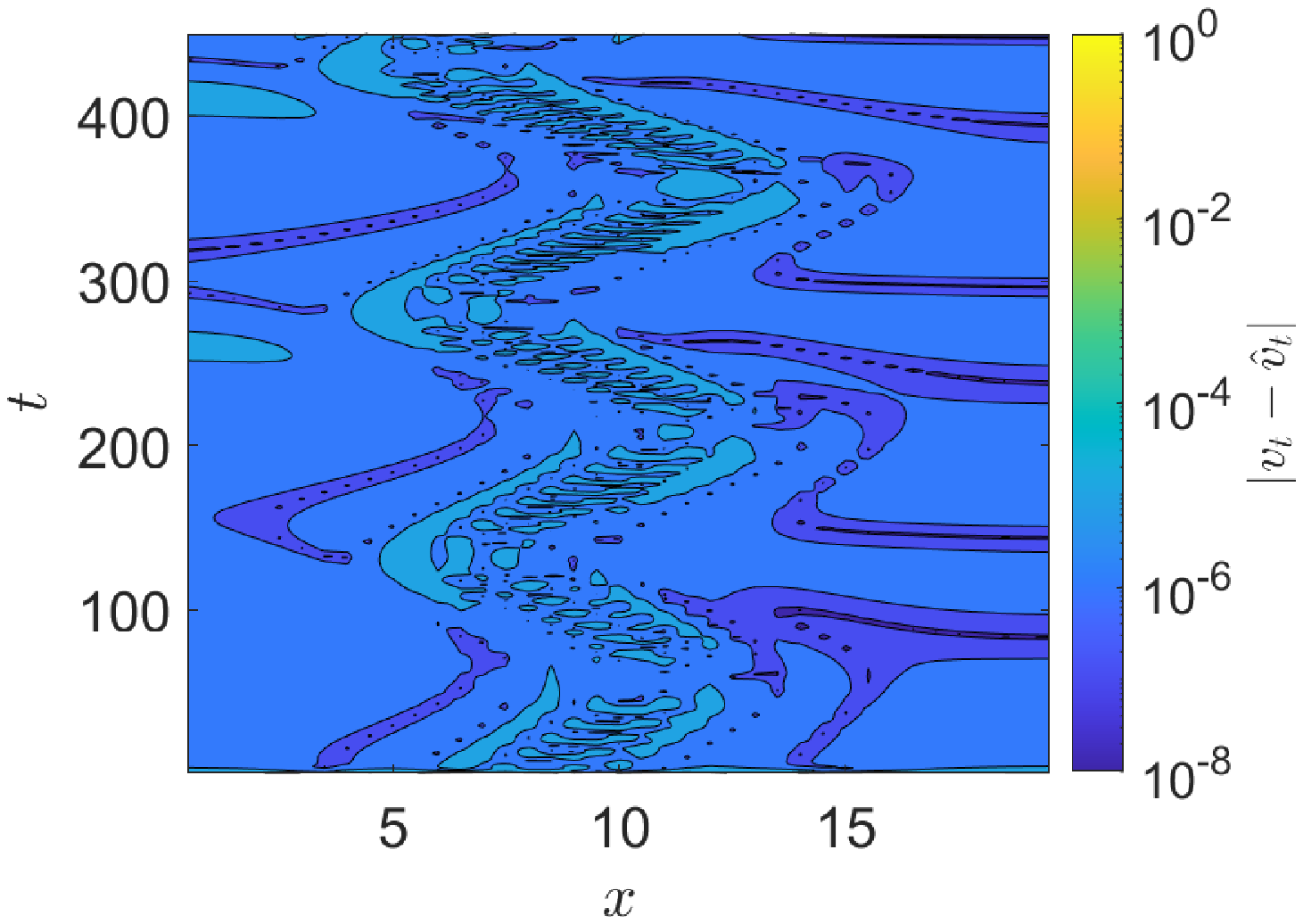}
    }
    \subfigure[]{
    \includegraphics[width=0.47 \textwidth]{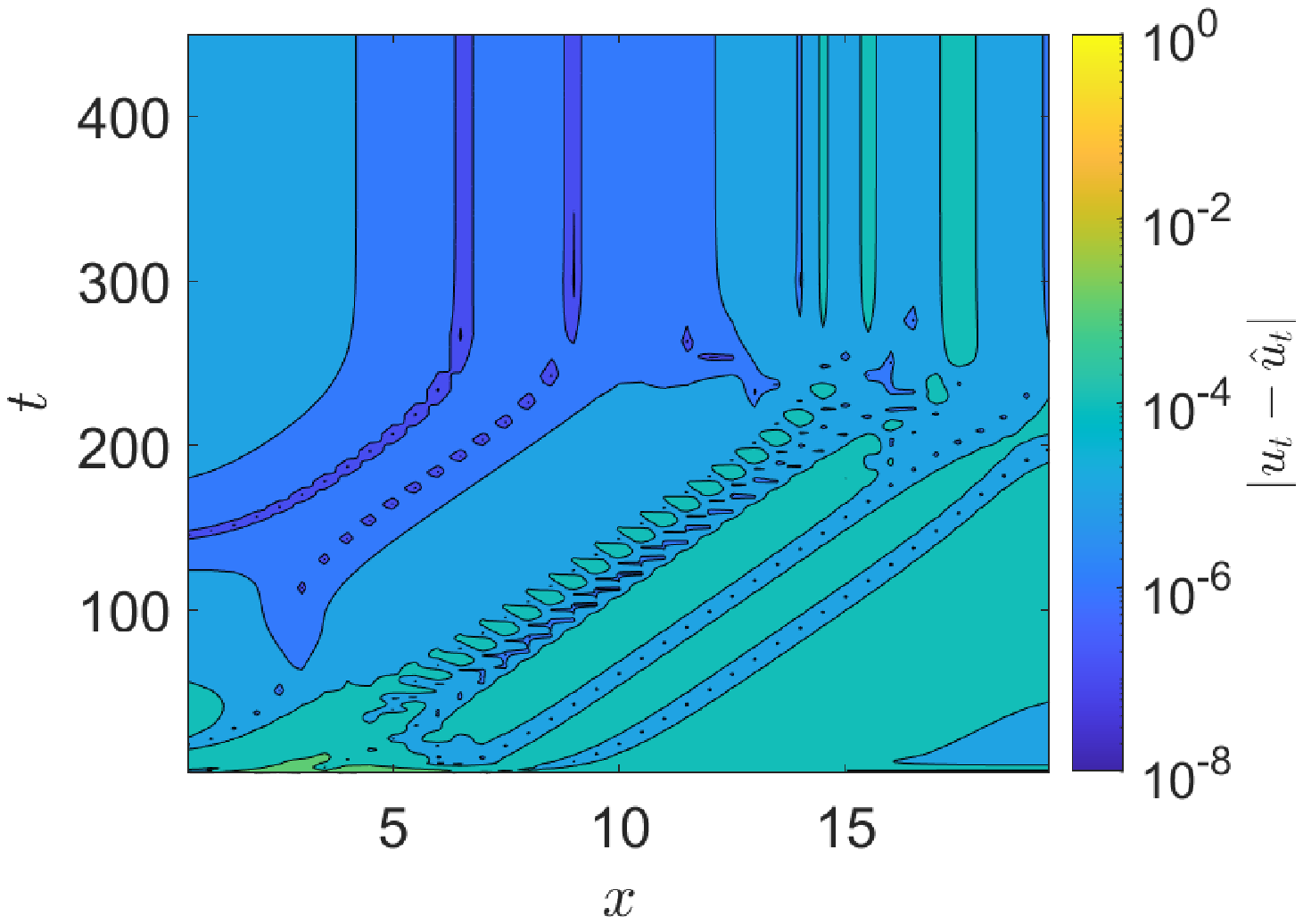}
    }
    \subfigure[]{
    \includegraphics[width=0.47 \textwidth]{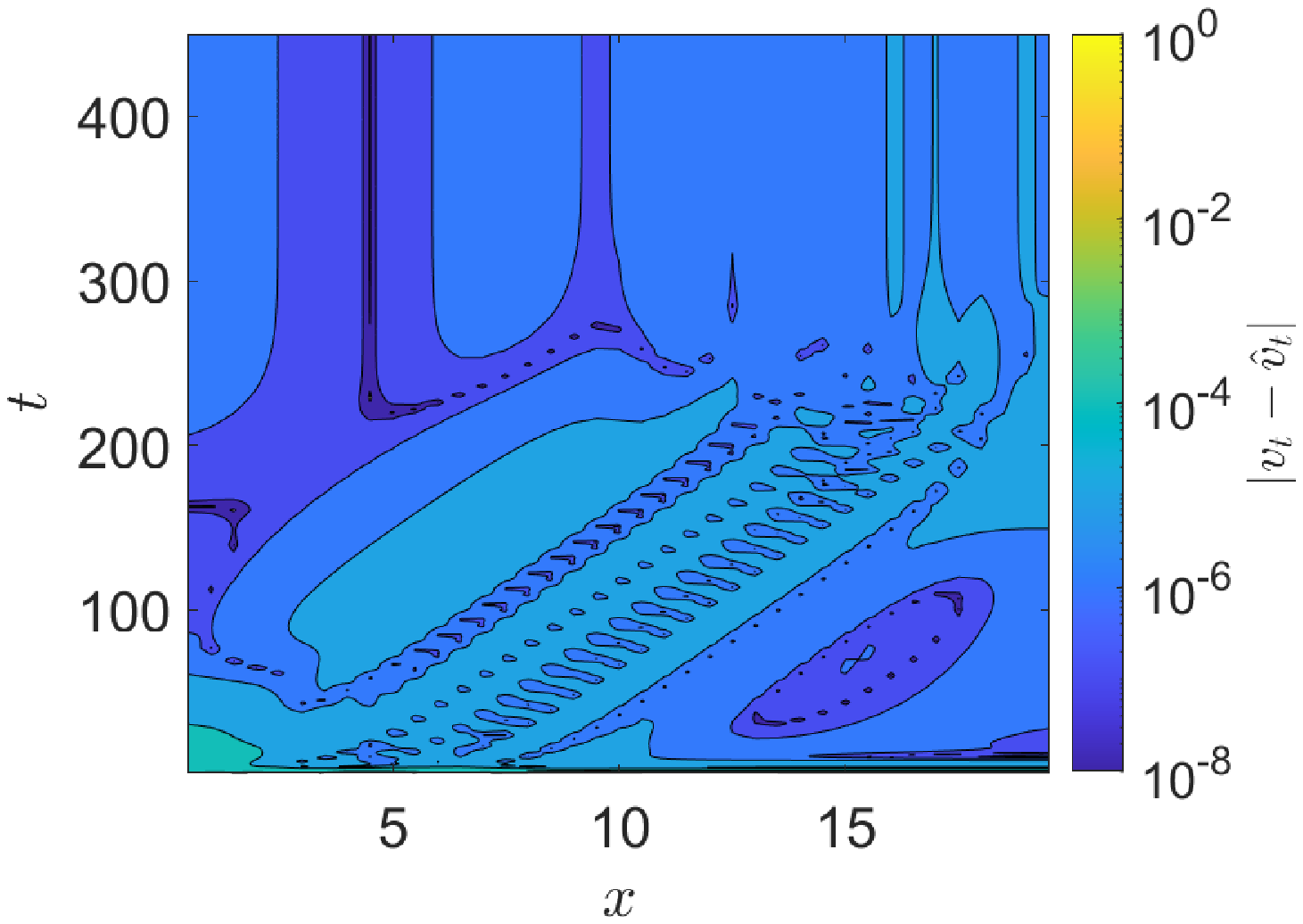}
    }
    \subfigure[]{
    \includegraphics[width=0.47 \textwidth]{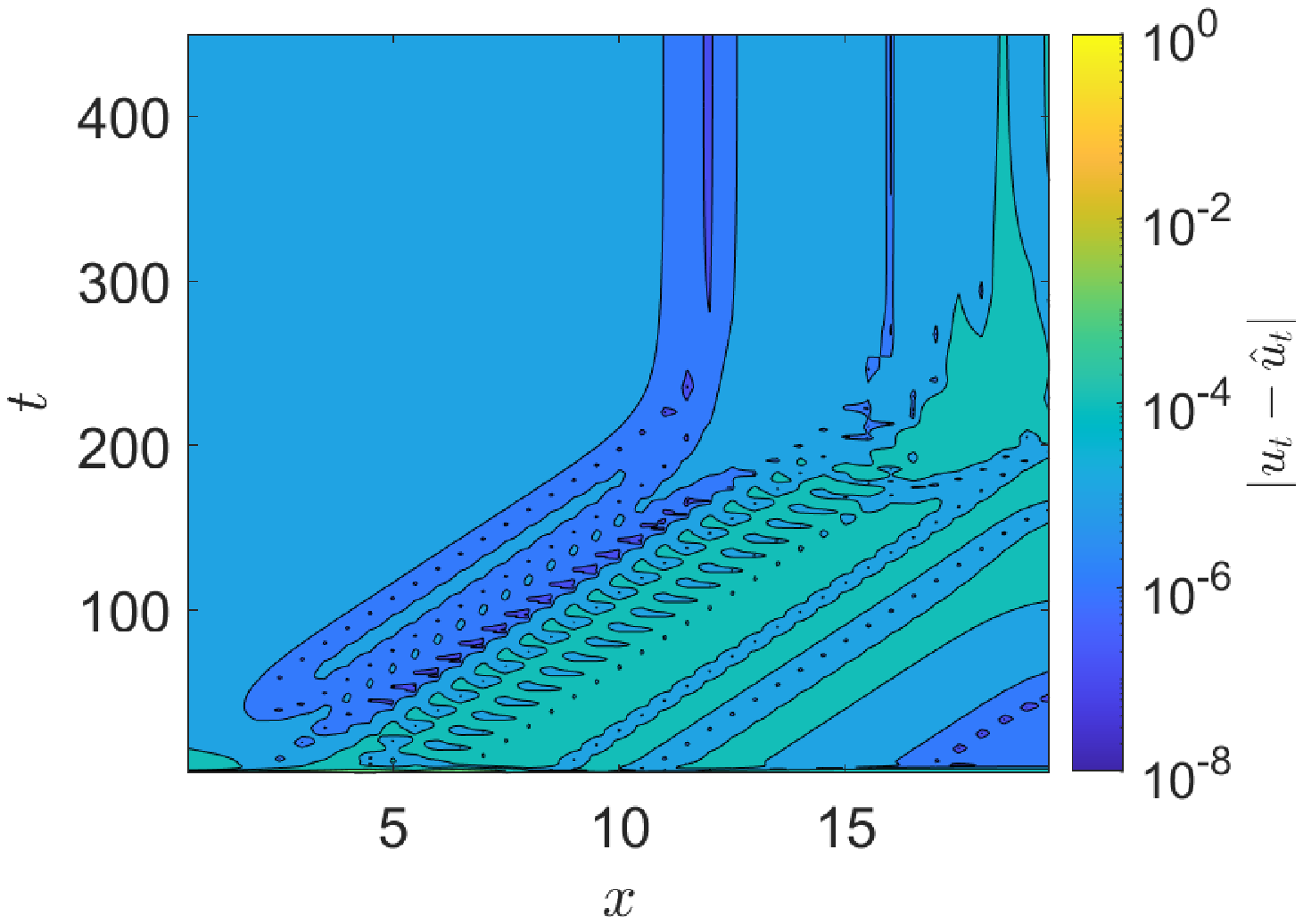}
    }
    \subfigure[]{
    \includegraphics[width=0.47 \textwidth]{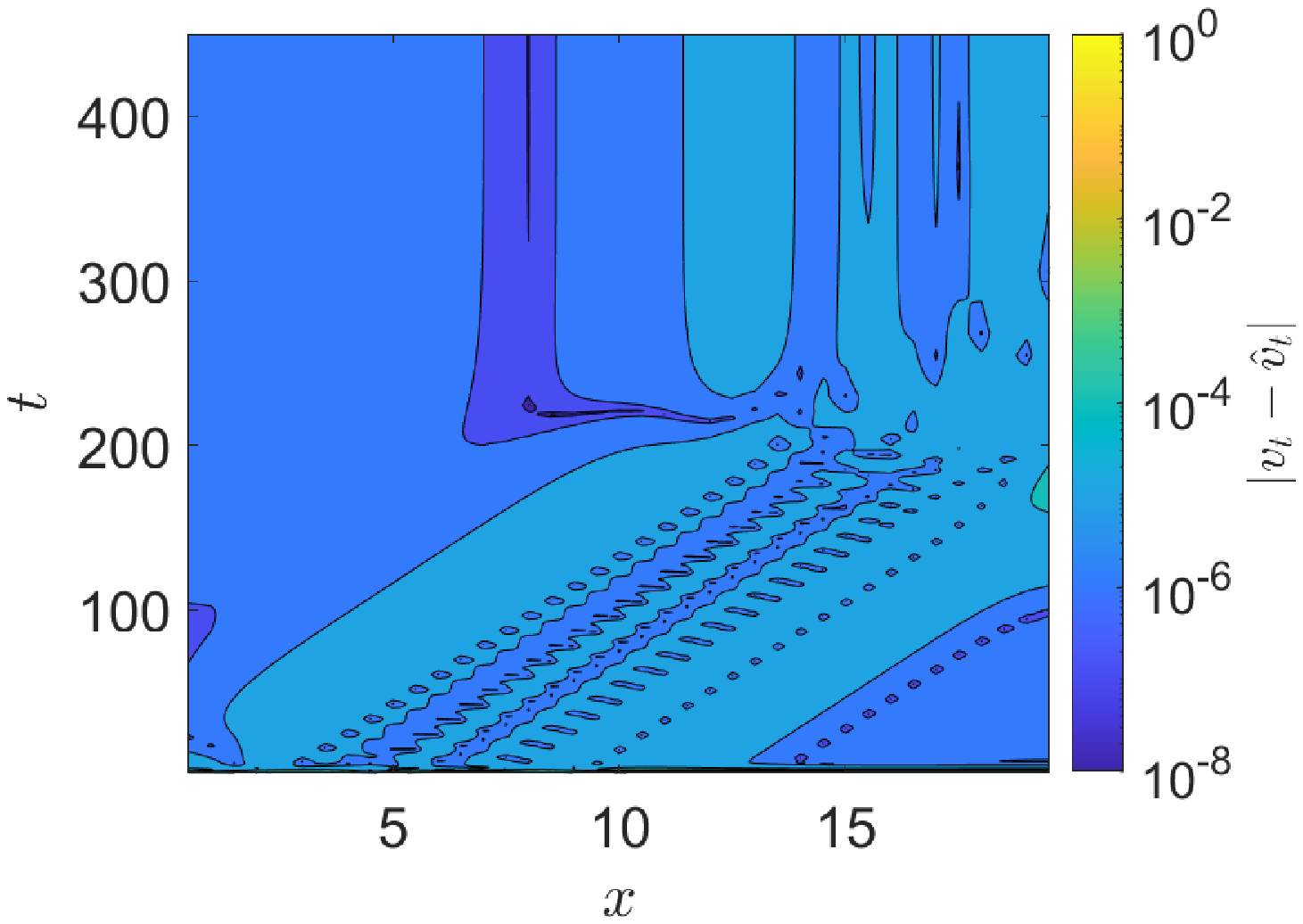}
    }
    \caption{Approximation accuracy in the test data with feature selection as obtained with the RPNNs. Contour plot of the absolute values of differences in space and time, of $|u_t(x,t)-\hat{u}_t(x,t)|$ ((a), (c), and (e)) and of $|v_t(x,t)-\hat{v}_t(x,t)|$ ((b), (d), and (f)) for different values of $\varepsilon$: (a) and (b) $\varepsilon=0.0114$ near the Andronov-Hopf point, (c), (d) $\varepsilon=0.4$, (e) and (f) $\varepsilon=0.9383$ near the turning point.}
    \label{fig:errors_ELM_FS}
\end{figure}

\begin{table}[ht!]
\begin{tabular}{lcccc}
\hline
            & \multicolumn{2}{c}{\textbf{$u_t=(\phi^u_1,\phi^u_2,\phi^u_3)$}} & \multicolumn{2}{c}{\textbf{$v_t=(\phi^v_1,\phi^v_2,\phi^v_3)$}} \\ \hline
            & \textbf{Features}           & \textbf{Total Loss}         & \textbf{Features}           & \textbf{Total Loss}         \\ \hline
\textbf{1d} & $(u)$                       & 4.3E-03                          & $(u)$                       & 7.6E-03                          \\
\textbf{2d} & $(u,v)$                     & 6.37E-06                          & $(u,v)$                     & 1.91E-05                          \\
\textbf{3d} & $(u,v,u_{xx})$              & 2.77E-07                          & $(u,v,v_{xx})$              & 6.29E-07                          \\
\textbf{4d} & $(u,v,u_x,u_{xx})$          & 1.03E-07                          & $(u,v,v_x,v_{xx})$          & 1.34E-07                          \\ \hline
\end{tabular}
\caption{The ``best'' set of variables that effectively parametrize the intrinsic coordinates (($\phi_1^u,\phi_2^2,\phi_3^u$) and ($\phi_1^v,\phi_2^v,\phi_3^v$)) and the corresponding sums of total losses across all the values of the bifurcation parameter $\varepsilon$.}
\label{reglos}
\end{table}
Finally, we repeated the same steps but now using as inputs in the FNNs and RPNNs the reduced input domain as obtained from the feature selection process. 
%Figures \ref{uregr}, \ref{vregr}, \ref{uelmr} and \ref{velmr} show the regression results of the two schemes. 
Table \ref{errors_networks} summarizes the performance of the schemes on the training and the test sets. 
Figures \ref{difr} and \ref{fig:errors_ELM_FS} illustrate the norms of the differences between the predicted from the FNNs and RPNNs  and the actual time derivatives of both variables.

Hence, as it is shown, the proposed feature selection approach based on the parsimonious Diffusion Maps revealed correctly the structure of the embedded PDEs in the form of:
\begin{equation}
\begin{aligned}
\frac{\partial u(x,t)}{\partial t}=\hat{F}^u(u(x,t),v(x,t),u_{xx}(x,t),\varepsilon),\\
\frac{\partial v(x,t)}{\partial t}=\hat{F}^v(u(x,t),v(x,t),v_{xx}(x,t),\varepsilon)
\end{aligned}
\end{equation}
where $\hat{F}^u$ and $\hat{F}^v$ are the outputs of the FNNs (or the RPNNs). 
%Implementing the pseudo-arc-length conntiuation approach using $\varepsilon$; the Jacobian matrices of $F^u$ and $F^v$ were estimated numerically using central finite differences.   
The constructed bifurcation diagram with feature selection is shown in Figure \ref{bifdiag_final}. Using the FNNs, we estimated the Andronov-Hopf point at $\varepsilon\approx 0.0195$ and the turning point at $\varepsilon\approx 0.9762$. Using the RPNNs, we estimated the Andronov-Hopf point at $\varepsilon\approx 0.0192$ and the turning point at $\varepsilon\approx 0.9752$.

\section{Conclusions}
Building on previous efforts \cite{Lee2020}, we present a machine learning methodology for solving the inverse problem in complex systems modelling and analysis, thus identifying effective PDEs from data and constructing the coarse-bifurcation diagrams. The proposed approach is a three tier one. In the first step, we used non-linear manifold-learning and in particular Diffusion Maps to select an appropriate set of coarse-scale observables that define the low-dimensional manifold on which the emergent dynamics evolve in the parameter space. At the second step, we learned the right-hand-side of the effective PDEs with respect to the coarse-scale observables; here we used shallow FNNs with two hidden layers and single layer RPNNs which basis functions were constructed using appropriately designed random sampling. Finally, based on the constructed black-box machine learning models, we constructed  the coarse-grained bifurcation diagrams exploiting the arsenal of numerical bifurcation toolkit. To demonstrate the approach, we used $D1Q3$ Lattice Boltzmann simulations of the FitzHugh-Nagumo PDEs and compared the machine learning constructed bifurcation diagram with the one obtained by discretizing the PDEs with central finite differences.\\
The results show that the proposed machine learning framework was able to identify the ``correct" set of coarse-scale variables that are required to model the emergent dynamics in the form of PDEs and based on them systematically study the coarse-scale dynamics by constructing the emerging macroscopic bifurcation diagram. In terms of approximation accuracy of the macroscopic dynamics, both schemes (the two hidden- layers FNNs and the single-hidden layer RPNNs) performed for all practical purposes equivalently. However, in terms of the computational cost in the training phase, the RPNNs were 20 to 30 times faster than the two hidden layers FNNs. Hence, the proposed RPNN scheme is a promising alternative approach to deep learning for solving the inverse problem for high-dimensional PDEs from big data \cite{karniadakis2021physics,raissi2018deep,raissi2018hidden}.\\
Here, we have focused on the construction of black-box machine learning models of the emergent dynamics of complex systems that can be described by (parabolic) PDEs over the parametric space. The proposed approach can be extended to construct ``gray-box'' models by incorporating information from the physics into the machine learning architecture \cite{lovelett2019partial,karniadakis2021physics}. Furthermore, based on previous efforts aiming at extracting normal forms of ODEs from data \cite{yair2017reconstruction}, the proposed approach can be extended to discover normal forms of high-dimensional PDEs.

\section*{Acknowledgments}
This work was supported by the Italian program ``Fondo Integrativo Speciale per la Ricerca (FISR)'' - FISR2020IP 02893/ B55F20002320001. Y.K. acknowledges partial support from US Department of Energy and the US Air Force Office of Scientific Research.

%Bibliography
%\bibliographystyle{unsrt}  
%\bibliography{references}  

\end{document}